\documentclass[12pt]{report}
\usepackage{german}
\usepackage{makeidx}
\usepackage{amssymb}
\usepackage{stmaryrd}
\usepackage{wasysym}
\usepackage{fontenc}
\makeindex
\author{Andreas Johann Raab\\
Luisenstrasse 60, 80798 M"unchen,\\
Federal Republic of Germany\\
E-mail:andreas@andreasjohannraab.de}
\title{Konzepte der abstrakten Ergodentheorie\\
Erster Teil:\\
Abstraktion und Analyse des Attraktorbegriffes}
\begin{document}
\maketitle
\newpage
{\Large {\bf The Subject}}\newline
\newline
Our conception of a generalized ergodic theory shall exceed the 
generality of general topology:
In this first part of the generalized ergodic theory 
we shall investigate the logical constitution
of the term of the attractor of a given flow.
Our investigations of the logical constitution
of the term of the attractor 
start with an examination of invariant topologies in the first chapter and 
our investigations
lead to the constructive phase, which is the second chapter:
There we find 
a pluralism of
generalized attractors.
However the theorems 3.3. and 3.8 give 
some emphasis to 
certain forms of generalized attractors
within this pluralism: At the beginning of the second chapter 
we construct the set $\mbox{{\bf @}}(\xi,\mathcal{A})$ of attractors according to (2.1) and (2.2)
for a given set $\mathcal{A}$ of subsets of the phase-space $Y=\bigcup \mathcal{A}$ and
for the evolution
\begin{displaymath}
\begin{array}{c}
\xi:Y\times \mathbb{R}\to Y\ ,\\
(x,t)\mapsto \xi(x,t):=\xi^{t}(x)\ ,
\end{array}
\end{displaymath}
which is determined by a given one-dimensional bundle of flows $\{\xi^{t}\}_{t\in\mathbb{R}}$. 
However we do not exhaust the topic only presenting constructions of 
generalized attractors:
We shall show the following very general insensitive ergodic theorem, which combines 
the generalized insensitive ergodic theorem 3.3 and the explication-theorem 3.8:
Let 
$$\mathcal{A}^{c}:=\{Y\setminus X:X\in\mathcal{A}\}$$
be the complementary set of subsets of the phase-space $Y=\bigcup \mathcal{A}=\bigcup\mathcal{A}^{c}$
and let
$$\overline{\mathcal{A}}:=\{\mathbf{cl}_{\mathcal{A}}(Z):Z\subset Y\}\ $$
be the set of closed subsets of the phase-space.
Whereat
for any subset $Z\subset Y$  
$$\mathbf{cl}_{\mathcal{A}}(Z):=\bigcap\{X\in\mathcal{A}^{c}\setminus\{\emptyset\}:X\supset Z\}$$
denotes the generalized closure of $Z\subset Y$.
The implication
\begin{equation}\label{caton}
P\in \mathcal{A}\setminus\{\emptyset\}\Rightarrow \exists\ Q\in \mathcal{A}\setminus\{\emptyset\}:Q\subset(\xi^{t})^{-1}(P)
\end{equation}
modells a generalized form of continuosity
of the flow $\xi^{t}$, which 
is related to the given set of subsets of the phase-space $\mathcal{A}$.
Exactly if the implication (\ref{caton}) is true, we 
call the bijection $\xi^{t}$ Cantor-stetig. In the end of the third chapter we investigate this form 
of generalized continuosity as far as ergodic aspects are concerned: We show, that the explication-theorem 3.8 is true.
If we combine 
the theorem 3.3 and the explication-theorem 3.8 we find the following:
\newline\newline
If 
for any real number $t\in\mathbb{R}$
the implication (\ref{caton})
is true, then the set
$$[[\xi]]_{\mathcal{A}}:=\{\mathbf{cl}_{\mathcal{A}}(\xi(x,\mathbb{R})):x\in Y\}\subset \mbox{{\bf @}}(\xi,\overline{\mathcal{A}})$$
is a partiton of the phase-space $Y$. Note, that $\mbox{{\bf @}}(\xi,\overline{\mathcal{A}})\not=\mbox{{\bf @}}(\xi,\mathcal{A})$.
Also note, that
this assertion is an ergodic theorem indeed.

\newpage
\tableofcontents
\newpage
{\bf {\Large Einleitung}}\newline\newline
{\small Die Ergodentheorie kam nicht in der Mathematik zur Welt.\newline\newline
Denn sie stammt aus einem gleichermassen 
archetypischen wie
fundamentalen Dualismus, 
der ein wenig im Schatten der wohl ber"uhmtesten
Dualismen steht, 
n"amlich im Schatten des neuzeitlichen Leib-Seele-Dualismus
und des modernen Welle-Teilchen-Dualismus, an welchen
beiden wir den Wesenszug des Dualismus
ablesen k"onnen:
Ein Dualismus besteht zwischen zwei zueinander komplement"aren Denkans"atzen.
Deren Komplementarit"at ist es dabei,
dass kein dritter Ansatz denkbar zu sein scheint, wobei
beide Ans"atze unter jeweiligen Gesichtspunkten als plausibel oder als evident erscheinen
und wobei beide Ans"atze
aber dennoch gedanklich unvereinbar sind.
Ein Dualismus "ahnelt also dem Paradoxon, 
weil er mit konkreten Paradoxa verbunden ist, die
ihn illustrieren. Die einzelen
Paradoxa machen die
Momente der Spannung sp"urbar, die jeder 
Dualismus hat.
Diese Spannung treibt denjenigen Menschen
um, dessen sie
sich bem"achtigt.\newline
Jener alte und fundamentale Dualismus, aus dem die 
Ergodentheorie kommt, ist derjenige,
der zwischen dem Ansatz des Determinismus und dem Ansatz
der Zuf"alligkeit zun"achst jenseits der Mathematik besteht und der
ohne die physische Zeitlichkeit nicht denkbar ist.  
De Laplace d"urfte die Spannung des 
deterministisch-indeterministischen Dualismus
nicht unbekannt gewesen sein,
hat de Laplace doch Schwerpunkte seiner Forschung sowohl dem Paradigma 
des Determinismus, der Mechanik des Himmels,
als auch dem Zufall gewidmet. 
W"ahrend der Leib-Seele-Dualismus, so wie der 
Welle-Teilchen-Dualismus, keine Asymmetrie 
zeigt, ist es die Version des asymmetrischen deterministisch-indeterministischen Dualismus
die bevorzugt angenommen wird und die in dem Sinn asymmetrisch ist, dass 
beim asymmetrischen deterministisch-indeterministischen Dualismus\index{deterministisch-indeterministischer Dualismus}
der Zufall als ein perspektivischer 
Effekt des verdeckten Blickes auf den 
Determinismus verstanden wird: Die
Zuf"alligkeit gilt hierbei als ein Ph"anomen des
Determinismus.\newline
W"ahrend der Zufall als Ph"anomen galt, der Determinismus hingegen zur allgemeinen Zufriedenheit
in gedanklich nur augenscheinlich vollst"andig 
durchdrungenen Modellen
angenommen war und der Determinismus
sowohl gedanklich als auch
realiter als gut greifbar und als begriffen galt: Da
war es Ludwig Boltzmann, der 
der sich an die Arbeit machte, die Ph"anomenalit"at 
des Zufalles im Detail zu rekonstruieren.
\newline
Boltzmann war es, der schliesslich die Ergodenhypothese formulierte,
und der den Begriff der Ergodik pr"agte, der
sich auf die Physik bezog. Daher kann der Ergodenhypothese nur dann ihre
topologische Unm"oglichkeit vorgehalten werden,
wenn 
davon ausgegangen wird, dass der physikalische Raum
durch eine nat"urliche Topologie strukturiert ist.
Diese Annahme 
der durch 
eine nat"urliche Topologie gegebenen Struktur des physikalischen Raumes 
stand auch f"ur 
E. u. T. Ehrenfest nicht etwa nur ausser Frage; die
Fraglichkeit der Annahme 
der durch 
eine nat"urliche Topologie gegebenen Struktur des physikalischen Raumes 
war als solche vielmehr unerh"ort.\footnote{Auszusprechen, dass f"ur den physikalischen 
Raum eine andere Topologisierung angemessener sein k"onnte als die nat"urliche
Topologisierung durch die euklidische Metrik, dies darf erst gewagt werden, nachdem 
sich etabliert, was Einstein 1916 pr"asentierte; n"amlich, dass
der physikalische Raum eine physikalische Entit"at ist und, dass dessen Kr"ummung seine 
Wirklichkeit in Form wirksamer Gravitation beschreibt. Und selbst nach Einstein
bedarf es noch des quantentheoretischen Hintergrundes, vor dem der Gedanke
aufkommt,
dass der Raum quantisiert sein k"onnte, ehe
die Angemessenheit der nat"urlichen Topologie f"ur den physikalischen Raum in Frage steht.}
Und daher wiesen E. u. T. Ehrenfest in der Arbeit \cite{ehre} darauf 
hin, dass die Quasiergodizit"at die ph"anomenelle Rekonstruktion des
Zufalles in Aussicht stellt, wobei 
Quasiergodizit"at nicht ausdr"ucklich
f"ur die stillschweigend unterstellte nat"urliche Topologie 
des Zustandsraumes formuliert war:
Die Ehrenfestsche Idee der {\em Quasiergodizit"at} formuliert den im wesentlichen vollst"andigen Keim der
mathematischen
Ergodentheorie, 
die wir hier weit fassen und
zu der wir auch die Theorie der Sensitivit"at 
z"ahlen, welche wir im zweiten Teil der abstrakten Ergodentheorie hoffentlich pr"asentieren d"urfen.
Wir fassen die Ergodentheorie weit und
wir schreiben die abstrakte Ergodentheorie ganz in die Mathematik:
Die abstrakte Ergodentheorie ist die Theorie, die das Verh"altnis 
von Phasenflussgruppen $(\Xi,\circ)$
zu jeweiligen Zustandsraumstrukturierungen $\mathcal{A}\subset 2^{\mathbf{P}_{1}\xi}$
untersucht. Letztere sind "uberdeckende Mengensysteme der Potenzmenge $2^{\mathbf{P}_{1}\xi}$
des jeweiligen Zustandsraums $$\bigcup\mathcal{A}=\mathbf{P}_{1}\xi=\mathbf{P}_{2}\xi$$ einer  
Phasenflussgruppe $(\Xi,\circ)$, der die gemeinsame
Definitions- und Wertemenge der Fl"usse $\xi\in \Xi$ ist.
Die abstrakte Ergodentheorie ist also durch eine Objektklasse 
des Cantorschen Universums\index{Cantorsches Universum} gegeben: Die 
abstrakte Ergodentheorie
ist die Beschreibung der Klasse der beschriebenen Paare 
$$\Bigl((\Xi,\circ),\mathcal{A}\Bigr)\ ,$$
exakt welche wir als die Cantor-Systeme\index{Cantor-System}
bezeichnen. Die Klasse der Phasenflussgruppen objektiviert die 
dynamischen Systeme. Die Klasse der Cantor-Systeme umfasst die Klasse
der 
nach \cite{raab} exakt in dem Fall, dass $\mathcal{A}$ eine Topologie ist,
als Smale-Systeme bezeichneten Paare
$((\Xi,\circ),\mathcal{A})$.\index{Smale-System}
Damit sprengt die abstrakte Ergodentheorie die allgemeine Topologie.
Dass die Allgemeinheit dieser Konzeption der 
abstrakten Ergodentheorie
angemessen ist, f"uhrt diese Abhandlung und ein geplantes weiteres 
Traktat vor, das die allgemeine Theorie der Sensitivit"at thematisiert:
Wir finden n"amlich, dass
von der Gebenheit der
Vereinigungs- und der endlichen Schnittabgeschlossenheit der 
Topologien losgel"ost f"ur
die jeweiligen Zustandsraumstrukturierungen $\mathcal{A}$
beliebiger 
Cantor-Systeme $((\Xi,\circ),\mathcal{A})$ ergodentheoretische Sachverhalte bestehen. Dass dieselben
ergodentheoretische Sachverhalte sind,
erkennen wir
vor dem Hintergrund des skizzierten genealogischen Dualismus der Ergodenlehre daran, dass wir
jene Sachverhalte spezifizieren k"onnen, zu Aussagen,
die ergodentheorischen Bedarf beispielsweise der Physik decken.\newline
Die Klasse der Cantor-Systeme objektiviert ein Konzept.\index{Konzept}
Denn den Begriff des Konzeptes objektivieren
wir dadurch, dass wir unter einem Konzept
eine Klasse verstehen, welche eine Kokatenation jeweiliger Klassen
ist. Demnach ist beispielsweise die Klasse der Massr"aume
ebenso ein Konzept, wie es die Klasse der topologischen R"aume ist. Da auch die Klasse
der Gruppen ein Konzept ist, wohingegen die Klasse der blossen Mengensysteme 
kein Konzept ist, ist die Klasse der Cantor-Systeme ein Konzept, welches
das Konzept der Gruppen erweitert. Der Philosoph mag an dieser 
Verfassung des Begriffes des Konzeptes vermissen, dass in
dessen Objektivierung die jeweilige Intention nicht ausgesprochen wird,
welche mit den Komponenten eines jeweiligen Konzeptes gemeint ist.
Intention ist aber vielleicht kein Gegenstand, sondern 
Wollen, das eine Repr"asentation als neuronale Aktivit"at haben mag. 
\footnote{
Neben der Schwierigkeit
hochgradig hypothetischen Vermutens dabei,
das
Pendant des jeweiligen 
Wollens zu treffen,
ist die Zuordung einer jeweiligen neuroaktiven Repr"asentation zu 
dem Wollen, dem sie jeweils Nachweis sein soll, 
mit folgender noch gravierenderen
Schwierigkeit verbunden: Mysterium 
genug ist es, 
wenn das danach gefragte Kind, ob es 
ein Bonbon wolle,
ermisst, dass dem so sei.
Was die bewusste Identifizierung eigenen Wollens
jenseits vorgelegter Bonbons betrifft, tut
sich ein ge"ubtes Bewusstsein nicht schwer, 
wenn es sich selbst --  erfindet. Wollen ist, salopp gesagt, ziemlich 
transzendent und vielleicht nur 
durchaus unter Sozialisierungsdr"ucken gut einge"ubtes
W"ahnen, welches das Gemeinwesen synchronisiert. 
Entdecken wir die auf dem Planeten heimische
Religion des Zweckes und -- Sinns --  und Gemeinsinns?
So manchen mannhaft hochaufgewachsenen Studenten oder Absolventen 
der Ingenieurswissenschaft
m"ogen diese Worte sehr
verst"oren, nicht aber jeden. Von der Beschw"orung des \glqq Konzeptes\grqq in Zeiten der 
\glqq Krise\grqq wendet sich der Denker ab, weil er seine Kritik den um Steuerung 
bem"uhten
Menschen in der Not lieber
erspart.}
Die Absicht und den Zweck mathematischer Konstruktionen gibt es 
freilich -- jenseits der Mathematik. 
Absichten liegen im
Irgendwo jenseits der Objektivierbarkeit. Absichtlichkeit ist uns 
allen vertraut als das Lebenszeichen der Geistigkeit.
\newline 
Wo Absicht und Zweck mathematischer Konstruktionen  
jenseits der Objektivierbarkeit liegen,  
wird der Zweck mathematischer Konstruktionen in den jeweiligen
Anwendungen erlebbar.
Welche Nachfrage unsere Abstraktionen wecken k"onnen, dies
muss allerdings erst noch eigens illustriert werden.
Zun"achst aber wollen wir vor allem konzipieren, was nicht ohne K"uhnheit ist.
Die Spezifizierungen der abstrakten Ergodentheorie und die ihnen entsprechenden  
Interpretationen jeweils durchzuf"uhren, ist die Anwendung der
abstrakten Ergodentheorie, die wir in den ersten beiden Teilen der
abstrakten Ergodentheorie
noch nicht vorstellen. 
\newline
Wir schreiben die abstrakte Ergodentheorie ganz in die Mathematik. Dass
die Ergodentheorie nicht in der Mathematik zur Welt kam, wollen wir dabei aber nicht vergessen.}
\newline
\newline
Was ist ein Attraktor? 
Attraktoren werden "ublicherweise als Kompakta aufgefasst.
Und zwar sind es diejenigen flussinvarianten Kompakta $\chi\subset Z$ bez"uglich der Topologie $\mathbf{T}(Z)$
des Zustandsraumes $Z$, auf dem
der Phasenfluss $\{\Psi^{t}\}_{t\in\mathbb{R}}$ 
gegeben ist,
f"ur den f"ur alle $t\in\mathbb{R}$ 
\begin{equation}\label{veorm}
\Psi^{t}(\chi)=\chi
\end{equation} 
gilt, welche die Attraktoren des jeweiligen 
Phasenflusses $\{\Psi^{t}\}_{t\in\mathbb{R}}$ bez"uglich der Topologie $\mathbf{T}(Z)$
sind -- falls
diese flussinvarianten Kompakta $\chi$ "uberdies das Kriterium\index{Koh\"arenzkriterium} erf"ullen,
dass f"ur alle bez"uglich der Relativtopologie $\mathbf{T}(Z)\cap\chi$ offenen und nicht leeren Teilmengen $a,b\subset\chi$ eine Zahl $t\in \mathbb{R}$
existiert, f"ur welche die Reichhaltigkeitsaussage
\begin{equation}\label{verm}
\Psi^{t}(a)\cap b\not=\emptyset
\end{equation} 
gilt. 
Dies ist der weit verbreitete Begriff des Attraktors.
\newline\newline
{\small Was hierbei in Gestalt der    
Reichhaltigkeitsaussage (\ref{verm}) formuliert ist, fasst die Koh"arenz innerhalb eines Attraktors.
Zusammen mit der Forderung der Kompaktheit ist die 
in diversen Abwandlungen formulierbare
Koh"arenzbedingung (\ref{verm}) an eine
flussinvariante Teilmengen 
des Zustandsraumes, die als ein Attraktor des Phasenflusses $\{\Psi^{t}\}_{t\in\mathbb{R}}$ gelten soll,
f"urwahr ein erstaunlicher Handgriff,\index{entscheidender Handgriff} ein auff"alliges Glanzlicht der an Glanzlichtern nicht
armen 
Geschichte mathematischer
Begriffsbildungen!
Das erstaunlichste an der Begriffsbildung des Attraktors 
ist nicht die 
Flexibilit"at der Koh"arenzbedingung an einen
Attraktor, die
in verschiedenen Modifizierungen formulierbar ist, auch wenn
alleine schon diese auf der Ebene der Begriffsverfassung gegebene 
Metainvarianz des Begriffes des Attraktors gegen"uber
Ab"anderung der
Koh"arenzbedingung schon sehr erstaunlich
ist.\newline
Die Begriffsbildung des Attraktors will von vorne herein
auf eine Elementarisierung der 
flussinvarianten Mengen hinaus:
Jede Vereinigung einer Menge flussinvarianter Mengen des Zustandsraumes 
ist offensichtlich ebenfalls flussinvariant. Die Attraktoren
sollen aber nicht irgendwelche flussinvarianten Mengen sein. Sie sollen
in einem zun"achst noch unbestimmten Sinn minimal sein.
Die n"achstliegende Elementarisierung der 
flussinvarianten Mengen
ist es,
einfach die kleinsten flussinvarianten Mengen
zu nehmen, welche die sehr wohl 
pr"azise fassbaren integeren Mengen jeweiliger 
flussinvarianter Topologien sind, wie wir im Abschnitt 1.1 darlegen. Indess:
\newline 
\newline
{\em Der Begriff des Attraktors transzendiert den Begriff der Integrit"at. }\footnote{Diese
Formulierung muss dunkel erscheinen, weil wir zum einen den Begriff der Integrit"at noch 
nicht definierten, immerhin jedoch bereits sagten, dass die kleinsten flussinvarianten Mengen
integere Mengen seien; was hier noch mehr irritiert, mag sein, dass
von Transzendenz die Rede ist. Es besteht hier aber folgende 
konkrete
Analogie zu dem 
mit dem Begriff der Transzendenz in der Mathematik konotierten 
Sachverhaltsklassiker:
Die abgeschlossenen H"ullen der Folgen rationaler Zahlen geben uns reelle Zahlen und
die abgeschlossenen H"ullen der integeren Mengen
flussinvarianter Topologien sind nach dem Satz von der Existenz der Zimmer
des Traktates "uber den elementaren Quasiergodensatz \cite{raab} 
gerade Attraktoren.}
\newline
\newline
Das erstaunlichste an der Begriffsbildung des Attraktors ist,
dass sie auf eine Elementarisierung der 
flussinvarianten Mengen hinaus will,
jedoch
an der n"achstliegenden Elementarisierung der 
flussinvarianten Mengen vorbei geht und
en passant -- ins Schwarze trifft.
\newline
Wohingegen wir hier vorerst wohl noch ins Obskure
verbannt sind. Stellen wir uns naiv und fragen so:
Inwieweit erf"ullt 
diese Vereinbarung,
von gerade denjenigen flussinvarianten Mengen $\chi$, die
als Attraktoren des Phasenflusses $\{\Psi^{t}\}_{t\in\mathbb{R}}$ gelten sollen, 
die gem"ass (\ref{verm}) verfasste Koh"arenz 
zu verlangen,
folgenden Zweck: Den Zweck, diejenigen flussinvarianten Mengen der invarianten Topologie der jeweiligen Flussfunktion $\Psi$
auszuzeichen, die gerade die \glqq kleinsten Mengen\grqq der invarianten Topologie sind? 
So k"onnen wir fragen, weil das Basiskriterium bei der Konstituierung des Begriffes des 
Attraktors das Kriterium der Flussinvarianz (\ref{veorm}) ist, f"ur das die
Transitivit"atsaussage in Form der Implikation (\ref{verimp}) gilt.
Es steht uns dabei frei, in unsere Zukunft zu investieren und uns dumm stellen:\footnote{Bei der Jagd auf Ignoranz ist der Grundschritt der selbstreferente, sich 
dumm zu -- stellen. Siehe E.K"astners Feuerzangenbowle.}\newline
Die erhobene Frage, inwieweit 
die \glqq kleinsten Mengen\grqq der invarianten Topologie Attraktoren sind,
umfasst auch die Frage, ob es sein kann, dass 
die \glqq kleinsten Mengen\grqq der invarianten Topologie dem Koh"arenzkriterium nicht
gen"ugen.  
Wir bezeichnen dieses Kriterium der Flussinvarianz (\ref{veorm}) deshalb als Basiskriterium,\index{Basiskriterium} weil die 
Konstituierung des Begriffes des 
Attraktors mit ihm anf"angt.
Das Basiskriterium stellt insofern eine Bedingung, die sich gleichsam an den Rohstoff
richtet, in dem Attraktoren zu finden sind und es ist insofern das Kriterium
bei der Auswahl des Rohstoffes, bei der folgende Implikation
f"ur je zwei Teilmengen $\chi_{1}\subset Z$ und $\chi_{2}\subset Z$ des Zustandsraumes $Z$ gilt:
\begin{equation}\label{verimp}
\begin{array}{c}
\Bigl(t\in\mathbb{R}\Rightarrow\Psi^{t}\chi_{1}=\chi_{1}\ \land\ \Psi^{t}\chi_{2}=\chi_{2}\Bigr)\\
\Rightarrow\\
\Bigl(t\in\mathbb{R}\Rightarrow\Psi^{t}(\chi_{1}\cup\chi_{2})=\chi_{1}\cup\chi_{2}\Bigr)\ ,
\end{array}
\end{equation} 
gleich, welche Topologisierung $\mathbf{T}(Z)$ des Zustandsraumes $Z$
vorliegt. Und ebenfalls unabh"angig von der Topologisierung $\mathbf{T}(Z)$ desselben
gilt hierbei: Die Vereinigung speziell 
zweier Kompakta bez"uglich der Topologie $\mathbf{T}(Z)$
ist allemal ein Kompaktum bez"uglich derselben.
Ist
der Phasenfluss $\{\Psi^{t}\}_{t\in\mathbb{R}}$ stetig, d.h ein Phasenfluss, f"ur den f"ur alle $t\in\mathbb{R}$ die
Bijektionen $\Psi^{t}$ des Zustandsraumes $Z$ auf denselben
stetig bez"uglich der Topologie $\mathbf{T}(Z)$ sind, so gilt f"ur alle $t\in\mathbb{R}$ 
und f"ur alle Teilmengen $\chi\subset Z$ des Zustandsraumes $Z$
die Implikation
\begin{equation}\label{veriimp}
\Psi^{t}\chi=\chi\Rightarrow\\
\Psi^{t}(\mathbf{cl}_{\mathbf{T}(Z)}(\chi))=\mathbf{cl}_{\mathbf{T}(Z)}(\chi)\ ,
\end{equation} 
wobei
$\mathbf{cl}_{\mathbf{T}(Z)}(\chi)$ den Abschluss der Menge $\chi\subset Z$
bez"uglich der Topologie $\mathbf{T}(Z)$ bezeichne. 
Falls $(Z,\mathbf{T}(Z))$ ein kompakter topologischer Raum ist, ist $\mathbf{cl}_{\mathbf{T}(Z)}(\chi)$ 
kompakt.\newline
Wir k"onnen daher die Frage stellen, ob das an die 
flussinvarianten Mengen gestellte Koh"arenzkriterium (\ref{verm}) dem Elementarisierungszweck
dienen soll. Soll
es bezwecken, die kleinsten flussinvarianten Mengen -- oder auf eine
raffiniertere Weise die in einem noch unbestimmten Sinn kleinsten flussinvarianten Mengen 
auszuzeichen? Wobei jene in einem noch unbestimmten Sinn kleinsten flussinvarianten Mengen dann, falls sie 
bez"uglich der Topologie $\mathbf{T}(Z)$
kompakt sind, als Attraktoren gelten d"urfen?
Diese 
Fragestellung sei
durch folgende, verdeutlichende Karrikatur der Naseweisheit
dramatisiert:\newline
Die Naseweisheit nimmt an, die V"ater des Begriffes des Attraktors h"atten, wohl wissend, dass
jede Vereinigung  
einer Menge flussinvarianter Mengen ebenfalls flussinvariant ist,  
beabsichtigt,
besondere, elementare, flussinvariante Mengen auszuzeichen.
Zun"achst verlegen um eine Objektivierung der Elementarit"at der 
gemeinten elementaren flussinvarianten Mengen,
h"atten sie aber, wohlerfahren 
insbesondere
mit dem kontinuierlichen Phasenfluss endlichdimensionaler reeller oder 
komplexer Zustandsr"aume, auch um folgenden Sachverhalt gewusst:
Diejenigen 
flussinvarianten Mengen $\chi$ endlichdimensionaler reeller oder 
komplexer Zustandsr"aume erweisen sich als mit den gemeinten elementaren
flussinvarianten Mengen auf zufriedenstellende Weise "ubereinstimmend,
welche ein gewisses Koh"arenzkriterium
erf"ullen; n"amlich das Koh"arenzkriterium, dass f"ur alle 
in der jeweiligen flussinvarianten Menge $\chi$ liegenden und
bez"uglich der nat"urlichen Topologie 
jeweiliger endlichdimensionaler reeller oder 
komplexer Zustandsr"aume offenen
Mengen $a$ und $b$ die Reichhaltigkeitsaussage
(\ref{verm}) gilt.
\newline
Was wir offenbar nicht im Schilde f"uhren, ist, den Ruhm der V"ater des Begriffes des Attraktors
schm"alern zu wollen,\footnote{Manche derer leben ja noch. Und der Verschiedenen S"ohne
w"aren da auch noch...}  
sondern stattdessen f"uhren wir im Schilde, die 
Naseweisheit vorzuf"uhren:
Die Naseweisheit meint n"amlich nun, sie spielte den Sch"opfern des Attraktorbegriffes einen Streich, wenn
sie zeigte, dass sich, was ein Attraktor ist, doch auch ganz simpel so sagen 
liesse: \glqq Ein Attraktor ist eine kleinste flussinvariante Menge.\grqq
\newline
Nun darf gelacht werden:
Denn bekanntlich sind die kleinsten flussinvarianten Mengen
kontinuierlicher dynamischer Systeme endlichdimensionaler reeller oder 
komplexer Zustandsr"aume deren Trajektorien oder Fixpunkte und 
dieselben sind gerade lediglich die trivialen Attraktoren.
\newline 
Diese Dramatisierung mit Hilfe der Figur der Naseweisheit macht uns deutlich, dass
die gestellte naive Frage, ob die 
Koh"arenzbedingung an einen Attraktor dem Zweck diene,
die kleinsten flussinvarianten Menge 
zu adressieren, mit Nachdruck verneint 
werden muss.
Den durch (\ref{veorm}) und (\ref{verm}) festgelegten Begriff des Attraktors
bezeichnen wir 
in dem Fall als den Urbegriff des Attraktors,
dass der Zustandsraum eine Teilmenge eines endlichdimensionalen reellen oder 
komplexen Raumes ist, den dessen jeweilige nat"urliche Topologie
strukturiert.\index{Urbegriff des Attraktors}
Dann l"asst sich die Freiheiten offenlassende Frage stellen, unter welchen Umst"anden 
den Urbegriff des Attraktors ausdehnende, denselben gem"ass dem Permanenzprinzip
verallgemeinernde 
Versionen,
so beschaffen sind, dass sie "uber den engen Horizont der trivialen Attraktoren, der
kleinsten flussinvarianten Mengen hinausblicken.
\newline
Die dieser Freiheiten offenlassenden Frage entsprechende Aufgabe ist die,
den besagten Urbegriff des Attraktors gem"ass dem Permanenzprinzip zu erweitern,
ohne dass wir dabei bloss triviale Attraktoren in verallgemeinerter Version verfassen.
Diese Aufgabe
gehen wir an, indem
wir im n"achsten Abschnitt begrifflich pr"azisieren, dass die kleinsten flussinvarianten Mengen
integere Mengen invarianter Topologien sind.
Im darauffolgenden Abschnitt 1.2 er"ortern
wir dieselben im Zuge einer knappen Diskussion der integeren Mengen gemeinsamer invarianter
Topologien homogener Autobolismenmengen.
Anschliessend untersuchen wir im zweiten Kapitel, wie konservative, dem Permanenzprinzip gem"asse Generalisierungen
des Urbegriffes des Attraktors von den terminologisch bestimmten   
kleinsten flussinvarianten Mengen abweichen.
\chapter{Invariante Topologien}}
\section{Selbstduale Topologien und deren\\ integere Mengen}
{\small In diesem ersten Kapitel besprechen wir einfache Sachverhalte der allgemeinen Topologie
im Zuschnitt auf die Thematik der abstrakten Ergodentheorie.
Diese einfachen Sachverhalte werden
in der klassischen Lehrbuchliteratur der allgemeinen Topologie,
z.B. in den nichtsdestotrotz empfehlenswerten Lehrb"uchern \cite{alex}, \cite{manh}, \cite{kell} und \cite{ward}
nicht in diesem Zuschnitt auf die abstrakte Ergodentheorie diskutiert.}
\newline
\newline
Es sei die zum jeweiligen Phasenfluss $\{\Psi^{t}\}_{t\in\mathbb{R}}$ geh"orende, zu demselben 
"aquivalente\footnote{Jeder Phasenfluss $\{\Psi^{t}\}_{t\in\mathbb{R}}$ ist zu einer 
Flussfunktion, n"amlich zu
der Flussfunktion (\ref{veoorm}) "aquivalent. Aber nicht jede Flussfunktion $\Phi$ ist phasisch,\index{phasische Flussfunktion} 
d.h von der Art, dass 
$(\{\Phi(\mbox{id},t):t\in\mathbb{R}\},\circ)$ eine Gruppe und
$\{\Phi(\mbox{id},t)\}_{t\in\mathbb{R}}$ ein Phasenfluss ist;
obgleich
jede Flussfunktion den Determinismus modelliert. Diesen Sachverhalt erl"autern wir gegen Ende des Abschnittes 1.2.} Flussfunktion\index{Flussfunktion}
\begin{equation}\label{veoorm}
\begin{array}{c}
\Psi:Z\times\mathbb{R}\to Z\ , \\
(z,t)\mapsto\Psi^{t}(z)\ .
\end{array}
\end{equation} 
Die Menge aller flussinvarianten Mengen, die der Gleichung (\ref{veorm}) gen"ugen,
ist ein  
Mengensystem, 
das offensichtlich sowohl gegen"uber der Vereinigung "uber Teilmengen dieses Mengensystemes als auch
gegen"uber dem Schnitt "uber Teilmengen desselben abgeschlossen ist.
Da f"ur alle $t\in\mathbb{R}$ 
$$\Psi(\emptyset,t)=\emptyset$$ 
ist, ist dieses Mengensystem 
eine Topologie,
die invariante Topologie\index{invariante Topologie einer Flussfunktion}
der jeweiligen Flussfunktion $\Psi$ 
\begin{equation}\label{flinv}
\widehat{\mathbf{T}}(\{\Psi^{t}:t\in \mathbb{R}\})=\Bigl(\bigcup_{t\in \mathbb{R}}\Psi^{t}\Delta\mbox{id}\Bigr)^{-1}(\{\emptyset\})\ .
\end{equation}
Wobei f"ur je zwei Mengen $A$ und $B$
\begin{equation}
A\Delta B:=(A\setminus B)\cup (B\setminus A)
\end{equation}
deren sogenannte symmetrische Differenz\index{symmetrische Differenz zweier Mengen} bezeichnet, was keine allzu un"ubliche Bezeichnungsweise ist;
was eine beispielsweise in dem
verbreiteten 
Lehrbuch Elstrods \cite{elst} gebrauchte Notation ist. Die
Bezeichnung $\Psi^{t}:=\Psi(\mbox{id},t)$ f"ur alle $t\in \mathbb{R}$ und
f"ur den jeweiligen Phasenfluss f"ur $\{\Psi^{t}\}_{t\in\mathbb{R}}$
ist dabei eine fest
etablierte Schreibweise.
Die Bezeichnung der invarianten Topologie einer Bijektion $\xi\in A^{A}\ni\xi^{-1}$ einer Menge $A$
auf sich selbst
\begin{equation}
\mathbf{T}(\xi):=\{a\subset A:\xi(a)=a\}
\end{equation}
$$=(\xi\Delta\mbox{id})^{-1}(\{\emptyset\})$$
ist eine spezielle Festlegung hier, die auch schon in dem Traktat "uber den elementaren 
Quasiergodensatz \cite{raab} Anwendung fand, in dem die
Bijektionen einer Menge auf sich selbst als Autobolismen\index{Autobolismus} bezeichnet werden.\footnote{Die
Bijektionen einer Menge auf sich selbst einfach
Automorphismen zu nennen, 
ginge ja vielleicht noch durch,
und w"are dem einen oder anderen mit der verfestigten Bedeutung des Begriffes des Automorphismus
vereinbar;
diese Benennung
passte dennoch nicht nur insofern nicht zur Logik der Morpheme, als
bei der Bijektion einer Menge auf sich
keinerlei Bezug zu einer 
eventuell gegebenen oder auch nicht vorhandenen
Strukturierung der jeweiligen Menge besteht,
die insofern amorph ist;
jene Benennung k"onnte in dem Fall, dass auf einer jeweiligen Menge eine Struktur vorhanden ist,
Irritationen hervorrufen. Wir k"onnten Autobolismen allerdings auch als Autojektionen 
bezeichnen. Dies f"ugte sich einerseits besser in die Reihe der Benennungen der
Bijektionen, Injektionen und Surjektionen ein; andererseits ist 
-- anders als bei den Worten \glqq Bijektion \grqq, \glqq Injektion \grqq und \glqq Surjektion \grqq
das Wort Autojektion
eine heterogene F"ugung eines griechischen Pr"afix
an einen lateinischen Stamm.}
$\mathbf{T}(\xi)$ ist also die 
invariante Topologie des Autobolismus $\xi$.\index{invariante Topologie eines Autobolismus}
Wir haben in $\mathbf{T}(\mbox{id})$ einen universellen, auf der Klasse der 
Autobolismen definierten
Operator, dessen jeweiliger Wert die jeweilige invariante Topologie seines Argumentes ist.
Den Operator $\mathbf{T}(\mbox{id})$ 
mit Hilfe des Synonyms $(\xi\Delta\mbox{id})^{-1}(\{\emptyset\})$ f"ur $\mathbf{T}(\xi)$ zu verfassen,
verlangt von uns,
den Definitionsbereich der identischen Abbildung mitanzugeben:
Es gibt F"allen, in denen unklar ist, auf welchen Definitionsbereich sich die 
identische Abbildung $\mbox{id}$ bezieht.
Es sei vereinbart, dass in diesen F"allen die Indzierung mit dem jeweiligen Definitionsbereich Klarheit
schaffe: Ist $A$ eine Menge, so sei $\mbox{id}_{A}$ die auf $A$ definierte identische Funktion;
Ist ${\rm A}$ eine Klasse, so sei $\mbox{id}_{{\rm A}}$ die auf ${\rm A}$ definierte identische Abbildung.
Damit ist die Paraphrasierung
\begin{equation}\label{merc}
\mathbf{T}(\mbox{id})=\{a\subset \mathbf{P}_{1}\mbox{id}:\xi(a)=a\}\\
\end{equation}
$$=(\mbox{id}\Delta\mbox{id}_{\mathbf{P}_{1}\mbox{id}})^{-1}(\{\emptyset\})$$
m"oglich. Sie ist aber wohl auch 
f"ur denjenigen 
schwer auszulegen, dem die Schreibweise
$\mathbf{P}_{1}f$ bzw. $\mathbf{P}_{2}f$ f"ur die Definitionsmenge bzw. Wertemenge einer jeweiligen Funktion $f$
vertraut ist, welche wir hier benutzen.\newline
Des weiteren ist
f"ur eine Menge $\Xi\subset A^{A}$ von lauter Autobolismen der Menge $A$ die
gemeinsame invariante Topologie der Menge Autobolismen $\Xi$ der Schnitt
\begin{equation}
\widehat{\mathbf{T}}(\Xi):=\bigcap_{\xi\in \Xi}\mathbf{T}(\xi)
=\Bigl(\bigcup_{\xi\in \Xi}\xi\Delta\mbox{id}\Bigr)^{-1}(\{\emptyset\})\ .
\end{equation}
Im Fall einer einelementigen Menge von Autobolismen ist
$\mathbf{T}(\xi)=\widehat{\mathbf{T}}(\{\xi\})$.
Der Operator
$$\widehat{\mathbf{T}}(\mbox{id})= \bigcap_{\xi\in\mbox{id}}\mathbf{T}(\xi)\ ,$$
ist auf der Klasse der homogenen
Autobolismenmengen $Q$ definiert. 
Als die Klasse homogenener
Autobolismenmengen bezeichnen wir hierbei die Klasse
derjeniger Autobolismenmengen $Q$,
f"ur die es jeweils
eine Menge gibt, welche die gemeinsame Definitions- und Wertemenge
aller Autobolismen der jeweiligen homogenenen Autobolismenmenge $Q$ ist.\index{homogenene Autobolismenmenge}  
Wir verbleiben ganz innerhalb verfestigter 
Notationstraditionen, wenn wir den gem"ass (\ref{merc})
f"ur einzelne Autobolismen definierten Operator $\mathbf{T}(\mbox{id})$
auch auf irgendwelche Autobolismenmengen $P$ anwenden: Es ist
dann 
\begin{equation}\label{ganke}
\mathbf{T}(P)=\{\mathbf{T}(p):p\in P\}\ .
\end{equation}
$\mathbf{T}(P)$ ist insofern 
f"ur beliebige Autobolismenmengen $P$ definiert,
wohingegen $\widehat{\mathbf{T}}(Q)$ nur f"ur homogenene Autobolismenmengen 
definiert ist.
Jede invariante Topologie $\mathbf{T}(\xi)$ eines Autobolismus $\xi\in A^{A}$ hat sowohl die 
Menge $A$ selbst als auch die leere Menge $\emptyset$ als Elemente, was ein Spezialfall des 
offensichtlichen
Sachverhaltes 
ist, dass 
\begin{equation}\label{mera}
\mathbf{T}(\xi)=\mathbf{T}(\xi)^{c}
\end{equation}
ist,
wobei f"ur jede Topologie $\mathbf{T}$ das Mengensystem
\begin{equation}\label{merb}
\mathbf{T}^{c}:=\Bigl\{\bigcup\mathbf{T}\setminus X:X\in\mathbf{T}\Bigl\}
\end{equation}
ein Mengensystem ist, das im allgemeinen keine Topologie ist. F"ur das Mengensystem $\mathbf{T}^{c}$ sind die
Rollen, die Vereinigungen und Schnitte f"ur Topologien spielen, miteinander vertauscht:
Alle Schnitte "uber Teilmengen aus $\mathbf{T}^{c}$ sind in dem Mengensystem $\mathbf{T}^{c}$; es ist jedoch
nur f"ur endliche Teilmengen aus $\mathbf{T}^{c}$ garantiert, dass die Vereinigung 
"uber diese endlichen Teilmengen in
$\mathbf{T}^{c}$ sind. Es liegt daher nahe, das Mengensystem $\mathbf{T}^{c}$ als das zu
$\mathbf{T}$ duale Topologiekomplement\index{duales Topologiekomplement} zu bezeichnen.\newline
Die 
Identit"at (\ref{mera}) jeder invarianten Topologie $\mathbf{T}(\xi)$ eines Autobolismus $\xi$
mit dem zu ihr dualen Topologiekomplement
ist demnach die Selbstdualit"at invarianter Topologien,\index{Selbstdualit\"at invarianter Topologien}
die sich unmittelbar auf die gemeinsamen invarianten Topologien jeweiliger homogener Autobolismenmengen $\Xi$
"ubertr"agt, f"ur die offensichtlich ebenfalls Selbstdualit"at in Form der G"ultigkeit der Gleichung
\begin{equation}
\widehat{\mathbf{T}}(\Xi)=\widehat{\mathbf{T}}(\Xi)^{c}
\end{equation}
vorliegt. Wenn eine Topologie $\mathbf{T}=\mathbf{T}^{c}$
selbstdual ist, so ist sie abgeschlossen sowohl
gegen"uber dem Schnitt als auch gegen"uber der Vereinigung "uber 
beliebige
Teilmengen der selbstdualen Topologie $\mathbf{T}$. Da im Fall der
Selbstdualit"at $\mathbf{T}=\mathbf{T}^{c}$
f"ur alle ${\rm U}_{1},{\rm U}_{2}\in \mathbf{T}$ das Komplement
${\rm U}_{2}^{c}$ eine offene Menge der selbstdualen Topologie $\mathbf{T}$ ist,
ist 
\begin{equation}\label{meerb}
{\rm U}_{1}\setminus{\rm U}_{2}={\rm U}_{1}\cap{\rm U}_{2}^{c}\in \mathbf{T}\ .
\end{equation}
Die invarianten Topologien jeweiliger Autobolismen oder die 
gemeinsamen invarianten Topologien jeweiliger homogener Autobolismenmengen\index{gemeinsame invariante Topologie einer homogenen Autobolismenmenge}
erinnnern uns insofern an die Niveaulinientopologien\index{Niveaulinientopologie einer Funktion} einer Funktion $f$ 
\begin{equation}\label{aulini}
\begin{array}{c}
\mathbf{T}_{N}(f):=\Bigr\{\bigcup X:X\in\{f^{-1}(\{y\}):y\in\mathbf{P}_{2}f\}\cup\{\emptyset\}\Bigl\}\\
=\Bigr\{f^{-1}(\{y\}):y\in\mathbf{P}_{2}f\Bigl\}^{\cup}\cup\{\emptyset\}\\
=\mathbf{T}_{N}(f)^{c}\ ,
\end{array}
\end{equation}
die es als das Mengensystem aller Vereinigungen "uber Wertefasern $f^{-1}(\{y\})$ f"ur ein
jeweiliges Wertemengenelement $y\in\mathbf{P}_{2}f$ 
zu jeder Funktion $f$ gibt: Jede Niveaulinientopologien ist
ebenfalls selbstdual; wobei es sich hier anbietet, die f"ur jedes Mengensystem definierte Schreibweise
\begin{equation}\label{mere}
\mathcal{A}^{\cup}:=\{\bigcup X:X\in\mathcal{A}\}
\end{equation}
einzuf"uhren; und was
sich damit fortsetzt, den entsprechenden Operator $\mbox{id}^{\cup}$ einzuf"uhren, der 
auf der Klasse der Mengensysteme definiert ist und dessen
Wert f"ur ein jeweiliges Mengensystem $\mathcal{A}$ das Mengensystem $\mathcal{A}^{\cup}$ ist: Genau dann, wenn
$\mathcal{B}(\mathbf{T})$ eine Basis einer Topologie $\mathbf{T}$ ist,\index{Basis einer Topologie} ist
\begin{equation}\label{meree}
\mathbf{T}=\mathcal{B}(\mathbf{T})^{\cup}\cup\{\emptyset\} ,
\end{equation}
weil $\mathcal{A}$ in $\mathcal{A}^{\cup}$ liegt. $(\mbox{id}^{\cup})^{-1}(\{\mathbf{T}\setminus\{\emptyset\}\})$ ist daher
das Mengensystem aller Basen der Topologie $\mathbf{T}$.\newline
Dass die Niveaulinientopologien $\mathbf{T}_{N}(f)$ jeweiliger Funktionen $f$ 
durch die Vereinigungen "uber alle Teilmengen
der
Partition\footnote{ 
Die Menge $\mathbf{part}(X)$ aller Partitionen der Menge $X$ ist eine Menge von Mengensystemen, eine 
Menge von Mengen von Mengen.} deren jeweiliger Wertemenge 
in Wertefasern, d.h durch die Vereinigungen "uber alle Teilmengen
der Partition
$$\Bigr\{f^{-1}(\{y\}):y\in\mathbf{P}_{2}f\Bigl\}\in\mathbf{part}(\mathbf{P}_{2}f)$$
hervorgehen, ist hier per constructionem transparent. 
Dieser Sachverhalt stellt uns aber folgende Frage:
Sind die invarianten Topologien jeweiliger Autobolismen
und damit auch die 
gemeinsamen invarianten Topologien jeweiliger homogener Autobolismenmengen
auf analoge Weise 
aufgebaut:
Sind auch die invarianten Topologien als die Bilder jeweiliger Partitionen  
der Abbildung durch den Operator $(\mbox{id}^{\cup}\cup\{\emptyset\})$ auffassbar, n"amlich als die
$(\mbox{id}^{\cup}\cup\{\emptyset\})$-Bilder
einer Partition, welche die Definitionsmenge des jeweiligen 
Autobolismus bzw. die 
gemeinsame Definitionsmenge der Autobolismen
der jeweiligen homogenen Autobolismenmengen in kleinste invariante Mengen einteilt?
\newline
Der naheliegende Gedanke, dass wir 
invariante Topologien als Niveaulinientopologien 
formulieren k"onnten, um die Frage nach der 
Existenz kleinster invarianter Mengen beantworten zu k"onnen,
stellt uns allerdings vor die Schwierigkeit, die
notwendige und hinreichende 
Bedingung daf"ur, dass
eine Menge $a$ eine invariante Menge eines Autobolismus $\xi$ ist, n"amlich
die G"ultigkeit der Gleichung
$$\xi(a)=a$$
in ein notwendiges und hinreichendes Kriterium daf"ur umzuwandeln,
dass $a$ die Wertfaser einer Funktion $f$ ist; zu diesem Zweck w"are eine Funktion $f_{\xi}$ zu $\xi$ zu konstruieren, f"ur die
die "Aquivalenz 
\begin{equation}
(x,y\in a \Leftrightarrow f_{\xi}(x)=f_{\xi}(y))\Leftrightarrow  \xi(a)=a 
\end{equation}
gilt. Vielleicht k"onnen wir aber allgemein zeigen, dass f"ur jede Topologie $\mathbf{T}$ folgende Implikation gilt:
\begin{equation}
\mathbf{T}=\mathbf{T}^{c}\Rightarrow\ 
(\mbox{id}^{\cup})^{-1}(\{\mathbf{T}\setminus\{\emptyset\}\})\cap\mathbf{part}(\bigcup\mathbf{T})\not=\emptyset\ .
\end{equation}
Dass es also zu jeder selbstdualen Topologie $\mathbf{T}=\mathbf{T}^{c}$ eine 
partitive Basis $\mathcal{B}(\mathbf{T})\in\mathbf{part}(\bigcup\mathbf{T})$ gibt.
\newline
\newline
{\bf Satz 1.1: Die Struktur selbstdualer Topologien}\newline
{\em Jede Topologie ist genau dann selbstdual, wenn sie eine Basis hat, die ihre Vereinigung partioniert.}
\newline
\newline
{\bf Beweis:}\newline
Dass jede Topologie, die eine Basis hat, die ihre Vereinigung partioniert, selbstdual ist,
ist evident.
\newline
Zum Beweis der nicht trivialen Implikationsrichtung des Satzes 1.1
benutzen wir die Konstruktion der Auswahl aus einer Menge:
Im Traktat "uber den elementaren Quasiergodensatz \cite{raab}
wurde f"ur jedes Mengensystem $\mathcal{W}$ und f"ur jede Menge $X$ die
Auswahl des Mengensystemes $\mathcal{W}$ durch die Menge $X$ als das Untermengensystem
\begin{equation}\label{wahl} 
\mathcal{W}_{X}:=\{\chi\in \mathcal{W}:\chi\cap X\not=\emptyset\} 
\end{equation} 
eingef"uhrt.\index{Auswahl durch eine Menge} 
Zu jeder 
Topologie $\mathbf{T}$ existiert der Operator 
\begin{equation}\label{wahal}
\bigcap\mathbf{T}_{\mbox{id}}=\bigcap\{{\rm U}\in \mathbf{T}:{\rm U}\cap \mbox{id}\not=\emptyset\}\ ,
\end{equation}
der nicht nur auf der Topologie $\mathbf{T}$ selbst, sondern 
der insofern auf dem topogischen Raum $(\bigcup\mathbf{T},\mathbf{T})$ gegeben ist, als die 
Definitionsmenge des Operators $\bigcap\mathbf{T}_{\mbox{id}}$ die Menge 
$$\mathbf{P}_{1}\bigcap\mathbf{T}_{\mbox{id}}=2^{\bigcup\mathbf{T}}$$
ist.
Wenn die Topologie $\mathbf{T}=\mathbf{T}^{c}$ selbstdual ist, sind f"ur den 
auf die Topologie $\mathbf{T}$ eingeschr"ankten 
Operator $\bigcap\mathbf{T}_{\mbox{id}}|\mathbf{T}$
an jeder Stelle nur zwei Werte m"oglich: Es
ist f"ur alle $X\in \mathbf{T}$ 
$$\bigcap\mathbf{T}_{X}\in\{\emptyset,X\}\ .$$
Diesen Sachverhalt sehen wir auf folgende Weise ein:
Wenn es eine offene Menge ${\rm W}\in \mathbf{T}$ gibt, f"ur die
$${\rm W}\not\supset X\ \land\ {\rm W}\cap X\not=\emptyset$$
zutrifft, so gibt es in der selbstdualen Topologie $\mathbf{T}$ auch das 
Element ${\rm W}^{c}\in \mathbf{T}$, sodass
$$\bigcap\mathbf{T}_{X}=\bigcap\{{\rm U}\in \mathbf{T}:{\rm U}\cap X\not=\emptyset\}$$
$$\subset\bigcap\{{\rm W}\cap X,{\rm W}^{c}\cap X\}=\emptyset$$
ist.
Falls also nicht f"ur alle offenen Mengen ${\rm U}\in \mathbf{T}$
die Alternative
\begin{equation}\label{meire}
{\rm U}\supset X\ \lor\ {\rm U}\cap X=\emptyset
\end{equation}
zutrifft, ist
$\bigcap\mathbf{T}_{X}=\emptyset$.
Andernfalls, falls (\ref{meire}) gilt, ist
$$\mathbf{T}_{X}=\bigcap\{{\rm U}\in \mathbf{T}:{\rm U}\cap X\not=\emptyset\}$$
ein Schnitt "uber lauter Mengen, die $X$ enthalten. Da $X$ selbst in dem Mengensystem 
$\{{\rm U}\in \mathbf{T}:{\rm U}\cap X\not=\emptyset\}$ ist, ist, falls die Alternative (\ref{meire}) gilt,
$\bigcap\mathbf{T}_{X}=X$.
Es gilt also f"ur alle Topologien $\mathbf{T}$ die
Implikation
\begin{equation}\label{meiere}
\begin{array}{c}
\mathbf{T}=\mathbf{T}^{c}\ \land\ {\rm U}\in\mathbf{T}\\
\Rightarrow\\
\Bigl(\exists\ {\rm W}\in\mathbf{T}:{\rm W}\not\supset {\rm U} \land 
{\rm W}\cap {\rm U}\not=\emptyset\ \Leftrightarrow\ \bigcap\mathbf{T}_{{\rm U}}=\emptyset\Bigr)\ .
\end{array}
\end{equation}
Wenn $\mathbf{T}$ eine selbstduale Topologie ist, ist dabei 
das Mengensystem 
\begin{equation}\label{meier}
\mathbf{E}(\mathbf{T}):=\Bigl\{\bigcap\mathbf{T}_{X}:X\in \mathbf{T}\Bigr\} \in\ 
\mathbf{part}\Bigl(\bigcup\mathbf{E}(\mathbf{T})\Bigr)
\end{equation}
deshalb partitiv, weil die Menge
$\bigcap\mathbf{T}_{X}$ gem"ass (\ref{meire}) nur in dem Fall nicht leer ist, in dem 
$X$ zu jedem anderen Element der Topologie disjunkt ist oder in demselben liegt:
Wenn $\bigcap\mathbf{T}_{X(1)}\not=\emptyset$ und $\bigcap\mathbf{T}_{X(2)}\not=\emptyset$
f"ur $X(1),X(2)\in \mathbf{T}$ gilt,
dann ist die Aussage
\begin{displaymath}
\begin{array}{c}
(X(1)\supset X(2)\  \lor\ X(1)\cap X(2)=\emptyset)\ \land\\ 
(X(2)\supset X(1)\  \lor\ X(2)\cap X(1)=\emptyset)
\end{array}
\end{displaymath}
wahr und
$$X(1)=X(2)\  \lor\  X(1)\cap X(2)=\emptyset\ .$$
Das Mengensystem (\ref{meier}) ist zwar in dem Sinn partitiv, dass es 
lauter paarweise disjunkte Mengen als Elemente hat; 
es ist aber denkbar, dass es den topologischen Raum $(\bigcup\mathbf{T},\mathbf{T})$ nicht "uberdeckt.
Die Menge
\begin{equation}\label{maier}
{\rm R}(\mathbf{T}):=(\bigcup\mathbf{T})\setminus{\rm E}(\mathbf{T})
\end{equation}
k"onnte nicht leer sein. In diesem Fall w"are die Spurtopologie $\mathbf{T}':=\mathbf{T}\cap {\rm R}(\mathbf{T})$
ebenfalls selbstdual, weil die leere Menge ihr Element w"are und es zu jeder 
ihrer offenen Mengen
${\rm U}'\in \mathbf{T}'$ eine offene Menge
${\rm U}\in \mathbf{T}$
g"abe, f"ur die
$${\rm U}'={\rm U}\setminus\bigcup\mathbf{E}(\mathbf{T})\ ,$$
mithin
$$(\bigcup\mathbf{T}')\setminus {\rm U}'=  {\rm U}^{c}\setminus\bigcup\mathbf{E}(\mathbf{T})\in \mathbf{T}'$$
w"are; sodass f"ur alle Teilmengen $W'\subset\mathbf{T}'$
eine Teilmenge $W\subset\mathbf{T}$ existierte, f"ur die 
$$\bigcup W'=(\bigcup W)\setminus\bigcup\mathbf{E}(\mathbf{T})\ \in \mathbf{T}'\ ,$$
$$\bigcap W'=(\bigcap W)\setminus\bigcup\mathbf{E}(\mathbf{T})\ \in \mathbf{T}'$$
w"are. Das Mengensystem 
$$\mathbf{E}(\mathbf{T}'):=\Bigl\{\bigcap\mathbf{T}'_{X}:X\in \mathbf{T}'\Bigr\}$$
gen"ugte dabei einerseits der Inklusion
$$\mathbf{E}(\mathbf{T}')\subset\mathbf{E}(\mathbf{T})$$
und andererseits der vermittelten Disjunktion, dass
f"ur alle 
$$\Bigl(\bigcap\mathbf{T}_{X(1)},\bigcap\mathbf{T}'_{X(2)}\Bigr)\in\mathbf{E}(\mathbf{T})\times\mathbf{E}(\mathbf{T}')$$
die Disjunktion
$$\bigcap\mathbf{T}_{X(1)}\cap\bigcap\mathbf{T}'_{X(2)}=\emptyset$$
zutrifft, sodass
$$\mathbf{E}(\mathbf{T}')=\{\emptyset\}$$
w"are.
Nach der Implikation (\ref{meiere}) g"abe es daher f"ur alle 
${\rm U}'\in\mathbf{T}'$ eine Menge
${\rm W}'\in\mathbf{T}'$,
f"ur welche die Aussage 
${\rm W}'\not\supset {\rm U}' \land {\rm W}'\cap {\rm U}'\not=\emptyset$
wahr w"are, sodass
keine offene Menge der selbstdualen Topologie $\mathbf{T}'$ einelementig 
sein d"urfte. 
Dabei existierte
wegen der Schnittabgeschlossenheit 
der selbstdualen Topologie $\mathbf{T}'$ das echt in ${\rm U}'$ enthaltene und nicht leere 
Element
${\rm W}'\cap {\rm U}'\in \mathbf{T}'\setminus\{\emptyset\}$. F"ur jedes Element 
$x\in {\rm U}'$ w"are daher genau eine der beiden Mengen
$${\rm V}',\ (\bigcup\mathbf{T}'\setminus{\rm V}')\in\mathbf{T}'$$
eine offene Menge der selbstdualen Topologie $\mathbf{T}'$, 
in der $x$ nicht w"are, exakt welche ${\rm V}(x)$ bezeichne und f"ur die 
$${\rm U}'\cap {\rm V}(x)\in\mathbf{T}'\setminus\{\emptyset\}\ $$
g"alte. Das Mengensystem
$$\Bigl\{{\rm U}(x)\in \mathbf{T}'\setminus\{\emptyset\}:{\rm U}(x)\subset{\rm U}'\land x\not\in {\rm U}(x)\Bigr\}\not=\emptyset$$
w"are also 
f"ur alle Punkte der Menge ${\rm U}'$
nicht leer und die 
Vereinigung "uber es w"are die 
offene Menge der Topologie $\mathbf{T}'$
$$\hat{{\rm U}}(x):=\bigcup\Bigl\{{\rm U}(x)\in \mathbf{T}'\setminus\{\emptyset\}:{\rm U}(x)\subset{\rm U}'
\land x\not\in {\rm U}(x)\Bigr\}\in\mathbf{T}'\setminus\{\emptyset\}\ ,$$
die in der Umgebung ${\rm U}'$ enthalten ist.
Im Widerspruch dazu allerdings, dass
jede offene Menge der selbstdualen Topologie $\mathbf{T}'$ nicht einelementig w"are,
w"are dann f"ur alle $\alpha\in {\rm U}'$ die Menge
$$\bigcap\Bigl\{\hat{{\rm U}}(x):x\in{\rm U}'\setminus\{\alpha\}\Bigr\}=\{\alpha\}\in \mathbf{T}'\setminus\{\emptyset\}$$
einelementig.
Daher ist die 
zun"achst denkbare M"oglichkeit, dass 
die in
(\ref{maier}) festgelegte Differenz $R(\mathbf{T})$ nicht leer ist, inkonsistent, sodass
die Identit"at 
\begin{equation}
\bigcup\mathbf{E}(\mathbf{T})=\bigcup\mathbf{T}
\end{equation}
gilt.
Das Mengensystem $\mathbf{E}(\mathbf{T})$ "uberdeckt also die Menge $\bigcup\mathbf{T}$, sodass
jede offene Menge ${\rm U}\in \mathbf{T}$ in ihrer "Uberdeckung
$$\bigcup\mathbf{E}(\mathbf{T})_{{\rm U}}\supset{\rm U}$$
enthalten ist. G"abe es einen Punkt $y$ in der Differenz
$$(\bigcup\mathbf{E}(\mathbf{T})_{{\rm U}})\setminus{\rm U}\ ,$$
so g"abe es 
wegen der Partitivit"at (\ref{meier})
genau eine Menge
$\bigcap\mathbf{T}_{X}(y)\in\mathbf{E}(\mathbf{T})_{\{y\}}$,
in der $y\not\in {\rm U}$ ist.
Daher zerfiele $\bigcap\mathbf{T}_{X}(y)$ in 
die beiden disjunkten Teilmengen
$(\bigcap\mathbf{T}_{X}(y)\cap{\rm U})\in \mathbf{T}$ und
$(\bigcap\mathbf{T}_{X}(y)\cap{\rm U}^{c})\in \mathbf{T}$, was im Widerspruch
zur Beschaffenheit der Mengen des Mengensystemes
$\mathbf{E}(\mathbf{T})$ st"unde.
Es gilt also
\begin{equation}
\mathbf{E}(\mathbf{T})^{\cup}\cup\{\emptyset\}=\mathbf{T}\ .
\end{equation}
{\bf q.e.d.}
\newline
\newline
Im Beweis des Satzes 1.1 erwies sich der gem"ass (\ref{wahal}) konstruierte
Operator 
$$\bigcap\mathbf{T}_{\mbox{id}}: 2^{\bigcup\mathbf{T}}\to 2^{\bigcup\mathbf{T}}$$
als ein n"utzliches Instrument,\footnote{Mathematiker spielen f"urwahr viele Instrumente. Und "uberdies bauen sie sie
sogar noch selbst. -- Wie?}
das hierbei 
eine 
Spezifizierung
der allgemeinen 
Konstruktion
$$\bigcap\mathcal{A}_{\mbox{id}}:2^{\bigcup\mathcal{A}}\to 2^{\bigcup\mathcal{A}}$$
f"ur ein beliebiges Mengensystem $\mathcal{A}$ ist,
n"amlich die Spezifizierung f"ur den Fall, dass
$\mathbf{T}$ eine Topologie ist.
Wir bezeichnen dabei
jede Menge $Y$ der Potenzmenge 
$2^{\bigcup\mathcal{A}}$ genau dann als eine 
bez"uglich des Mengensystemes $\mathcal{A}$
integere Menge, wenn $Y$ die Gleichung
\begin{equation}
\bigcap\mathcal{A}_{Y}=Y
\end{equation}
erf"ullt, wenn $Y$ also ein Fixpunkt des
auf einem Mengensystem definierten 
Operators $\bigcap\mathcal{A}_{\mbox{id}}$ ist.\index{integere Menge eines Mengensystemes}
Den Satz "uber die Struktur selbstdualer Topologien k"onnen wir damit so paraphrasieren:
Eine Topologie $\mathbf{T}$ ist also genau dann selbstdual, wenn die bez"uglich ihr 
integeren Mengen -- das sind gerade die Elemente des Mengensystemes $\mathbf{E}(\mathbf{T})$ --
ihre Basis sind.
\newline
Auch der Operator $\mathbf{E}$ kann universalisiert 
aufgefasst werden als einer, der auf der Klasse der Mengensysteme definiert ist und 
der jedem Mengensystem $\mathcal{A}$ den Wert
\begin{equation}\label{elema}
\mathbf{E}(\mathcal{A}):=\Bigl\{\bigcap\mathcal{A}_{A}:A\in\mathcal{A}\}
\end{equation}
zuweist. Es gilt demnach
$$\emptyset\in\mathbf{E}(\mathcal{A})\Leftrightarrow \emptyset\in\mathcal{A}\ ,$$
da 
f"ur alle $A\in \mathcal{A}$
$$\bigcap\mathcal{A}_{A}\supset A$$ 
gilt.
$\mathbf{E}$ bildet das jeweilige Mengensystem $\mathcal{A}$ auf das
Mengensystem aller bez"uglich 
$\mathcal{A}$ integeren Mengen ab. Wir bezeichnen $\mathbf{E}$ als den 
Elementarisierungsoperator.\index{Elementarisierungsoperator} 
F"ur jede 
Topologie $\mathbf{T}$ ist der 
Operator $\bigcap\mathbf{T}_{\mbox{id}}$
genau in dem Fall der Selbstdualit"at $\mathbf{T}=\mathbf{T}^{c}$
mit dem H"ullenoperator $\mathbf{cl}_{\mathbf{T}}$ identisch:
F"ur jedwede
Topologie $\mathbf{T}$ sei
\begin{equation}\label{jedwed}
\mathbf{cl}_{\mathbf{T}}:=\bigcap\{{\rm U}\in\mathbf{T}:{\rm U}\supset\mbox{id}\}
\end{equation}
der auf der jeweiligen Potenzmenge $2^{\bigcup\mathbf{T}}$ definierte, gem"ass (\ref{jedwed}) festgelegte
H"ullenoperator, den wir hier den topologischen\index{topologischer H\"ullenoperator}
H"ullenoperator nennen und dessen jeweilige Werte $\mathbf{cl}_{\mathbf{T}}(X)$ wir als die topologische H"ulle der 
Menge $X$ bez"uglich der Topologie $\mathbf{T}$ bezeichnen; manchmal allerdings, wenn die
Kontextbindung dies erlaubt, bezeichnen wir $\mathbf{cl}_{\mathbf{T}}(X)$ auch einfach als 
die H"ulle der 
Menge $X$ oder schlicht als deren Abschluss.
\footnote{Sodass $\mathbf{cl}_{\mbox{id}}$ als ein
universeller, auf der Klasse der Topologien formulierter, abbildungswertiger Operator aufassbar ist.}
Dass 
f"ur jede Topologie $\mathbf{T}$ f"ur alle Teilmengen $X,Y\subset\bigcup\mathbf{T}$ des zu ihr 
geh"origen topologischen Raumes 
die Aussagen 
\begin{equation}\label{ijedwed}
\begin{array}{c}
\mathbf{cl}_{\mathbf{T}}(X\cap Y)\subset\mathbf{cl}_{\mathbf{T}}(X)\cap\mathbf{cl}_{\mathbf{T}}(Y)\ ,\\
\mbox{bzw.}\\
\mathbf{cl}_{\mathbf{T}}(X\cup Y)\supset\mathbf{cl}_{\mathbf{T}}(X)\cup\mathbf{cl}_{\mathbf{T}}(Y)
\end{array}
\end{equation}
gelten, ist viel mehr bekannt, als leicht verifizierbar; 
leicht ist es hingegen, sich anhand naheliegender Beispiele selbstdualer Topologien
transparent zu machen,
dass keine st"arkere Aussage "uber
die H"ulle des Schnittes je zweier Elemente einer beliebigen   
Topologie gemacht werden kann als in (\ref{ijedwed}), wohingegen f"ur 
die H"ulle des Vereinigung je zweier Elemente einer beliebigen die Gleichung
Topologie $\mathbf{T}$
\begin{equation}\label{wwwed}
\mathbf{cl}_{\mathbf{T}}(X\cup Y)=\mathbf{cl}_{\mathbf{T}}(X)\cup\mathbf{cl}_{\mathbf{T}}(Y)
\end{equation}
gilt. Dass diese Gleichung f"ur 
alle 
nat"urlichen 
Topologien gilt, ist hierbei
ein recht vertrauter Sachverhalt. Deren 
allgemeine
G"ultigkeit werden wir im Zuge der Untersuchungen zeigen, die wir im Abschnitt \ref{ddcln} anstellen, in dem 
wir die 
Generalisierungen vorstellen, zu denen uns die Abschlusskommutation\index{Abschlusskommutation} 
leiten wird.
\newline
Es gilt f"ur jede Topologie $\mathbf{T}$, ebenfalls als eine Paraphrase des Satzes "uber die Struktur selbstdualer Topologien, die 
"Aquivalenz 
\begin{equation}
\mathbf{T}=\mathbf{T}^{c}\Leftrightarrow\mathbf{cl}_{\mathbf{T}}=\mathbf{E}(\mathbf{T})_{\mbox{id}}\ .
\end{equation}
Inwiefern selbstduale Topologien einen Grenzfall der
Topologiehaftigkeit eines Mengensystemes darstellen,
bringt
das
folgende Korollar des 
Satzes 1.1 "uber die Struktur selbstdualer Topologien abschliessend
auf den Punkt:
\newline
\newline
{\bf Korollar 1.2:}\newline
{\em Jede selbstduale Topologie $\mathbf{T}=\mathbf{T}^{c}$
ist genau dann keine ${\rm T}_{0}$-Topologie, wenn
sie nicht die Potenzmenge $2^{\bigcup \mathbf{T}}$
ist.}
\newline\newline 
{\bf Beweis:}\newline 
Es w"are
$$\mathbf{T}=\mathbf{T}^{c}=2^{\bigcup \mathbf{T}}\ ,$$
wenn
es keine nicht leere  
integere Menge 
$${\rm U}\in \mathbf{E}(\mathbf{T})\setminus
\Bigl\{\{x\}:x\in\bigcup \mathbf{T}\Bigr\}$$
g"abe, die 
demnach mindestens zwei verschiedene Elemente $x,y\in {\rm U}$ hat,
f"ur die aber das ${\rm T}_{0}$-Axiom verletzt ist.\newline 
{\bf q.e.d.}
\section{Integere Mengen invarianter Topologien}\label{imit}
Wir wissen nun, was die 
kleinsten
flussinvarianten Mengen
sind, weil
wir nun wissen, dass es sie in Gestalt der integeren Mengen
selbstdualer 
Topologien gibt; und dabei wissen wir auch, dass die 
flussinvarianten 
Topologien selbstduale Topologien sind.
Deren 
integere Mengen verf"uhren dazu, zu meinen, dass sie 
Kandidaten 
f"ur Attraktoren im allgemeinen darstellen. Die flussinvarianten Topologien gem"ass (\ref{flinv})
$$\widehat{\mathbf{T}}(\{\Psi^{t}:t\in \mathbb{R}\})$$
sind speziell die gemeinsamen
invarianten Topologien der Autobolismenmenge $\{\Psi^{t}:t\in \mathbb{R}\}$, die 
die Menge der Phasenflussmitglieder $\Psi^{t}$ eines jeweiligen
Zustandsraumes f"ur alle $t\in \mathbb{R}$ ist. 
Diese Autobolismenmenge $\{\Psi^{t}:t\in \mathbb{R}\}$ ist insofern
algebraisch abgeschlossen, als 
$(\{\Psi^{t}:t\in \mathbb{R}\},\circ)$ eine Gruppe ist.
\newline
Wir untersuchen deshalb zun"achst die invarianten Topologien 
von Autobolismen und von homogenen Autobolismenmengen, weil der Begriffsbildung des Attraktors
von den Fixmengen von Autobolismen ausgeht, um dann mit intuitiver Zauberhand 
das Koh"arenzkriterium in Form der Reichhaltigkeitsaussage
(\ref{verm}) hinzuzunehmen.
Anschliessend untersuchen wir im Hinblick auf das 
Koh"arenzkriterium, wie sich zwei Teilmengen $a\subset \chi$ und $b\subset \chi$
eines integeren Elementes
$\chi\in\mathbf{E}(\widehat{\mathbf{T}}(\Xi))$
f"ur eine homogene Autobolismenmenge $\Xi$
unter der Transformation durch Autobolismen der Menge $\Xi$ verhalten.
\newline
Die Fixmengen eines Autobolismus $\xi$ sind die Elemente seiner invarianten Topologie $\mathbf{T}(\xi)=\mathbf{T}(\xi)^{c}$,
die selbstdual ist, sodass es die bez"uglich $\mathbf{T}(\xi)$ integeren Fixmengen des Autobolismus $\xi$
geben muss, welche die kleinsten Fixmengen des Autobolismus $\xi$ sind. Die bez"uglich $\mathbf{T}(\xi)$ integeren Fixmengen des Autobolismus $\xi$
sind vergleichsweise "uberschaubare Konstrukte: Im Vergleich zu den 
bez"uglich $\widehat{\mathbf{T}}(\Xi)$
integeren Elementen der 
ebenfalls selbstdualen
gemeinsamen invarianten Topologie $\widehat{\mathbf{T}}(\Xi)=\widehat{\mathbf{T}}(\Xi)^{c}$
einer homogenen Autobolismenmenge $\Xi$ sind
n"amlich die sogenannten Orbits eines
Autobolismus $\xi$ einfach:
F"ur jedes Element $z\in\mathbf{P}_{1}\xi$ ist die Menge
\begin{equation}\label{orson}
\mathbf{or}_{\xi}(\{z\}):=\{x\in\mathbf{P}_{1}\xi:\exists\ j\in\mathbb{Z}: \xi^{j}(z)=x\}
\end{equation}
offenbar eine Fixmenge des Autobolismus $\xi$, wobei 
f"ur alle $j\in \mathbb{Z}$ die Komposition
$\xi\circ\xi^{j}=\xi^{j}\circ\xi$ der Autobolismus
$\xi^{j+1}$ sei und $\xi^{0}=\mbox{id}$, sodass
$\xi^{-1}$ die Inversion von $\xi$ ist. 
Nach dem 
Satz 1.1 "uber die Struktur selbstdualer Topologien ist demnach das Mengensystem 
\begin{equation}\label{eiere}
\mathbf{T}(\xi)=\Bigl\{\mathbf{or}_{\xi}(\{z\}):z\in\mathbf{P}_{1}\xi\Bigr\}^{\cup}\ \cup\ \{\emptyset\}
\end{equation}
die invariante Topologie $\mathbf{T}(\xi)$,
falls folgende beiden Bedingungen erf"ullt sind:\newline
1. F"ur alle $z\in \mathbf{P}_{1}\xi$ ist die Fixmenge $\mathbf{or}_{\xi}(\{z\})=\xi(\mathbf{or}_{\xi}(\{z\})$
integer bez"uglich
der invarianten Topologie $\mathbf{T}(\xi)$.\newline
2. Ferner hat "uberdies folgende Form der Vollst"andigkeit vorzuliegen: Alle Fixmengen des Autobolismus $\xi$ sind 
partionierbar in Orbits $\mathbf{or}_{\xi}(\{z\})=\xi(\mathbf{or}_{\xi}(\{z\})$ 
der Form (\ref{orson})
f"ur einen jeweiligen Zustand $z$ des Zustandsraumes $\mathbf{P}_{1}\xi$.
Diese zweite Bedingung ist offensichtlich erf"ullt:
Jeder Zustand $z\in\mathbf{P}_{1}\xi$ liegt in einem Orbit.\footnote{
Es ist doch denkbar, dass der Zustandsraum $\mathbf{P}_{1}\xi$ eine weitere
Topologie $\mathbf{T}(\mathbf{P}_{1}\xi)$ tr"agt, die beispielsweise
eine ${\rm T}_{2}$-Topologie sei, bez"uglich welcher der Zustandsraum $\mathbf{P}_{1}\xi$
kompakt ist, sodass zu jedem nicht endlichen Orbit $\mathbf{or}_{\xi}(\{z\})$ f"ur 
einen jeweiligen Zustand $z\in\mathbf{P}_{1}\xi$
Ber"uhrungsspunkte $z^{\star}$ existieren. Solche Ber"uhrungsspunkte $z^{\star}$ sind im allgemeinen 
wegen der Invertierbarkeit des Autobolismus $\xi$
keine Punkte des
Orbits $\mathbf{or}_{\xi}(\{z\})\not\ni z^{\star}$. Es lassen sich aber sogleich
Autobolismen $\xi^{\star}$ formulieren, bei denen ein Ber"uhrungsspunkt $z^{\star}$ eines Orbits
$\mathbf{or}_{\xi}(\{z\})\not\ni z^{\star}$ in dem Orbit des abgewandelten Autobolismus $\xi^{\star}$
$$\mathbf{or}_{\xi^{\star}}(\{z\})=\mathbf{or}_{\xi^{\star}}(\{z^{\star}\})\supset\mathbf{or}_{\xi}(\{z\})\cup
\mathbf{or}_{\xi}(\{z^{\star}\})$$
liegt. Der Orbit
$\mathbf{or}_{\xi^{\star}}(\{z\})$ l"asst sich insofern als quasiindeterministischer Orbit auffassen.\index{quasiindeterministischer Orbit}
Der Gedanke, dass jene Ber"uhrungsspunkte $z^{\star}$ eines Orbits
$\mathbf{or}_{\xi}(\{z\})$ Fixpunkte des 
Autobolismus $\xi$ seien, ist,
gelinde gesagt, Leichtsinn:
Nichtsdestotrotz verl"auft durch jeden derartigen Ber"uhrungsspunkt $z^{\star}$ 
ein Orbit, n"amlich der
Orbit $\mathbf{or}_{\xi}(\{z^{\star}\})$.}
Hierbei die mit der Wertemenge 
des Autobolismus $\xi$ 
identische Definitionsmenge desselben $\mathbf{P}_{1}\xi=\mathbf{P}_{2}\xi$ als Zustandsraum\index{Zustandsraum}
zu bezeichnen, ist dadurch gerechtfertigt, dass der Autobolismus $\xi$
zu dem diskreten dynamischen System\index{diskretes dynamisches System} 
$$(\xi^{\mathbb{Z}},\circ)$$
"aquivalent ist.
Das Integrit"atskriterium\index{Integrit\"atskriterium\} an ein Element ${\rm U}\in \mathbf{T}(\xi)$
der invarianten Topologie als einer selbstdualen Topologie
ist hierbei, dass gem"ass (\ref{meiere})
\begin{equation}\label{aeiere}
\forall\ {\rm W}\in\mathbf{T}(\xi)\ {\rm W}\supset {\rm U} \lor{\rm W}\cap {\rm U}=\emptyset 
\end{equation}
gilt.
Da f"ur alle Orbits $\omega_{1},\omega_{2}\in \{\mathbf{or}_{\xi}(\{z\}):z\in\mathbf{P}_{1}\xi\}$
wegen der Bijektivit"at des Autobolismus $\xi$ Partitivit"at in Form der "Aquivalenz
\begin{equation}
\omega_{1}\cap\omega_{2}\not=\emptyset\Leftrightarrow \omega_{1}=\omega_{2}
\end{equation}
vorliegt, ist das Integrit"atskriterium\index{Integrit\"atskriterium} (\ref{aeiere}) erf"ullt,
wenn
jede nicht leere, echte Teilmenge $q\subset\omega$ eines Orbits $\omega$ keine 
Fixmenge des Autobolismus $\xi$ ist. Falls weder $\omega\setminus q$ noch
$q$ leer ist, gibt
es das Element $z\in \mathbf{P}_{1}\xi$, f"ur das $\omega=\mathbf{or}_{\xi}(\{z\})$ ist.
Des weiteren gibt
es 
eine entweder einelementige 
oder eine endliche, aus aufeinanderfolgenden, ganzen Zahlen bestehende 
Teilmenge 
$X\subset \mathbb{Z}$, f"ur welche die Menge
$\{j\in X: \xi^{j}(z)=x\}$ in $q$ liegt
und f"ur die auserdem zutrifft, dass entweder
das Element
$\xi^{\max X+1 }(z)$ nicht in $q$ ist; oder es ist der Fall, dass das Element
$\xi^{\min X-1 }(z)$ nicht in $q$ ist.
Sonst w"are $q$ leer oder mit $\omega$ identisch. Also ist 
jede nicht leere, echte Teilmenge $q\subset\omega$ eines Orbits $\omega$ keine
Fixmenge des Autobolismus $\xi$ und $\omega$ ist daher integer.
Die Aufz"ahlung (\ref{eiere}) beschreibt die invariante Topologie $\mathbf{T}(\xi)$
also in der Tat.
\newline
Wie aber sieht die gemeinsame invariante Topologie $\widehat{\mathbf{T}}(\{\xi_{1},\xi_{2}\})$
einer homogenen Autobolismenmenge $\{\xi_{1},\xi_{2}\}$ zweier Autobolismen 
$\xi_{1}$ und $\xi_{2}$ aus, die auf der gemeinsamen Definitionsmenge
$\mathbf{P}_{1}\xi_{1}=\mathbf{P}_{1}\xi_{2}$ erkl"art sind?
Die gemeinsame Definitionsmenge
$\mathbf{P}_{1}\xi_{1}=\mathbf{P}_{1}\xi_{2}$ 
bezeichnen wir hierbei als Zustandsraum,\index{Zustandsraum}
obwohl 
die homogene Autobolismenmenge $\{\xi_{1},\xi_{2}\}$
exakt in dem Fall algebraisch
abgeschlossen ist, dass einer der beiden 
Autobolismen die Idenitit"at
ist und der jeweils andere Autobolismus
eine selbstinverse Involution.
Im allgemeinen liegt also mit der Menge 
$\{\xi_{1},\xi_{2}\}$ in dem Paar
$(\{\xi_{1},\xi_{2}\},\circ)$ noch keine Gruppe vor.
Das Paar $(\{\xi_{1},\xi_{2}\},\circ)$ k"onnen wir
im Allgemeinfall
daher beim besten Willen nicht als dynamisches System
ansehen.\newline  
Wir wollen n"amlich die
allgemeine Auffassung einer Phasenflussgruppe und eines dynamischen Systemes teilen:
Wenn $\Xi$ eine Menge von Bijektionen $\xi$ einer Menge $\zeta$
auf sich selbst ist, so ist das Paar $(\Xi,\circ)$ genau dann eine Gruppe,
wenn $\Xi$ abgeschlossen gegen"uber der Komposition $\circ$ je zweier 
seiner Elemente ist. Exakt jede Gruppe der Form
$(\Xi,\circ)$ nennen wir hier eine Phasenflussgruppe\index{Phasenflussgruppe} und bezeichnen 
dabei all deren Gruppenelemente als Phasenflussmitglieder, wohingegen wir
exakt die speziellen Phasenflussgruppen $(\{\Phi^{t}:t\in\mathbb{R}\},\circ)$
als kontinuierliche Phasenflussgruppen bezeichnen, die uns in diesem Trakat 
durchg"angig
die 
kontinuierlichen dynamischen Systeme 
objektivieren.\index{kontinuierliche Phasenflussgruppe}\index{kontinuierliches dynamisches System}
Da aber andererseits durch die homogene Autobolismenmenge $\{\xi_{1},\xi_{2}\}$
bereits deren algebraischer Abschluss $\langle\{\xi_{1},\xi_{2}\} \rangle$
gem"ass (\ref{auss})
festgelegt ist und $(\langle\{\xi_{1},\xi_{2}\} \rangle,\circ)$ sehr wohl
eine Gruppe ist, die wir 
konform der allgemeinen Auffassung als ein dynamisches System
ansehen d"urfen, ist die Bezeichnung der gemeinsamen Definitionsmenge
$\mathbf{P}_{1}\xi_{1}=\mathbf{P}_{1}\xi_{2}$
als Zustandsraum gerechtfertigt.
Auf die entsprechende Weise k"onnen wir auch rechtfertigen, die
gemeinsame Definitionsmenge jedweder homogenen Autobolismenmenge
als deren Zustandsraum gelten zu lassen.   
Setzen wir unsere Betrachtung der gemeinsamen invarianten 
Topologie $\widehat{\mathbf{T}}(\{\xi_{1},\xi_{2}\})$
einer homogenen Autobolismenmenge $\{\xi_{1},\xi_{2}\}$ zweier Autobolismen 
$\xi_{1}$ und $\xi_{2}$ fort, die 
den Zustandsraum 
$$\mathbf{P}_{1}\xi_{1}=\mathbf{P}_{2}\xi_{1}=
\mathbf{P}_{1}\xi_{2}=\bigcup\widehat{\mathbf{T}}(\{\xi_{1},\xi_{2}\})$$
der Autobolismenmenge $\{\xi_{1},\xi_{2}\}$ auf partitve Weise
topologisiert:
$\widehat{\mathbf{T}}(\{\xi_{1},\xi_{2}\})$ ist zwar einfach der Schnitt 
$$\widehat{\mathbf{T}}(\{\xi_{1},\xi_{2}\})=\mathbf{T}(\xi_{1})\cap\mathbf{T}(\xi_{2})$$
$$\Bigl\{\mathbf{or}_{\xi_{1}}(\{z\}):z\in\mathbf{P}_{1}\xi_{1}\Bigr\}^{\cup}\cap\Bigl\{\mathbf{or}_{\xi_{2}}(\{z\}):z\in\mathbf{P}_{1}\xi_{2}\Bigr\}^{\cup}\ \cup\ \{\emptyset\}\ .$$
Dessen Elemente sind aber Widerspiegelungen der 
algebraischen Beschaffenheit der von den beiden
Autobolismen $\xi_{1}$ und $\xi_{2}$ abz"ahlbar erzeugten 
Autobolismengruppe
\begin{equation}
(\langle\{\xi_{1},\xi_{2}\} \rangle,\circ)\ ,
\end{equation}
deren erste Komponente 
\begin{equation}\label{auss}
\langle\{\xi_{1},\xi_{2}\} \rangle
\end{equation}
$$:=\Bigl\{\alpha\in\mathbf{P}_{1}\xi_{1}^{\mathbf{P}_{1}\xi_{1}}:\exists\ n\in\mathbb{N},(k_{1},k_{2},\dots k_{n}),
(j_{1},j_{2},\dots j_{n})\in\mathbb{Z}^{n}:$$
$$\alpha=\xi_{1}^{k_{1}}\circ\xi_{2}^{j_{1}}\circ\xi_{1}^{k_{2}}\circ\xi_{2}^{j_{2}}
\circ\dots\xi_{1}^{k_{n}}\circ\xi_{2}^{j_{n}}\Bigr\}$$
wohldefiniert ist, weil in dieser Festlegung 
$\xi_{1}$ und $\xi_{2}$ offensichtlich vertauschbar sind.
Die nicht leeren Elemente 
der gemeinsamen invarianten Topologie $\widehat{\mathbf{T}}(\{\xi_{1},\xi_{2}\})$
sind, weil diese der Schnitt $\mathbf{T}(\xi_{1})\cap\mathbf{T}(\xi_{2})$ ist, gerade folgende 
Teilmengen ${\rm U}$ der gemeinsamen Definitionsmenge
$\mathbf{P}_{1}\xi_{1}=\mathbf{P}_{1}\xi_{2}$, 
die der Zustandsraum ist: 
Exakt die Teilmengen ${\rm U}$ des Zustandsraumes,
die gleichzeitige Partionierbarkeit in dem Sinn haben, dass sie sowohl 
in Orbits 
$\mathbf{or}_{\xi_{1}}(\{z_{1}\})$ des Autobolismus $\xi_{1}$ f"ur Zust"ande $z_{1}\in Z_{1}\subset\mathbf{P}_{1}\xi_{1}$
als auch in
Orbits 
$\mathbf{or}_{\xi_{2}}(\{z_{2}\})$ des Autobolismus $\xi_{2}$
f"ur Zust"ande $z_{2}\in Z_{2}\subset\mathbf{P}_{1}\xi_{1}$ partionierbar sind,
sind die nicht leeren Elemente 
der gemeinsamen invarianten Topologie $\widehat{\mathbf{T}}(\{\xi_{1},\xi_{2}\})$.
Es ist
$$\widehat{\mathbf{T}}(\{\xi_{1},\xi_{2}\})=\Bigl\{Q\in 2^{\mathbf{P}_{1}\xi_{1}}:\exists
Z_{1}\times Z_{2}\subset \mathbf{P}_{1}\xi_{1}\times\mathbf{P}_{1}\xi_{1}:$$
$$Q=\bigcup\{\mathbf{or}_{\xi_{1}}(\{z_{1}\}):z_{1}\in Z_{1}\}=
\bigcup\{\mathbf{or}_{\xi_{2}}(\{z_{2}\}):z_{2}\in Z_{2}\}\Bigr\}\ .$$
Die gleichzeitig in Orbits 
des Autobolismus $\xi_{1}$ 
und in
Orbits 
des Autobolismus $\xi_{2}$ partionierbaren Mengen $Q$ der  
gemeinsamen invarianten Topologie $\widehat{\mathbf{T}}(\{\xi_{1},\xi_{2}\})$ sind
also 
f"ur alle Paare ganzer Zahlen $(j,k)\in \mathbb{Z}\times\mathbb{Z}$
invariant gegen"uber der Abbildung durch jeden 
Autobolismus
der Menge $\langle\{\xi_{1},\xi_{2}\} \rangle$.                                         )
Daher gibt es f"ur jede Menge $Q$ der  
gemeinsamen invarianten Topologie $\widehat{\mathbf{T}}(\{\xi_{1},\xi_{2}\})$
eine Menge von $Z(Q)\subset\mathbf{P}_{1}\xi_{1}$, f"ur die
$$Q=\bigcup\{\alpha(z):(\alpha, z)\in \langle\{\xi_{1},\xi_{2}\} \rangle\times Z(Q)\}$$
ist. Dabei ist $Q$
offenbar nicht bez"uglich $\widehat{\mathbf{T}}(\{\xi_{1},\xi_{2}\})$ integer, wenn es
einen Zustand $z\in \mathbf{P}_{1}\xi_{1}$ gibt, f"ur den zwar die Inklusion
\begin{equation}\label{muver}
\bigcup\langle\{\xi_{1},\xi_{2}\} \rangle(z):=\bigcup\{\alpha(z):\alpha\in \langle\{\xi_{1},\xi_{2}\} \rangle\}\subset Q
\end{equation}
gilt, jedoch ohne, dass deren linke Seite 
mit $Q$ "ubereinstimmt; wobei jede nicht leere, echt in $\bigcup\langle\{\xi_{1},\xi_{2}\} \rangle(z(Q))$ enthaltene Teilmenge
offensichtlich nicht invariant gegen"uber allen Autobolismen der Menge $\langle\{\xi_{1},\xi_{2}\} \rangle$ ist. Die
integeren Elemente der selbstdualen Topologie $\widehat{\mathbf{T}}(\{\xi_{1},\xi_{2}\})$ sind
also gerade die Mengen der Form gem"ass (\ref{muver}) und es ist
\begin{equation}\label{maver}
\widehat{\mathbf{T}}(\{\xi_{1},\xi_{2}\})=\Bigl\{\bigcup\langle\{\xi_{1},\xi_{2}\} \rangle(z):z\in \mathbf{P}_{1}\xi_{1}\Bigr\}\ \cup\ \{\emptyset\}
\end{equation}
$$=\bigcup\bigcup\langle\{\xi_{1},\xi_{2}\} \rangle(\mathbf{P}_{1}\xi_{1})\ \cup\ \{\emptyset\}\ .$$
Denn die Vollst"andigkeit 
dieses Mengensystemes liegt in diesem Sinn
vor: Alle Fixmengen des Autobolismus $\xi$ sind 
partionierbar in Biorbits der Form
\begin{equation}
\mathbf{or}_{\xi_{1},\xi_{2}}(\{z\})=\bigcup\langle\{\xi_{1},\xi_{2}\} \rangle(z)
\end{equation}
f"ur einen jeweiligen Zustand $z$ des Zustandsraumes $\mathbf{P}_{1}\xi$, weil jeder 
Zustand $z\in\mathbf{P}_{1}\xi$ in einem Biorbit ist.
\newline
Aus der expliziten Aufz"ahlung 
der Elemente einer gemeinsamen invarianten Topologie $\widehat{\mathbf{T}}(\{\xi_{1},\xi_{2}\})$
einer homogenen zweielementigen Autobolismenge
(\ref{maver}) ergibt sich auch schon die Beschaffenheit der gemeinsamen invarianten Topologien 
$\widehat{\mathbf{T}}(\Xi)$
endlicher homogener
Autobolismengen $\Xi=\{\xi_{1},\xi_{2},\dots \xi_{\nu}\}$ f"ur $\nu\in\mathbb{N}$: Wenn
\begin{equation}\label{ausss}
\langle \Xi\rangle=\langle\{\xi_{1},\xi_{2},\dots  \xi_{\nu}\}\rangle
\end{equation}
$$:=\Bigl\{\alpha\in\mathbf{P}_{1}\xi_{1}^{\mathbf{P}_{1}\xi_{1}}:\exists\ n\in\mathbb{N},(k_{1,1},k_{2,1},\dots k_{n,1}),
(k_{1,2},k_{2,2},\dots k_{n,2}),\dots$$ $$(k_{1,\nu},k_{2,\nu},\dots k_{n,\nu})\in\mathbb{Z}^{n}:$$
$$\alpha=\xi_{1}^{k_{1,1}}\circ\xi_{2}^{k_{1,2}}\circ\dots\xi_{\nu}^{k_{1,\nu}}\circ
\xi_{1}^{k_{2,1}}\circ\xi_{2}^{k_{2,2}}\circ\dots\xi_{\nu}^{k_{2,\nu}}\circ\dots
\xi_{1}^{k_{n,1}}\circ\xi_{2}^{k_{n,1}}\circ\dots\xi_{\nu}^{k_{n,\nu}}\Bigr\}$$
die Tr"agermenge der von den Autobolismen der Menge $\Xi$ abz"ahlbar erzeugten Gruppe ist, ist
\begin{equation}\label{maver}
\widehat{\mathbf{T}}(\Xi)=\bigcup\bigcup\langle\Xi\rangle(\mathbf{P}_{1}\xi_{1})\ \cup\ \{\emptyset\}
\end{equation}
die gemeinsame invariante Topologie der endlichen homogenen
Autobolismenge $\Xi$.
Die Pendants zu den abz"ahlbar erzeugten 
algebraischen 
H"ullen $\langle\{\xi_{1},\xi_{2},\dots  \xi_{\nu}\}\rangle$ k"onnen wir in dem Fall, dass
$\Xi$ keine endliche Menge homogener Autobolismen ist, nicht mehr wie in
(\ref{auss}) und in (\ref{ausss})
ausschreiben. Die abz"ahlbar erzeugten algebraischen 
H"ullen $\langle\{\xi_{1},\xi_{2},\dots  \xi_{\nu}\}\rangle$
stimmen mit den Werten
"uberein, die ihnen 
der allgemeine, auf der Klasse der homogenen Autobolismenmengen $\Xi\subset\mathbf{ab}(A)\subset A^{A}$ einer Menge 
$$A=\downarrow\mathbf{P}_{1}\Xi=\downarrow\{\mathbf{P}_{1}\xi:\xi\in\Xi\}$$
definierte 
H"ullenoperator
\begin{equation}\label{maaver}
\langle\mbox{id}\rangle:=\bigcap
\Bigl\{X\in\mathbf{ab}\ \circ\downarrow\mathbf{P}_{1}:
X\supset\mbox{id}\ \land X\circ X=X\Bigr\}
\end{equation}
zuordnet, exakt welchen wir 
als den algebraischen H"ullenoperator\index{algebraischer H\"ullenoperator} einer Autobolismenmenge bezeichnen. Dabei ist $\downarrow\mbox{id}$ der auf der Klasse einelementiger Mengen erkl"arte
Operator,
dessen jeweiliger Wert f"ur jede einelementige Menge ${\rm E}=\{e\}$
deren Element
\begin{equation}\label{wahn}
\downarrow{\rm E}:=e
\end{equation} sei. 
$\mathbf{ab}$ ist als ein mengenwertiger Operator zu lesen, der 
jeder Menge
$A$ als jeweiligen Wert $\mathbf{ab}(A)$
die Menge aller Autobolismen der Menge $A$ zuordnet.
Den gem"ass (\ref{maaver}) verfassten algebraischen H"ullenoperator $\langle\mbox{id}\rangle$ k"onnen wir
auf beliebige homogene Autobolismenmengen $\Xi$ jeglicher Kardinalit"at 
anwenden. Er gibt uns offensichtlich allemal 
die Tr"agermenge der von den Autobolismen der Menge $\Xi$ abz"ahlbar erzeugten Gruppe.
Wir verifizieren nun auch leicht mit Hilfe des Satzes 1.1 "uber die Struktur selbstdualer Topologien, dass f"ur jede homogene Autobolismenmenge $\Xi$ deren
invariante Topologie $\widehat{\mathbf{T}}(\Xi)$ die Gleichung (\ref{maver}) erf"ullt.
\newline
Wir orientierten uns hier bei der Aufz"ahlung
der Elemente einer gemeinsamen invarianten Topologie $\widehat{\mathbf{T}}(\Xi)$ 
einer homogenen Autobolismenge an dem Befund 
des Satzes 1.1 "uber die Struktur selbstdualer Topologien,
dass jede selbstduale Topologie eine Basis bez"uglich ihr
integerer Mengen hat. Diese Orientierung sowohl an dem Sachverhalt, dass
es die integeren Mengen einer jeweiligen selbstdualen Topologie "uberhaupt gibt, als auch an dem
Wissen, dass die jeweiligen integeren Mengen das 
Integrit"atskriterium\index{Integrit\"atskriterium} (\ref{aeiere}) erf"ullen m"ussen,  
erwies sich hierbei zwar als durchaus hilfreich. Wir sehen
aber, dass wir in den F"allen 
der 
integeren Mengen der 
gemeinsamen invarianten Topologien
endlicher homogener
Autobolismengen auch ohne
die besagte Orientierungshilfe ausgekommen w"aren.
Deren integere Menge
sind dabei abz"ahlbar; auch die integeren Mengen der 
gemeinsamen invarianten Topologien
abz"ahlbarer homogener
Autobolismengen 
sind abz"ahlbar.
\newline
Wir bemerken dabei, dass
f"ur uns die Strukturbestimmung 
selbstdualer Topologien
durch den Satz 1.1
um so hilfreicher ist, je un"uberschaubarer die von 
jeweiligen homogenen
Autobolismengen festgelegten gemeinsamen invarianten Topologien sind.\newline
Damit haben wir die integeren Mengen der invarianten Topologien
jeweiliger homogener
Autobolismengen bestimmt. Die Beschaffenheit deren integerer Mengen ist
offenbar weitaus weniger schlicht als der Aufbau
der invarianten Topologien
aus deren integeren Mengen.\newline
Unserer Aufgabe folgend fragen wir 
nun im Hinblick auf das Koh"arenzkriterium (\ref{verm}),
f"ur welche Autobolismenmengen $\Xi$ sich zwei Teilmengen $a\subset \chi$ und $b\subset \chi$
eines integeren Elementes
$$\chi\in\mathbf{E}(\widehat{\mathbf{T}}(\Xi))$$
durch einen Autobolismus der Menge $\Xi$ so transformieren
lassen, dass sie einander schneiden.
$\mathbf{E}$ ist hierbei
der Elementarisierungsoperator gem"ass (\ref{elema}).\index{Elementarisierungsoperator} 
Das folgende Lemma gibt uns 
hier Aufschluss und eine Antwort auf die eingangs gestellte \glqq dumme Frage\grqq,\index{dumme Frage}
die uns an E.K"astners Feuerzangenbowle gemahnte; n"amlich,
ob es sein kann, dass die nunmehr als integere Mengen der gemeinsamen invarianten Topologien
objektivierten
kleinsten Mengen der invarianten Topologien dem Koh"arenzkriterium nicht
gen"ugen:
\newline 
\newline
{\bf Lemma 1.3: Indifferente Koh"arenz}\newline
{\em Es sei $\Xi$ eine homogene Autobolismenmenge und $\widehat{\mathbf{T}}(\Xi)$ deren
gemeinsame invariante Topologie und 
$\chi\subset(\bigcup\widehat{\mathbf{T}}(\Xi))\setminus\{\emptyset\}$ eine nicht leere
Teilmenge des Zustandsraumes. 
Es gilt die "Aquivalenz}
\begin{equation}\label{wanhn}
\begin{array}{c}
\chi\in\mathbf{E}(\widehat{\mathbf{T}}(\Xi))\\
\Leftrightarrow\\
\Bigl(a,b\subset\chi\Rightarrow\ \exists\ \xi\in\langle\Xi\rangle:\ 
\xi(a)\cap b\not=\emptyset\Bigr)\ .\\
\end{array}
\end{equation} 
\newline
{\bf Beweis:}\newline
Die Implikationsrichtung von links nach rechts ergibt sich folgendermassen:
F"ur alle integeren Mengen $\chi\in\mathbf{E}(\widehat{\mathbf{T}}(\Xi))$ der 
Topologie $\widehat{\mathbf{T}}(\Xi)$ gibt es
nach den vorangegangenen "Uberlegungen, die zu dem Resultat leiteten, dass
jede gemeinsame invariante Topologie $\widehat{\mathbf{T}}(\Xi)$
einer homogenen Autobolismenmenge $\Xi$
der Gleichung (\ref{maver}) gen"ugt,
einen Zustand $x\in \bigcup\widehat{\mathbf{T}}(\Xi)$, f"ur den
$$\chi=\bigcup\langle\Xi\rangle(x)=\{\xi(x):\xi\in \langle\Xi\rangle\}$$
ist. F"ur jedes Punktepaar
$$(\alpha,\beta)\in\chi\times\chi$$
gibt es daher ein Paar zweier Autobolismen
$$(\xi_{\alpha}, \xi_{\beta})\in \langle\Xi\rangle\times\langle\Xi\rangle\ ,$$
f"ur das
$$\alpha=\xi_{\alpha}(x)\ ,$$
$$\beta=\xi_{\beta}(x)$$
ist, sodass
$$\alpha=\xi_{\alpha}\xi_{\beta}^{-1}(\beta)$$
ist, wobei
$$\xi_{\alpha}\xi_{\beta}^{-1}\in \langle\Xi\rangle$$
ist.
Die andere Implikationsrichtung zeigen wir indirekt:
Ist $\chi\subset\bigcup\widehat{\mathbf{T}}(\Xi)$ keine integere
Menge der Topologie $\widehat{\mathbf{T}}(\Xi)$ und nicht leer, so
ist wegen der 
von dem Satz 1.1 "uber die Struktur selbstdualer Mengen behaupteten 
Partititivit"at des Mengensystemes der integeren Mengen
die Ungleichung
$$\mathbf{card}(\mathbf{E}(\widehat{\mathbf{T}}(\Xi))_{\chi})>1$$
wahr und es gibt 
zwei verschiedene integere Mengen
$e_{1},e_{2}\in \mathbf{E}(\widehat{\mathbf{T}}(\Xi))$,
f"ur die
$$\chi\cap e_{1}\not=\emptyset\ ,$$
$$\chi\cap e_{2}\not=\emptyset\ ,$$
$$(\chi\cap e_{1})\cap(\chi\cap e_{2})=\emptyset$$
gilt. Also gibt es zwei Teilmengen $a$ und $b$ der Menge $\chi$, 
die in $\chi\cap e_{1}\supset a$ und in $\chi\cap e_{2}\supset b$
liegen. Wenn es einen Autobolismus $\xi\in \langle\Xi\rangle$ g"abe, f"ur den das
Bild $\xi(a)$
zu $b$ nicht disjunkt w"are,
g"abe es zwei Zust"ande $\beta\in e_{2}$ und $\alpha\in e_{1}$,
f"ur die $$\xi(\alpha)=\beta$$
w"are, was der Eigenschaft der integeren Menge $e_{1}$,
Fixmenge aller Autobolismen der Menge $\Xi$ zu sein, widerspricht.\newline
{\bf q.e.d.}
\newline
\newline
Immerhin,\index{dumme Frage}
das wissen wir nun: Wenn es nicht die integeren Mengen der Menge $\mathbf{E}(\widehat{\mathbf{T}}(\Xi))$
f"ur die gemeinsame invariante Topologie $\widehat{\mathbf{T}}(\Xi)$ sind, die
dem Koh"arenzkriterium gen"ugen, so gibt es keine Attraktoren, die der homogen Autobolismenmengen
$\Xi$ zugewiesen werden k"onnten.
\newline
Die "Aquivalenz (\ref{wanhn}) impliziert offensichtlich, 
dass f"ur alle integeren Mengen $\chi_{1},\chi_{2}\in\mathbf{E}(\widehat{\mathbf{T}}(\Xi))$
der Topologie $\widehat{\mathbf{T}}(\Xi)$ die im Vergleich zu ihr schw"achere "Aquivalenz
\begin{equation}\label{wwnhn}
\begin{array}{c}
\Bigl(a\subset\chi_{1}\ \land\ b\subset\chi_{2}\ \land \exists\ \xi\in\langle\Xi\rangle:\ \xi(a)\cap b\not=\emptyset\Bigr)\\
\Leftrightarrow\\
\chi_{1}=\chi_{2}
\end{array}
\end{equation} 
gilt. Ferner beinhaltet die "Aquivalenz (\ref{wanhn})
offensichtlich auch die im Vergleich zu ihr schw"achere
Implikation
\begin{equation}\label{anhn}
\begin{array}{c}
\neg\Bigl(a,b\subset\chi\ \land \chi\in\mathbf{E}(\widehat{\mathbf{T}}(\Xi))\setminus\{\emptyset\}\Rightarrow\ \exists\ 
\xi\in\Xi:\ \xi(a)\cap b\not=\emptyset\Bigr)\\
\Rightarrow  \Xi\not=\langle\Xi\rangle\ ,
\end{array}
\end{equation} 
die nicht umkehrbar ist, wie das Beispiel der 
zweielementigen Permutationsgruppe $\mathbf{S}_{2}=({\rm S}_{2},\circ)$ zeigt, die
auf dem Zustandsraum $\{1,2\}$ konkretisiert sei, welcher also f"ur beide Permutationen 
$\xi\in{\rm S}_{2}$ mit
$\mathbf{P}_{1}\xi$ "ubereinstimmt: F"ur 
die nicht leere 
Menge der invarianten Topologie  $\{\emptyset,\{1,2\}\}=\mathbf{T}({\rm S}_{2}\setminus\{\mbox{id}\})$ 
gibt es die selbstinverse
Permutation in ${\rm S}_{2}\setminus\{\mbox{id}\}$, die
die Teilmenge $\{1\}\subset\{1,2\}$ des
Orbits $\{1,2\}$ auf die Teilmenge $\{2\}\subset\{1,2\}$
abbildet. Sogleich werden wir das klassische Beispiel f"ur die 
Unumkehrbarkeit der Implikation
(\ref{anhn}) vorf"uhren.
\newline
Nichtsdestotrotz bietet sich
die 
in der "Aquivalenz (\ref{wanhn}) beinhaltete, schw"achere
Implikation als negatives Kriterium f"ur die algebraische 
Abgeschlossenheit jeweiliger homogener Autobolismenmengen an.
Betrachten wir dazu diejenigen Funktionen, die uns hier als Wellenfunktionen gelten sollen im Vergleich zu den
Flussfunktionen:\newline 
Damit dabei keine Unklarheit "uber den Begriff der Flussfunktion bestehe, erl"autern wir, dass wir hier
jede Abbildung $\Phi$ eine Wellenfunktion nennen, f"ur die es zwei nicht leere
Mengen 
$A$ und $B$ gibt,
die so beschaffen sind, dass
\begin{equation}\label{well}
\Phi\in B^{A\times \mathbb{R}}
\end{equation} 
gilt. Wir sollen hier warnen: Dies ist n"amlich eine sehr allgemeine Auffassung der 
Wellenfunktion,\index{Wellenfunktion} die nicht allgemein verfestigt ist.
Auf der Klasse der Wellenfunktionen sei 
der Kollektivierungsoperator\index{Kollektivierungsoperator} $[\mbox{id}]$
erkl"art, dessen jeweiliger Wert f"ur die Wellenfunktion $\Phi$
das Mengensystem 
\begin{equation}
[\Phi]:=\{\Phi(x,\mathbb{R}):x\in A\}\subset 2^{B}
\end{equation} 
sei. F"ur jede Wellenfunktion gibt es die Menge der Abbildungen
\begin{equation}
\Gamma_{\Phi}:=\{\underbrace{\Phi(\mbox{id},t)}_{:=\Phi^{t}}:t\in \mathbb{R})\}\subset B^{A}\ ,
\end{equation}
exakt deren Elemente wir genau dann als die Fl"usse der Wellenfunktion $\Phi$
bezeichnen, falls $B$ in $A$ liegt.\index{Fl\"usse einer Wellenfunktion}
H"ochstens dann kann die Menge $\Gamma_{\Phi}$ das "Aquivalent einer homogenen Autobolismenmenge der Menge $A$
sein, wenn
es eine Bijektion $v$ der Menge $A$ auf die Wertemenge $\mathbf{P}_{2}\Phi\subset B$ gibt.
In diesem Fall sind n"amlich
die Mengen
$$\Gamma_{\Phi}\circ v^{-1}=\{\Phi^{t}\circ v^{-1}:t\in \mathbb{R})\}\subset \mathbf{P}_{2}\Phi^{\mathbf{P}_{2}\Phi}\ ,$$
$$v^{-1}\circ\Gamma_{\Phi}=\{v^{-1}\circ\Phi^{t}:t\in \mathbb{R})\}\subset A^{A}$$
genau dann homogene Autobolismenmengen, wenn
\begin{equation}\label{part}
[\Phi]\in\mathbf{part}(\mathbf{P}_{2}\Phi)
\end{equation} 
eine Partition ist. 
Weder die in der "Aquivalenz (\ref{wanhn}) beinhaltete 
Implikation (\ref{anhn}), noch 
die "Aquivalenz (\ref{wanhn}) des Lemma 1.3 selbst erm"oglicht den 
Schluss darauf, dass   
$$\Gamma_{\Phi}\circ v^{-1} =\langle\Gamma_{\Phi}\circ v^{-1}\rangle$$
algebraisch abgeschlossen ist, was genau dann der Fall ist, wenn
$v^{-1}\circ\Gamma_{\Phi}=\langle v^{-1}\circ\Gamma_{\Phi}\rangle$
ist; obwohl dann sowohl $[\Gamma_{\Phi}\circ v^{-1}]\in \mathbf{part}(\mathbf{P}_{2}\Phi)$
als auch $[\Gamma_{\Phi}\circ v^{-1}]\in \mathbf{part}(A)$
partitiv sind.
Hinsichtlich dieser Partitivit"at sind sie gerade so beschaffen, wie wir 
Flussfunktionen im allgemeinen anlegen.
Denn
jede Flussfunktion verstehen wir
als die folgende spezielle Form einer Wellenfunktion:  
Jede Wellenfunktion $\Phi$ ist genau dann eine Flussfunktion, wenn f"ur alle Definitionsmengenelemente
$(x,t)\in\mathbf{P}_{1}\Phi$\index{Flussfunktion} 
\begin{equation}\label{fevaf}
\begin{array}{c}
\Phi(\mbox{id},t)\in \mathbf{ab}(A)\ \land\\
\Phi(x,\mbox{id})\in\mathbf{ab}(\Phi(x,\mathbb{R}))
\end{array}
\end{equation} 
gilt und die Funktionen $\Phi(\mbox{id},t)$ und $\Phi(x,\mbox{id})$
Bijektionen sind; die Funktionen $\Phi(\mbox{id},t)$ 
sind f"ur alle $t\in \mathbb{R}$ Autobolismen
der Wertemenge $B=A$, was zu dem
Sachverhalt "aquivalent ist,
dass $[\Phi]$
eine Partition der Menge $A$ ist.\footnote{ 
Wobei uns hier exakt jede Abbildung $\Phi$ als 
eine reelle bzw. komplexe Wellenfunktion 
gelte, f"ur die es zwei 
nat"urliche Zahlen $\alpha$ und $\beta$ und zwei  
Teilmengen 
$A\subset\mathbb{R}^{\alpha}$ bzw. $A\subset\mathbb{C}^{\alpha}$ 
und 
$B\subset\mathbb{R}^{\beta}$ bzw. $B\subset\mathbb{C}^{\beta}$
gibt,
die so beschaffen sind, dass (\ref{well})
gilt. Analog bezeichnen wir 
exakt jede 
Flussfunktion, die eine reelle bzw. komplexe Wellenfunktion\index{reelle Wellenfunktion} ist, 
als  
eine reelle bzw. komplexe 
Flussfunktion.\index{reelle Flussfunktion}
\index{komplexe Flussfunktion}
Jede Wellenfunktion $\Phi\in B^{A\times \mathbb{R}}$ nennen wir eine metrische 
Wellenfunktion bez"uglich $d_{A}$ und $d_{B}$, wenn sowohl $d_{A}$ als auch $d_{B}$ Metriken 
auf $A$ bzw. $B$ sind und entsprechend nennen wir jede Flussfunktion, die eine metrische 
Wellenfunktion bez"uglich $d_{A}=d_{B}$
ist, eine metrische 
Flussfunktion bez"uglich $d_{A}$.\index{metrische Wellenfunktion bez\"uglich zweier Metriken}\index{metrische Flussfunktion bez\"uglich
einer Metrik} Die Redeweise von metrischen Wellenfunktionen bzw. von der metrischen Flussfunktionen
ist weniger streng festgelegt, als die Sprechweise von der 
metrisierten Wellenfunktion bzw. von der metrisierten Flussfunktion,
welche wir im zweiten Kapitel als bestimmte Tripel objektivieren.}
Kann es somit sein,
dass nicht
jede Flussfunktion $\Phi$ f"ur alle $t\in\mathbb{R}$ einen
Phasenfluss\index{Phasenfluss}
$\{\Psi^{t}\}_{t\in\mathbb{R}}=\{\Phi(\mbox{id},t)\}_{t\in\mathbb{R}}$ festlegt, von dem 
wir verlangen, dass die homogene Autobolismenmenge
\begin{equation}\label{fluss}
\Gamma_{\Phi}=\langle\Gamma_{\Phi}\rangle
\end{equation}
die Gruppe $(\Gamma_{\Phi},\circ)$ bildet? Das ist sehr wohl m"oglich!
Die Klasse der Flussfunktionen objektiviert lediglich den
Determinismus.\index{Determinismus} Phasizit"at und Determinismus sind aber zweierlei.
Mit Phasizit"at meinen wir dabei
das "Aquivalent algebraischer Abgeschlossenheit der Flussmengen.
Phasizit"at ist 
die Existenz der jeweiligen Phasenflussgruppe und sie stellt
eine spezielle Form des Determinismus dar, wie die folgende 
"Uberlegung transparent macht:
\newline
Weil f"ur jede Flussfunktion $\Phi$ gem"ass (\ref{fevaf}) f"ur alle Wertemengenelemente $x$
die Funktion
$\Phi(x,\mbox{id})\in\mathbf{ab}(\Phi(x,\mathbb{R}))$ invertierbar ist,
gibt es f"ur alle $t_{1},t_{2}\in\mathbb{R}$
eine Zahl $t_{\Phi}(x,t_{1},t_{2})$, die so beschaffen ist, dass
\begin{equation}\label{unz}
\Phi(\Phi(x,t_{1}),t_{2})=\Phi(x,t_{\Phi}(x,t_{1},t_{2}))
\end{equation}
ist und damit ist zu jeder Flussfunktion $\Phi$ das reellwertige
Feld
\begin{equation} 
\begin{array}{c}
t_{\Phi}:\mathbf{P}_{1}\Phi\times \mathbb{R}\to \mathbb{R}\ ,\\
(x,t_{1},t_{2})\mapsto t_{\Phi}(x,t_{1},t_{2})
\end{array}
\end{equation}
festgelegt, f"ur das f"ur alle $(x,t_{1},t_{2})\in \mathbf{P}_{1}t_{\Phi}$
\begin{equation}
\Phi^{t_{1}}\circ\Phi^{t_{1}}(x)=\Phi^{t_{\Phi}(x,t_{1},t_{2})}(x)
\end{equation}
ist: Offenbar gilt f"ur Flussfunktionen im allgemeinen nicht
das Unabh"angigkeitskriterium, dass f"ur alle Wertemengenelemente
$x\in \mathbf{P}_{2}\Phi$
\begin{equation}
t_{\Phi}=t_{\Phi}\circ (x,\mathbf{P}_{2},\mathbf{P}_{3})
\end{equation}
ist, dessen Erf"ulltheit genau dann vorliegt,
wenn die algebraische Abgeschlossenheit der Flussmenge (\ref{fluss}) gegeben ist und eine 
zur Flussfunktion $\Phi$ "aquivalente 
Phasenflussgruppe existiert.
Die 
Klasse der Phasenflussgruppen gilt 
verbreitetermassen als diejenige mathematische Objektklasse, die
als die Objektivierung des Begriffes des dynamischen Systemes aufgefasst wird;
auch das, was wir Phasenfluss nennen, deckt sich mit der etablierten 
Terminologie. Der Phasenfluss\index{Phasenfluss} $\{\Psi^{t}\}_{t\in\mathbb{R}}$ ist eine Familie
von Autobolismen des gemeinsamen Zustandsraumes $\downarrow\mathbf{P}_{1}\{\Psi^{t}:t\in\mathbb{R}\}$, die
die Phasenflussgruppe $(\{\Psi^{t}:t\in\mathbb{R}\},\circ)$ bilden. Die einzelnen
Autobolismen $\Psi^{t}$ sind dabei die Mitglieder der Familie $\{\Psi^{t}\}_{t\in\mathbb{R}}$.
\newline
Bekanntlich gibt es 
Wellenfunktionen $\Phi$,
die eine Flussmenge $\Gamma_{\Phi}$ haben, deren
Kollektivierung $[\Phi]$ partitiv ist,
denen nicht nur
Phasizit"at in dem Sinn fehlt,
dass die $\Phi$ zugeordnete 
Menge der Fl"usse 
$\Gamma_{\Phi}$ keine Gruppe bildet, sondern
denen der Determinismus der Flussfunktionen fehlt:
Der Determinismus kann nicht nur dadurch
verletzt sein, dass die Kollektivierung $[\Phi]$ der 
Flussmenge $\Gamma_{\Phi}$ nicht partitiv ist, sondern 
auch dadurch, dass es f"ur Elemente
des ersten kartesischen Faktors der  
Definitionsmenge 
$x\in \mathbf{P}_{1}\mathbf{P}_{1}\Phi$ Funktionen $\Phi(x,\mbox{id})$ gibt, die nicht invertierbar sind.
Gem"ass (\ref{fevaf}) ist f"ur 
Wellenfunktionen $\Phi$, die
Flussfunktionen sind, diese 
M"oglichkeit, den Determinismus des Phasenflusses zu verletzen, per definitionem
ausgeschlossen.\index{phasischer Fluss}\index{aphasischer Fluss}
\index{phasische Flussmenge}\index{aphasische Flussmenge}\newline
Fl"usse der Menge
$\Gamma_{\Phi}$ nennen 
wir genau dann aphasisch,
wenn $\Gamma_{\Phi}\not=\langle\Gamma_{\Phi}\rangle$ ist;
wohingegen wir exakt diejenigen
Flussmengen als phasisch bezeichnen, f"ur die 
(\ref{fluss}) gilt.
Flussfunktionen wie Wellenfunktionen $\Phi$, f"ur die $\Gamma_{\Phi}$
aphasisch bzw. phasisch ist, nennen wir ebenfalls aphasisch bzw. 
phasisch.\index{aphasische Flussfunktion}\index{phasische Flussfunktion}
Die aphasischen 
Fl"usse $\Gamma_{\Phi}$ aphasischer Flussfunktionen $\Phi$
stellen 
die angek"undigte klassische Exemplifizierung des Sachverhaltes dar, dass
die Implikation
(\ref{anhn}) nicht umkehrbar ist:
Da die Flussfunktion $\Phi$ in dem Sinn partitiv ist, dass
$[\Phi]\in\mathbf{part}(\mathbf{P}_{2}\Phi)$ 
ist, gibt es f"ur alle $\chi\in[\Phi]$ ein Element $x\in\mathbf{P}_{2}\Phi$,
f"ur das $\chi=\Phi(x,\mathbb{R})$ ist.
Wegen
der Invertierbarkeit der Funktionen $\Phi(y,\mbox{id})$ f"ur alle
$y\in \mathbf{P}_{2}\Phi$, gibt es f"ur alle nicht leeren Teilmengen $A$ und $B$ des jeweiligen 
Kollektivelementes $\chi\supset A,B$
reelle Zahlen $t$, f"ur die Koh"arenz in der Form konkretisiert ist, dass
\begin{equation}\label{konkoh}
\Phi^{t}(A)\cap B\not=\emptyset
\end{equation}
gilt.  
Und, dass 
selbst bei dem Determinismus der Flussfunktionen Aphasizit"at der generische Fall ist, ist dabei
das, was die Unumkehrbarkeit der Implikation (\ref{anhn}) dramatisiert. 
\newline
Allemal haben wir in der Menge der Fl"usse $\Gamma_{\Phi}$ einer Flussfunktion $\Psi$
eine homogene Autobolismenmenge, deren Abschluss $\langle\Gamma_{\Phi}\rangle$
existiert. Daher bezeichnen wir auch f"ur jede Flussfunktion $\Phi$ deren Wertemenge
$\mathbf{P}_{2}\Phi$ als deren Zustandsraum.
Die Elemente der Differenz $\Phi^{\vartheta}\in\langle\Gamma_{\Phi}\rangle\setminus \Gamma_{\Phi}$
k"onnen wir dabei nicht mehr mittels reeller Zahlen indizieren. Die Frage,
wie diese Differenzen beschaffen sind, f"uhrt von unserem Thema 
und der weiteren Er"orterung des Lemma 1.3 weg. 
\newline
Es gibt zwar 
f"ur jede aphasische Flussfunktion $\Phi$
die Elemente der Differenz $\Phi^{\vartheta}\in\langle\Gamma_{\Phi}\rangle\setminus \Gamma_{\Phi}$.
Da f"ur alle Kollektivelemente $\chi\in[\Phi]$ gilt, dass f"ur alle Teilmengen $A,B\subset \chi$
reelle Zahlen $t$ existieren, welche die Koh"arenz (\ref{konkoh}) konkretisieren,
gilt nach der "Aquivalenz (\ref{wanhn}) des Lemmas 1.3, dass
\begin{equation}\label{iedkoh}
\mathbf{E}(\widehat{\mathbf{T}}(\Gamma_{\Phi}))\setminus[\Phi]
=\{\emptyset\}\not=\emptyset
\end{equation}
ist. F"ur phasische Flussfunktionen $\Phi$ gilt diese Identifizierung 
der integeren Elemente der gemeinsamen invarianten Topologie 
der Autobolismenmenge $\Gamma_{\Phi}$ mit den
Elementen der Kollektivierung $[\Phi]$ 
genauso.
\chapter{Attraktoren als rein topologische Epikonstrukte}
{\small Die durchdringende Analyse des Attraktorbegriffes und seiner konzeptionellen Idee leitet uns
in diesem Kapitel zu Abstraktionen des Attraktorbegriffes. Der Schwerpunkt dieser Abstraktionen befindet sich
innerhalb der klassischen allgemeinen Topologie, die offenbar bislang, wie z.B. die
Lehrb"ucher \cite{alex}, \cite{manh}, \cite{kell} und \cite{ward},
fernab 
dieser Abstraktionen des Attraktorbegriffes liegt. \newline
Einerseits vertieft sich im Zuge dieser Begriffsausdehnung der Eindruck, dass die Lehre von 
Attraktoren wesentlich innerhalb der allgemeinen Topologie anzusiedeln sei.
Abseits des Interessenfokus dieses Kapitels, der innerhalb der 
allgemeinen Topologie zentriert ist, nehmen wir aber auch
zur Kenntnis, dass sich
der Attraktorbegriff selbst gegen eine Verallgemeinerung 
nicht streubt, die
"uber die Allgemeinheitsstufe der allgemeinen Topologie hinaus
geht: Wir k"onnen den Begriff des Attraktors 
auf der Generalit"atsstufe der Mengenlehre 
f"ur beliebige Mengensysteme verfassen.}
\section{Freie und topologische Attraktoren} 
Was ist geschehen, inwiefern sind wir nun am anderen Pol der am Eingang des Kapitels karrikierten Naseweisheit?
\newline
F"ur jede Flussfunktion $\Phi$ und f"ur jedes Mengensystem 
$\mathcal{A}\subset 2^{\mathbf{P}_{2}\Phi}$
"uber dem Zustandsraum $\mathbf{P}_{2}\Phi$, das denselben
"uberdeckt, sodass
$$\bigcup\mathcal{A}= \mathbf{P}_{2}\Phi$$
ist, k"onnen wir den freien Attraktor 
der Flussfunktion $\Phi$
relativ zu dem
Mengensystem $\mathcal{A}$ 
folgendermassen\index{freier Attraktor bez\"uglich eines Mengensystemes} definieren: 
Jede nicht leere Menge $\Theta$ der gemeinsamen invarianten Topologie 
$\widehat{\mathbf{T}}(\Gamma_{\Phi})\setminus\{\emptyset\}$
der Autobolismenmenge $\Gamma_{\Phi}$ gelte genau dann
als ein freier   
Attraktor 
der Flussfunktion $\Phi$
bez"uglich
eines Mengensystemes $\mathcal{A}$ "uber dem Zustandsraum $\mathbf{P}_{2}\Phi$, wenn
die Implikation
\begin{equation}\label{idkoh}
\begin{array}{c}
A,B\subset\Theta\ \land\ \emptyset\not\in\{A,B\}\ \land\ A,B\in\mathcal{A}\cap\Theta\\\
\Rightarrow\\
\exists\ t\in\mathbb{R}:\
\Phi^{t}(A)\cap B\not=\emptyset
\end{array}
\end{equation}
wahr ist. $\mathcal{A}\cap\Theta=\{A\cap\Theta:A\in\mathcal{A}\}$ 
ist dabei als Spurmengensystem zu lesen.
Damit generalisieren wir den entscheidenden 
Handgriff bei der herk"ommlichen Verfassung des Begriffes 
des Attraktors, das Koh"arenzkriterium zu stellen, auf triviale Weise:
Je restriktiver die durch die Wahl eines jeweiligen\index{relativierendes Mengensystem eines freien Attraktors}
relativierenden Mengensystemes $\mathcal{A}$
gestellte
Bedingung 
der Pr"amisse der Implikation (\ref{idkoh}) ist, dass
$$A,B\in(\mathcal{A}\cap\Theta)\setminus\{\emptyset\}$$
sei, desto
mehr Teilmengen $\Theta$ der gemeinsamen invarianten Topologie 
$\widehat{\mathbf{T}}(\Gamma_{\Phi})\setminus\{\emptyset\}$
gelten als freie Attraktoren. Wir generalisieren den entscheidenden
Handgriff, das Koh"arenzkriterium\index{Koh\"arenzkriterium}
ins Spiel zu bringen; diesen Handgriff zu lokalisieren, w"are bloss, denselben
zu zitieren. 
Wir wollen herausfinden, was es ist, was die spezielle Wahl ausgerechnet einer 
Topologie
zum relativierenden Mengensystem $\mathcal{A}$ so treffend 
macht: Wenn n"amlich ein Attraktor  
bez"uglich der Zustandsraumtopologie 
$\mathbf{T}(\mathbf{P}_{2}\Psi)$ im herk"ommlichen Sinn
als freier Attraktor formuliert werden soll,
w"ahlen wir die Topologie $\mathcal{A}=\mathbf{T}(\mathbf{P}_{2}\Psi)$ als relativierendes Mengensystem.
Wir sehen:
Es ist nicht so leicht, entscheidende  
Handgriffe zu erfassen.\index{entscheidender Handgriff}
\newline
Das Mengensystem all dieser freien Attraktoren, das Mengensystem
aller Mengen $\Theta$ der gemeinsamen invarianten Topologie $\widehat{\mathbf{T}}(\Gamma_{\Phi})$
der Autobolismenmenge $\Gamma_{\Phi}$ also, f"ur welche die Implikation (\ref{idkoh}) 
gilt, bezeichne 
\begin{equation}\label{ienatt}
\mbox{{\bf @}}(\Phi,\mathcal{A})\ .
\end{equation}
Kann diese Menge leer sein?
Das Lemma 1.3 sagt uns, dass immer
\begin{equation}
\mbox{{\bf @}}(\Phi,\mathcal{A})\supset\mathbf{E}(\widehat{\mathbf{T}}(\Gamma_{\Phi}))
\end{equation}
gilt, was die negative\index{dumme Frage}
Antwort auf die am Kapitelanfang gestellte \glqq dumme Frage\grqq paraphrasiert,
ob es sein kann, dass die integeren Mengen der gemeinsamen invarianten Topologien
dem Koh"arenzkriterium nicht
gen"ugen.
Allerdings ist es m"oglich, dass
integere Mengen gemeinsamer invarianter Topologien
bez"uglich einer jeweiligen Zustandsraumtopologie nicht kompakt und 
dann bez"uglich derselben keine Attraktoren sind.
Wenn $\mathcal{A}$ eine Topologie $\mathbf{T}^{\star}(\mathbf{P}_{2}\Phi)$ ist, die
den 
Zustandsraum $\mathbf{P}_{2}\Phi$ topologisiert,
ist das Mengensystem aller freien
Attraktoren 
der Flussfunktion $\Phi$ bez"uglich der Topologie $\mathbf{T}^{\star}(\mathbf{P}_{2}\Phi)$
gerade das Mengensystem, dessen  
bez"uglich der Topologie $\mathbf{T}^{\star}(\mathbf{P}_{2}\Phi)$
kompakte Elemente die
Attraktoren
des Phasenflusses $\{\Psi^{t}\}_{t\in\mathbb{R}}$ bez"uglich der 
Topologie $\mathbf{T}^{\star}(\mathbf{P}_{2}\Phi)$ im herk"ommlichen Sinn sind, sofern die Flussfunktion
$\Psi$ phasisch ist:\newline
F"ur jede Topologie $\mathbf{T}$
bezeichne $\mathbf{C}(\mathbf{T})$ die Menge aller Kompakta des topologischen 
Raumes $(\bigcup\mathbf{T},\mathbf{T})$.
Der verallgemeinerte Begriff des freien Attraktors als Element des
jeweiligen Mengensystemes $\mbox{{\bf @}}(\Phi,\mathcal{A})$ umfasst noch nicht
den herk"ommlichen Begriff des Attraktors 
eines Phasenflusses $\{\Psi^{t}\}_{t\in\mathbb{R}}$
bez"uglich einer Topologisierung 
dessen Zustandsraumes. 
Aber wir haben schliesslich in dem Mengensystem jeweiliger Spurmengensysteme
$$\{\mathbf{C}(\mathbf{T}^{\star}(\mathbf{P}_{2}\Phi))\cap A:
A\in\mbox{{\bf @}}(\Phi,\mathcal{A})\}$$
$$=\{\mbox{{\bf @}}(\Phi,\mathcal{A})\}\cap C:
C\in \mathbf{C}(\mathbf{T}^{\star}(\mathbf{P}_{2}\Phi))\}$$ 
\begin{equation}
=:\mbox{{\bf @}}(\Phi,\mathcal{A})\ \cap\ \mathbf{C}(\mathbf{T}^{\star}(\mathbf{P}_{2}\Phi))=:
\mbox{{\bf @}}_{\mathbf{T}^{\star}(\mathbf{P}_{2}\Phi)}(\Phi,\mathcal{A})
\end{equation}
gerade das 
Mengensystem all derjeniger Mengen, die wir im 
herk"ommlichen Sinn
als 
die Attraktoren
des Phasenflusses $\{\Psi^{t}\}_{t\in\mathbb{R}}$ bez"uglich der 
Topologie $\mathbf{T}^{\star}(\mathbf{P}_{2}\Phi)$ bezeichnen, falls das Mengensystem $\mathcal{A}$
zu der Topologie
$$\mathcal{A}=\mathbf{T}^{\star}(\mathbf{P}_{2}\Phi)$$
spezifiziert ist -- sofern die Flussfunktion $\Phi$ 
phasisch ist. Denn das, was wir exakt
als den herk"ommlichen Begriff des Attraktors auffassen, ist
f"ur Phasenfl"usse und nicht f"ur Flussfunktionen festgelegt.
Inwiefern sind wir aber nun eigentlich am 
anderen Pol der Naseweisheit, die 
kurzhand meint, dass Attraktoren lediglich die integeren Mengen des Mengensystemes $\mathbf{E}(\widehat{\mathbf{T}}(\Gamma_{\Phi}))$ jeweiliger 
gemeinsamer invarianter Topologien $\widehat{\mathbf{T}}(\Gamma_{\Phi})$
seien? Das Mengensystem $\mathbf{E}(\widehat{\mathbf{T}}(\Gamma_{\Phi}))$ der integeren Mengen der
gemeinsamen invarianter Topologien $\widehat{\mathbf{T}}(\Gamma_{\Phi})$
ist
der Schnitt
\begin{equation} 
\bigcap\Bigl\{\mbox{{\bf @}}(\Phi,\mathcal{A}):\subset 2^{\mathbf{P}_{2}\Phi}\land\bigcup\mathcal{A}=\mathbf{P}_{2}\Phi\Bigr\}=
\mathbf{E}(\widehat{\mathbf{T}}(\Gamma_{\Phi}))\setminus\{\emptyset\} \ ,
\end{equation}
gleichsam die Objektivierung des kleinsten Begriffes des freien Attraktors.\newline\newline    
Auf der Klasse der Paare $(\Phi,\mathcal{A})$, deren jeweilige erste Komponente eine Flussfunktion
$\Phi$
ist 
und deren jeweilige zweite Komponente eine "uberdeckende Teilmenge $\mathcal{A}$ der 
Potenzmenge $2^{\mathbf{P}_{2}\Phi}$ des Zustandsraumes 
der jeweiligen ersten Komponente $\Phi$ ist,
k"onnen wir den Operator $\mbox{{\bf @}}$
definieren, dessen jeweilige Werte
die Mengensysteme freier Attraktoren $\mbox{{\bf @}}(\Phi,\mathcal{A})$ 
seien. Der Operator $\mbox{{\bf @}}_{\mathbf{P}_{3}}$
hingegen sei auf der Klasse der Tripel $(\Phi,\mathcal{A},\mathbf{T})$
definiert, deren erste beiden Komponenten wie bei den 
Paaren der Klasse beschaffen sind, auf welcher der Operator $\mbox{{\bf @}}$
definiert ist. Die dritte Komponente der Klasse der Tripel $(\Phi,\mathcal{A},\mathbf{T})$,
auf welcher der Operator $\mbox{{\bf @}}_{\mathbf{P}_{3}}$ erkl"art werden soll,
ist eine auf dem Zustandsraum
der ersten Komponente $\Phi$ gegebene Topologie.
Die Werte des Operators $\mbox{{\bf @}}_{\mathbf{P}_{3}}$
sind die 
Mengensysteme 
bez"uglich der Zustandsraumtopologie $\mathbf{T}$
kompakter Attraktoren $\mbox{{\bf @}}_{\mathbf{T}}(\Phi,\mathcal{A})$. 
Mit der Einf"uhrung 
der 
Operatoren $\mbox{{\bf @}}$ und $\mbox{{\bf @}}_{\mathbf{P}_{3}}$  
kamen wir also unserer am Eingang des Kapitels 
gestellten
Aufgabe nach, zu erkunden, unter welchen Umst"anden 
verallgemeinernde 
Versionen des
Urbegriffes des Attraktors\index{Urbegriff des Attraktors}
so beschaffen sind, dass sie den 
in Gestalt der integeren flussinvarianten Mengen
pr"azisiert formulierten Begriff kleinster
flussinvarianter Mengen transzendieren. F"ur die 
Operatoren $\mbox{{\bf @}}$ und $\mbox{{\bf @}}_{\mathbf{P}_{3}}$ gilt
\begin{equation} 
\mbox{{\bf @}}\supset\mbox{{\bf @}}_{\mathbf{P}_{3}}
\end{equation}
insofern, als f"ur alle Paare $(\Phi,\mathcal{A})\in\mathbf{P}_{1}\mbox{{\bf @}}$ und jede
Zustandsraumtopologie $\mathbf{T}$
$$\mbox{{\bf @}}(\Phi,\mathcal{A})\supset\mbox{{\bf @}}_{\mathbf{T}}(\Phi,\mathcal{A})$$
gilt.
Ferner k"onnen wir die aus dem Lemma 1.3 gefolgerte Identifizierung (\ref{iedkoh}) 
der integeren Elemente der gemeinsamen invarianten Topologie 
der Autobolismenmenge $\Gamma_{\Phi}$ mit den
Kollektivelementen aus $[\Phi]$ fortsetzen: F"ur alle Flussfunktionen $\Phi$ ist
\begin{equation} 
\mathbf{E}(\widehat{\mathbf{T}}(\Gamma_{\Phi}))\setminus\{\emptyset\} =[\Phi]=\mbox{{\bf @}}(\Phi,2^{\mathbf{P}_{2}\Phi})\ ,
\end{equation}
wobei
f"ur alle $\mathbf{P}_{2}\Phi$ "uberdeckenden Mengensysteme $\mathcal{A}\subset 2^{\mathbf{P}_{2}\Phi}$ die Implikationen
\begin{equation} 
\begin{array}{c}
\forall\ A\in\mathcal{A}\ \exists\ \gamma\in\lbrack\Phi\rbrack:A\supset\gamma\\
\Rightarrow\\
\mbox{{\bf @}}(\Phi,\mathcal{A})=\lbrack\Phi\rbrack\ ,\\
\quad\\
\mathcal{A}\supset\{Q\subset\mathbf{P}_{2}\Phi:\exists\ \gamma\in\lbrack\Phi\rbrack:\gamma\supset Q\}\\
\Rightarrow\\
\mbox{{\bf @}}(\Phi,\mathcal{A})=\lbrack\Phi\rbrack
\end{array}
\end{equation}
gelten. 
Seit der Einf"uhrung 
der 
Operatoren $\mbox{{\bf @}}$ und $\mbox{{\bf @}}_{\mathbf{P}_{3}}$ 
gehen wir der Aufgabe nach,
das Verh"altnis 
verallgemeinernder 
Versionen des
Urbegriffes des Attraktors\index{Urbegriff des Attraktors}
zu den integeren Mengen jeweiliger gemeinsamer invarianter Topologien zu erkunden.
Wir verfolgen unsere Aufgabenstellung noch ein wenig weiter
und verfassen
verallgemeinernde Versionen des
Urbegriffes des Attraktors als rein topologische Epikonstrukte.
Mit rein topologischen Epikonstrukten meinen wir die
Resultate von Konstruktionen, welche die
Werte solcher
universeller Operatoren
sind, die auf homogenen Klassen von 
Mengen von Topologien oder von
Tupeln von Topologien definiert sind. 
\newline
\newline
{\bf Definition 2.1: Topologischer Attraktor}\newline
{\em Es sei $\mathbb{T}$ eine Menge von Topologien und 
${\rm X}^{\star},{\rm X}_{\star}$ seien zwei Mengensysteme, welche
die Menge
\begin{equation}\label{aegy}
Z(\mathbb{T}):=\bigcup\bigcap\mathbb{T}
\end{equation}
"uberdecken und die in deren Potenzmenge liegen.
Exakt die Mengen des Mengensystemes
\begin{equation}\label{dikos} 
\begin{array}{c}
\Bigl\{{\rm K}\in (2^{Z(\mathbb{T})}\cap{\rm X}_{\star})\setminus\{\emptyset\} :\\
A,B\subset{\rm K}\ \land\ A,B\in {\rm X}^{\star}\Rightarrow\exists\ 
\theta\in\bigcap\mathbb{T}:\\
\vartheta\cap A\not=\emptyset\not=B\cap\vartheta\Bigl\}\\
=:\underline{\mbox{{\bf @}}}(\mathbb{T},{\rm X}^{\star},{\rm X}_{\star})
\end{array}
\end{equation}
bezeichnen wir als die topologischen Attraktoren des 
Tripels\index{topologischer Attraktor} 
\begin{equation}\label{diko}
(\mathbb{T},{\rm X}^{\star},{\rm X}_{\star}) 
\end{equation}
oder
als die
topologischen Attraktoren der Topologienmenge $\mathbb{T}$ mit
der Koh"arenz ${\rm X}^{\star}$ und mit der
Kadenz ${\rm X}_{\star}$.\index{Kadenz eines
topologischen Attraktors}\index{Koh\"arenz eines topologischen Attraktors}}
\newline
\newline
Dieses Mengensystem topologischer Attraktoren der Topologienmenge $\mathbb{T}$ mit
der Koh"arenz ${\rm X}^{\star}$ und mit der
Kadenz ${\rm X}_{\star}$ ist wohldefiniert.
Der Schnitt
$\bigcap\mathbb{T}$
existiert allemal und wenn er leer ist, ist 
$\underline{\mbox{{\bf @}}}(\mathbb{T},{\rm X}^{\star},{\rm X}_{\star})=\emptyset$ und es gibt dann
keinen topologischen Attraktor der Topologienmenge $\mathbb{T}$, gleich, welche
Koh"arenz ${\rm X}^{\star}$ und gleich, welche
Kadenz ${\rm X}_{\star}$ vorliegt. Allemal gilt aber 
auch
\begin{equation}
\{\emptyset\}\not\in\underline{\mbox{{\bf @}}}(\mathbb{T},{\rm X}^{\star},{\rm X}_{\star})\ .
\end{equation}
Die leere Menge ist kein topologischer Attraktor --
jedweder Topologienmenge $\mathbb{T}$, gleich, welche
Koh"arenz ${\rm X}^{\star}$ und gleich, welche
Kadenz ${\rm X}_{\star}$ gegeben ist.
Die leere Menge kann auch weder ein freier Attraktor noch ein kompakter Attraktor
einer Flussfunktion
bez"uglich irgendeines Mengensystemes und einer Topologie sein; und 
damit ist die leere Menge auch kein Attraktor im herk"ommlichen Sinn.
Wir merken uns: Die leere Menge ist nicht attraktiv.
\newline 
Der Fall, dass die 
Kadenz ${\rm X}_{\star}$ mit der Potenzmenge $2^{Z(\mathbb{T})}$ "ubereinstimmt,
sei durch die Vereinbarung, dass  
\begin{equation}\label{dikoha}
\underline{\mbox{{\bf @}}}(\mathbb{T},{\rm X}^{\star}):=
\underline{\mbox{{\bf @}}}(\mathbb{T},{\rm X}^{\star},2^{Z(\mathbb{T})})
\end{equation}
sein soll, notativ vereinfacht. 
Exakt die 
Elemente des Mengensystemes $\underline{\mbox{{\bf @}}}(\mathbb{T},{\rm X}^{\star})$ nennen
wir die freien Attraktoren der Topologienmenge $\mathbb{T}$ mit
der Koh"arenz ${\rm X}^{\star}$.\index{freier topologischer Attraktor}
Weitergehend sei in dem analogen Fall, dass die nicht effektive 
Koh"arenz 
$\mathbb{T}^{\star}=2^{Z(\mathbb{T})}$
vorliegt, schlicht
\begin{equation}\label{dikohb}
\underline{\mbox{{\bf @}}}(\mathbb{T}):=\underline{\mbox{{\bf @}}}(\mathbb{T},
2^{Z(\mathbb{T})})
\end{equation}
gesetzt. Exakt die Elemente dieses Mengensystemes nennen wir die inkoh"arenten topologischen
Attraktoren der
Topologienmenge $\mathbb{T}$.\index{inkoh\"arenter topologischer Attraktor}
Es kann dabei sein, dass ein und dieselbe Menge ein inkoh"arenter topologischer
Attraktor der
Topologienmenge $\mathbb{T}_{1}$ ist und gleichzeitig ein topologischer
Attraktor der 
Topologienmenge
$\mathbb{T}_{2}$ einer Koh"arenz ${\rm X}^{\star}$.
In dem Fall der
Spezifizierung 
\begin{equation}\label{dikoh}
\mathbb{T}=\mathbf{T}(\Gamma_{\Phi})=\{\mathbf{T}(\Phi^{t}):\Phi^{t}\in\Gamma_{\Phi}\}
\end{equation}
stimmt demnach
das Mengensystem topologischer Attraktoren\index{topologischer Attraktor}
$$\underline{\mbox{{\bf @}}}(\mathbf{T}(\Gamma_{\Phi}))=\underline{\mbox{{\bf @}}}(\mathbf{T}(\Gamma_{\Phi}),2^{\mathbf{P}_{2}\Phi})={\mbox{{\bf @}}}(\Phi,2^{\mathbf{P}_{2}\Phi})$$
$$=[\Phi]=\mathbf{E}(\widehat{\mathbf{T}}(\Gamma_{\Phi}))\setminus\{\emptyset\}$$
mit dem Mengensystem der bez"uglich der gemeinsamen invarianten 
Topologie $\mathbf{T}(\Gamma_{\Phi})$ integeren
Mengen "uberein, was
die aus dem Lemma 1.3 gefolgerte Identifizierung (\ref{iedkoh}) weiter fortsetzt.\newline
Hat denn die 
mit der Festlegung der jeweiligen Werte gem"ass (\ref{dikos}) bzw. (\ref{dikoha}) bzw. (\ref{dikohb})
bestimmte
Konstruktion der homonymen Operatoren $\underline{\mbox{{\bf @}}}$, die auf der Klasse der Tripel der Form (\ref{diko}) bzw. auf
der Klasse der Paare von Topologienmengen und von Topologien bzw. auf
der Klasse von Mengen von Topologien definiert sind,\index{topologischer Attraktor}
"uberhaupt mit Attraktoren zu tun?\newline
Ja, denn
der Begriff des topologischen Attraktors
generalisiert den herk"ommlichen Begriff des Attraktors: F"ur
jede Flussfunktion $\Phi$
ist f"ur die Spezifizierung $\mathbb{T}=\mathbf{T}(\Gamma_{\Phi})$ gem"ass (\ref{dikoh})
f"ur jede auf dem Zustandsraum der jeweiligen Flussfunktion $\Phi$
gegebene Topologie $\mathbf{T}(\mathbf{P}_{2}\Phi)$
das Mengensystem 
\begin{equation}
\underline{\mbox{{\bf @}}}(\mathbf{T}(\Gamma_{\Phi}),\mathbf{T}(\mathbf{P}_{2}\Phi),\mathbf{C}(\mathbf{T}(\mathbf{P}_{2}\Phi)))
\end{equation}
die Menge aller Attraktoren der Flussfunktion $\Phi$ bez"uglich der
Topologie $\mathbf{T}(\mathbf{P}_{2}\Phi)$, ganz im herk"ommlichen Sinn. $\mathbf{C}(\mathbf{T})$ ist hierbei, wie am Eingang dieses Abschnittes bereits 
festgelegt, die Menge aller Kompakta des topologischen Raumes 
$(\bigcup\mathbf{T},\mathbf{T})$.
\newline
Daher m"ussen wir auch die 
im Traktat "uber den elementaren Quasiergodensatz \cite{raab}
als Zimmer bezeichneten Attratoren endlichdimensionaler 
reeller oder komplexer R"aume im Wertebreich der
Operatoren $\underline{\mbox{{\bf @}}}$ wiederfinden k"onnen. Und in der Tat:  
Es sei $\Lambda$ eine 
stetige reelle Flussfunktion, d.h eine bez"uglich der nat"urlichen
Produkttopologie der Definitionsmenge $\mathbf{P}_{1}\Lambda=\mathbf{P}_{2}\Lambda\times \mathbb{R}$ und der
nat"urlichen Topologie $\mathbf{T}(\dim\Lambda)$ stetige Funktion;\index{Produkttopologie}
deren 
bez"uglich der nat"urlichen Topologie $\mathbf{T}(\dim\Lambda)$ des $\mathbb{R}^{\dim\Lambda}$
kompakter Zustandsraum $\mathbf{P}_{2}\Lambda$ habe lauter 
Zust"ande der Komponentenzahl $\dim\Lambda\in\mathbb{N}$.
Es ist
$$[\Lambda]=\mbox{{\bf @}}(\Lambda,\widehat{\mathbf{T}}(\Gamma_{\Lambda}))$$
$$\subset\mbox{{\bf @}}(\Lambda,\mathbf{T}(\dim\Lambda))=
\underline{\mbox{{\bf @}}}(\Lambda,\widehat{\mathbf{T}}(\Gamma_{\Lambda}),\mathbf{C}(\mathbf{T}(\dim\Lambda)))$$
$$=[\Lambda]\cup[[\Lambda]]$$
im allgemeinen von dem Mengensystem
$$\underline{\mbox{{\bf @}}}(\mathbf{T}(\Gamma_{\Lambda}),2^{\mathbf{P}_{2}\Lambda},\mathbf{C}(\mathbf{T}(\dim\Lambda)))
=[\Lambda]\cap\mathbf{C}(\mathbf{T}(\dim\Lambda))$$
verschieden,
weil die gemeinsame invariante Topologie $\widehat{\mathbf{T}}(\Gamma_{\Lambda}))=[\Lambda]\cup\{\emptyset\}$ der 
Flussmenge $\Gamma_{\Lambda}$ im allgemeinen eklatanterweise nicht lauter
bez"uglich der nat"urlichen Topologie $\mathbf{T}(\dim\Lambda)$
kompakte Mengen hat: Die freien topologischen Attraktoren der 
gemeinsamen invarianten Topologie $\widehat{\mathbf{T}}(\Gamma_{\Lambda})$ 
der Menge 
$$\mbox{{\bf @}}(\Lambda,\mathbf{T}(\dim\Lambda))=
\underline{\mbox{{\bf @}}}(\Lambda,\widehat{\mathbf{T}}(\Gamma_{\Lambda}),\mathbf{C}(\mathbf{T}(\dim\Lambda)))=[[\Lambda]]$$
sind mit Ausnahme der leeren 
Menge die Abschl"usse 
der integeren Mengen der invarianten Topologie $\widehat{\mathbf{T}}(\Gamma_{\Lambda})$ 
bez"uglich der nat"urlichen Topologie $\mathbf{T}(\dim\Lambda)$.
Dies erkennen wir aber nicht etwa soeben. Wir wissen dies vom Traktat "uber den elementaren Quasiergodensatz \cite{raab} 
her,
in dem uns der Satz von der Existenz 
der Zimmer $(2.1.2)^{0}$ sagt, dass die Menge topologischer Attraktoren\index{topologischer Attraktor}
$\underline{\mbox{{\bf @}}}(\mathbf{T}(\Gamma_{\Lambda}),\mathbf{T}(\dim\Lambda),\mathbf{T}(\dim\Lambda))$
die Menge der Zimmer der Flussfunktion $\Lambda$, die wir als
$$[[\Lambda]]=\{\mathbf{cl}(\gamma):\gamma\in[\Lambda]\}$$
bezeichnen, umfasst. Die Kollektivelemente $\gamma\in[\Lambda]$, die keine
Zimmer sind, sind bei der betrachteten Flussfunktion $\Lambda$ nicht kompakt und daher nicht
in der Menge topologischer Attraktoren
in $\underline{\mbox{{\bf @}}}(\mathbf{T}(\Gamma_{\Lambda}),\mathbf{T}(\dim\Lambda),\mathbf{C}(\mathbf{T}(\dim\Lambda)))$.
\newline
Nicht erst, wenn der Zustandsraum $\mathbf{P}_{2}\Lambda$ 
bez"uglich der nat"urlichen Topologie $\mathbf{T}(\dim\Lambda)$
offene 
Teilmengen enth"alt, kommt dabei die Stetigkeit der Flussfunktion $\Lambda$ zum tragen:
Die Stetigkeit der Flussfunktion $\Lambda$ impliziert,
dass die Menge vereinzelter Zust"ande 
$$Z_{\Lambda}\subset \{z:\{z\}\in [\Lambda]\}$$
in der Menge der Fixelemente aller Fl"usse aus $\Gamma_{\Lambda}$ 
enthalten ist.
Wir erl"autern, dass
nicht nur
ein im Zustandsraum isolierter Punkt als ein vereinzelter Punkt gilt. 
Ein im Zustandsraum isolierter Punkt liegt lediglich in keiner Umgebung, die
im Zustandsraum enthalten ist.\index{isolierter Punkt}\index{vereinzelter Punkt}  
W"ahrend jeder im Zustandsraum vereinzelte Punkt $z\in Z_{\Lambda}$ ein Punkt ist, 
f"ur den es keine zusammenh"angende 
Umgebung der Relativtopologie $\mathbf{T}(\dim\Lambda)\cap\mathbf{P}_{2}\Lambda$ gibt, deren Element er ist.
\newline
Wie aber ist es um die Partitivit"at des Zustandsraumes in Trajektorienh"ullen in dem allgemeinen Fall bestellt, dass eine
Flussfunktion $\tilde{\Lambda}$ stetig ist bez"uglich derjeniger 
Produkttopologie der Definitionsmenge $\mathbf{P}_{1}\tilde{\Lambda}=\mathbf{P}_{2}\tilde{\Lambda}\times \mathbb{R}$, die
sich f"ur die Zustandsraumtopologie $\mathbf{T}(\mathbf{P}_{2}\tilde{\Lambda})$ und
die nat"urliche Zahlenstrahltopologie $\mathbf{T}(1)$ 
ergibt, und bez"uglich der Zustandsraumtopologie $\mathbf{T}(\mathbf{P}_{2}\tilde{\Lambda})$?
Auskunft dar"uber gibt uns der insensitive topologische Ergodensatz.
\section{Topologische Ergodizit"at als Konsequenz\\ der Abschlusskommutation}
Sei $(X,\mathbf{T})$ ein topologischer Raum und $f$ ein Hom"oomorphismus, der den topologischen Raum $(X,\mathbf{T})$
auf sich selbst abbildet;
$\mathbf{cl}_{\mathbf{T}}$ sei der gem"ass (\ref{jedwed}) verfasste topologische H"ullenoperator.\index{topologischer H\"ullenoperator}
Die allgemein vertraute Aussage, dass dann f"ur alle Teilmengen $Y\subset X$ die Identit"at
\begin{equation}
f(\mathbf{cl}_{\mathbf{T}}(Y))=\mathbf{cl}_{\mathbf{T}}(f(Y))
\end{equation}
gilt, bezeichnen wir als die Abschlusskommutation.\index{Abschlusskommutation}
Die Abschlusskommutation k"onnen wir als Spezifizierung des Explizierungssatzes 3.8\index{Explizierungssatz} auffassen.\footnote{Somit findet sich in dieser Abhandlung
auch der Beweis der Abschlusskommutation.} 
\newline
\newline
{\bf Insensitiver topologischer Ergodensatz 2.2:}\index{insensitiver topologischer Ergodensatz}\newline
{\em Es sei $\Phi$ eine Flussfunktion, all deren 
Fl"usse $\Phi^{t}\in \Gamma_{\Phi}$ stetig seien bez"uglich der Topologie des Zustandsraumes
$\mathbf{T}(\mathbf{P}_{2}\Phi)$, welchen eine Teilmenge 
\begin{equation}
{\rm P}\in 2^{\widehat{\mathbf{T}}(\Gamma_{\Phi})}\cap\mathbf{part}(\mathbf{P}_{2}\Phi)
\end{equation}
der gemeinsamen invarianten Topologie
partioniere. Dann gilt} 
\begin{equation}
\mathbf{cl}_{\mathbf{T}(\mathbf{P}_{2}\Phi)}({\rm P})\in 2^{\widehat{\mathbf{T}}(\Gamma_{\Phi})}\cap\mathbf{part}(\mathbf{P}_{2}\Phi)\ .
\end{equation}
\newline
Wobei erl"autert sei, dass
$$\mathbf{cl}_{\mathbf{T}(\mathbf{P}_{2}\Phi)}({\rm P})=\{\mathbf{cl}_{\mathbf{T}(\mathbf{P}_{2}\Phi)}(p):
p\in {\rm P}\}$$
ist.
\newline
\newline
{\bf Beweis:}\newline
Zu jedem der bez"uglich der Zustandstopologie
$\mathbf{T}(\mathbf{P}_{2}\Phi)$ als stetig vorausgesetzten 
Fl"usse $\Phi^{t}\in \Gamma_{\Phi}$ existiert dessen jeweilige Inversion $(\Phi^{t})^{-1}\in \Gamma_{\Phi}$,
die ebenfalls bez"uglich der Zustandstopologie stetig ist. Jeder 
Fluss $\Phi^{t}\in \Gamma_{\Phi}$ ist also ein Hom"oomorphismus des Zustandsraumes 
bez"uglich dessen Topologisierung $\mathbf{T}(\mathbf{P}_{2}\Phi)$, der den Zustandsraum $\mathbf{P}_{2}\Phi$ 
auf sich selbst abbildet.\newline
F"ur alle Fl"usse $\Phi^{t}\in \Gamma_{\Phi}$ und f"ur alle Partitionselemente $p\in {\rm P}$ gilt
$$\Phi^{t}(p)=p\ ,$$
weil die Partition ${\rm P}$ aus lauter gegen"uber allen Fl"ussen $\Phi^{t}\in \Gamma_{\Phi}$ invarianten
Mengen der gemeinsamen invarianten Topologie $\mathbf{T}(\Gamma_{\Phi})$ besteht.
Wegen der Stetigkeit aller Fl"usse $\Phi^{t}\in \Gamma_{\Phi}$ bez"uglich der Topologie des Zustandsraumes
$\mathbf{T}(\mathbf{P}_{2}\Phi)$ ist f"ur alle Fl"usse $\Phi^{t}\in \Gamma_{\Phi}$ und 
alle Partitionselemente $p\in {\rm P}$ wahr, dass 
$$\Phi^{t}(\mathbf{cl}_{\mathbf{T}(\mathbf{P}_{2}\Phi)}(p))=
\mathbf{cl}_{\mathbf{T}(\mathbf{P}_{2}\Phi)}(\Phi^{t}(p))$$
ist, wobei
wegen der Invarianz der Menge $p$
$$\mathbf{cl}_{\mathbf{T}(\mathbf{P}_{2}\Phi)}(\Phi^{t}(p))=\mathbf{cl}_{\mathbf{T}(\mathbf{P}_{2}\Phi)}(p)$$
ist: Alle Abschl"usse
$$\mathbf{cl}_{\mathbf{T}(\mathbf{P}_{2}\Phi)}(p)\in\mathbf{cl}_{\mathbf{T}(\mathbf{P}_{2}\Phi)}({\rm P})$$ 
sind demnach invariant gegen"uber 
allen Fl"ussen $\Phi^{t}\in \Gamma_{\Phi}$, sodass
das Mengensystem der Abschl"usse
$$\mathbf{cl}_{\mathbf{T}(\mathbf{P}_{2}\Phi)}({\rm P})\subset \mathbf{T}(\Gamma_{\Phi})$$
eine Teilmenge der gemeinsamen invarianten Topologie $\mathbf{T}(\Gamma_{\Phi})$ ist.
Jede Teilmenge der gemeinsamen invarianten Topologie $\mathbf{T}(\Gamma_{\Phi})$ hat lauter
paarweise disjunkte Elemente. Da 
das Mengensystem
${\rm P}$ den Zustandsraum $\mathbf{P}_{2}\Phi$ partioniert,
"uberdeckt es denselben auch und daher "uberdeckt das Mengensystem der 
Abschl"usse $\mathbf{cl}_{\mathbf{T}(\mathbf{P}_{2}\Phi)}({\rm P})$ den Zustandsraum
erst recht. Letzteres ist also ein Partition desselben, die eine
Teilmenge der gemeinsamen invarianten Topologie $\mathbf{T}(\Gamma_{\Phi})$ ist.\newline
{\bf q.e.d.}
\newline
\newline
Mit ein wenig K"uhnheit k"onnen wir also 
den insensitiven topologischen Ergodensatz in der Form paraphrasieren, 
zu behaupten, 
dass
f"ur jede eine Flussfunktion $\Phi$, all deren 
Fl"usse $\Phi^{t}\in \Gamma_{\Phi}$ stetig bez"uglich der Topologie des Zustandsraumes
$\mathbf{T}(\mathbf{P}_{2}\Phi)$ sind
\begin{equation}\label{cloe}
\mathbf{cl}_{\mathbf{T}(\mathbf{P}_{2}\Phi)}\Bigl( 
2^{\widehat{\mathbf{T}}(\Gamma_{\Phi})}\cap\mathbf{part}(\mathbf{P}_{2}\Phi)\Bigr)=
2^{\widehat{\mathbf{T}}(\Gamma_{\Phi})}\cap\mathbf{part}(\mathbf{P}_{2}\Phi)
\end{equation}
ist. 
Diese Identit"at identifizieren wir
mit dem insensitiven topologischen Ergodensatz.\index{insensitiver topologischer Ergodensatz} 
Die in dieser Identit"at
formulierte, deutlich allgemeinere und umfassendere Aussage als der spezifizierte Satz von der 
Existenz\index{Satz von der Existenz der Zimmer}
der Zimmer $(2.1.2)^{0}$ des Traktates \cite{raab} "uber den elementaren Quasiergodensatz 
ist in der speziellen Version, dass  
\begin{equation}\label{thisc}
\begin{array}{c}
\mathbf{T}(\mathbf{P}_{2}\Phi)=\mathbf{T}(n)\ \land\\
\Phi^{\mathbb{R}}\subset\mathcal{C}(\mathbf{T}(n))
\end{array}
\end{equation}
gilt,
also bereits bekannt;\index{Abschlusskommutation}
wobei $\mathcal{C}(\mathbf{T})$ f"ur alle Topologien $\mathbf{T}$ die Menge aller bez"uglich $\mathbf{T}$
stetiger Abbildungen der Menge $\bigcup\mathbf{T}$ auf dieselbe und wobei $\mathbf{T}(n)$ die 
jeweilige nat"urliche Zustandsraumtopologie
bezeichne: Denn wir zeigten in jener Abhandlung des elementaren Quasiergodensatzes, dass f"ur die Flussfunktionen 
$\Phi$, welche die Spezifizierung (\ref{thisc}) erf"ullen,
den Voraussetzungen des Satzes von der 
Existenz
der Zimmer gen"ugen, falls der Zustandsraum $\mathbf{P}_{2}\Phi$ kompakt ist.
Dabei illustrieren die Betrachtungen
jenes Traktates \cite{raab}, in welchem Sinn der insensitive topologische Ergodensatz 
als ein abstrahierter insensitiver Ergodensatz anzusehen ist.
\newline
Wer nun also meinte, dass der 
Satz von der Existenz der Zimmer reduzierbar sei auf eine spezielle
Variante des recht trivialen insensitiven topologischen Ergodensatzes, 
der m"usste sich wundern. Er m"usste sich dar"uber wundern, dass die
Abhandlung des elementaren Quasiergodensatzes die Argumentation pr"asentiert, die zum 
Satz von der Existenz der Zimmer leitet, welche im Vergleich zu dem 
gerade vorgetragenen Beweis des insensitiven topologischen Ergodensatzes
aufw"andiger ist.
Wer dies meint, verkennt die Differenz zwischen relativer 
Gegenstandshaltigkeit und Nicht-Gegenstandshaltigkeit, die
zwischen dem insensitiven topologischen Ergodensatz und dem Satz von der 
Existenz der Zimmer darin besteht, 
dass der insensitive topologische Ergodensatz
keine Partivit"atsaussage "uber chaotische Dynamiken macht, der
Satz von der Existenz der Zimmer indess sehr wohl:
\newline
Der Satz von der 
Existenz
der Zimmer $(2.1.2)^{0}$ der Abhandlung des elementaren Quasiergodensatzes
hat eine andere Voraussetzungstruktur als der insensitive topologische Ergodensatz.
Diesen Unterschied in der Voraussetzungstruktur
werden wir im zweiten Teil der Konzepte der abstrakten Ergodentheorie
als ein wesentlichen Unterschied
herausarbeiten: Bei dem Satz von der 
Existenz
der Zimmer ist  
die Immanenz\index{Immanenz} und  
die Komanenz\index{Komanenz} von viablen Flussfunktionen\index{viable Flussfunktion} vorausgesetzt.
Die Immanenz einer viablen Flussfunktion formuliert dabei insofern eine Voraussetzung, welche
indirekt
die Ausdehnung des Zustandsraumes betrifft, als die Immanenz
einer stetigen viablen Flussfunktion aus der Kompaktheit 
des endlichdimensionalen reellen Zustandsraumes folgt.
Bei dem insensitiven topologischen Ergodensatz sind keine 
Voraussetzungen "uber den Zustandsraum n"otig. Darin ist
der insensitive topologische Ergodensatz in seiner Reduktion auf die
nat"urliche Zustandsraumtopologie gem"ass (\ref{thisc}) umfassender als der
der Satz $(2.1.2)^{0}$ in der 
spezifizierten Form, dass (\ref{thisc}) f"ur einen kompakten
Zustandsraum $\mathbf{P}_{2}\Phi$
erf"ullt ist.\newline Die Komanenz allerdings 
ist eine Stetigkeitsforderung an eine 
Flussfunktion $\Phi$, die weniger fordert als dies, dass alle ihre 
Fl"usse $\Phi^{t}\in \Phi^{\mathbb{R}}$ Elemente der Menge 
$\mathcal{C}(\mathbf{T}(n))$ sind, wenn $\mathbf{T}(n)$ die
jeweilige nat"urliche Zustandsraumtopologie ist. Und zwar 
verlangt die Komanenz in folgender Hinsicht wesentlich weniger als 
die Inklusion 
$\Phi^{\mathbb{R}}\subset \mathcal{C}(\mathbf{T}(n))$: Diese 
Inklusion schliesst aus, dass sich f"ur $\Phi$ diejenige 
Form der Sensitivit"at\index{Sensitivit\"at} zeigen kann, 
die wir im zweiten 
Teil der Konzepte der abstrakten Ergodentheorie
als die undramatischste Form der Sensitivit"at bestimmen werden.
N"amlich als die
mengenweise 
Sensitivit"at von nulltem Grad, die wir im zweiten 
Teil neben anderen Formen der 
Sensitivit"at vorstellen werden. Die Formen der innerhalb von 
Attraktoren gegebenen
Sensitivit"at bestimmen die Formen des Chaos.\index{Chaos}
Die undramatischste Form des Chaos ist die durch die
mengenweise 
Sensitivit"at von nulltem Grad bestimmte Form des Chaos.
Daher behandelt der insensitive topologische Ergodensatz, wie sich im zweiten Teil
herausstellen wird, 
nur Szenarien
ohne Chaos.\newline Ganz anders dagegen ist es mit der
Komanenz. Die Komanenz schliesst keine Form der Sensitivit"at aus.
Der Satz von der Existenz der Zimmer ist zwar im vertrauten Rahmen der
endlichdimensionalen reellen Zustandsr"aume 
unter der Auflage deren jeweiliger Kompaktheit verfasst.
Er l"asst aber Chaos zu, das der 
insensitive topologische Ergodensatz ausschliesst.
Ferner ist der Satz von der Existenz der Zimmer offenbar
f"ur metrische Zustandsr"aume sowohl
analog formulierbar als beweisbar, wenn
wir metrische Immanenz und 
metrische Komanenz konzipieren.
\newline
Wir heben auch hervor, dass der insensitive topologische Ergodensatz insofern
etwas mehr als die am Ende des vorangegangenen Abschnittes f"ur eine 
Flussfunktion $\tilde{\Lambda}$ erhobene Frage beantwortet:
Diese Frage nach der Partivit"at des Zustandsraumes in Trajektorienh"ullen war ja f"ur
eine
Flussfunktion $\tilde{\Lambda}$ gestellt, die stetig ist bez"uglich derjeniger 
Produkttopologie der Definitionsmenge $\mathbf{P}_{1}\tilde{\Lambda}=\mathbf{P}_{2}\tilde{\Lambda}\times \mathbb{R}$, die
sich f"ur die Zustandsraumtopologie $\mathbf{T}(\mathbf{P}_{2}\tilde{\Lambda})$ und
die nat"urliche Zahlenstrahltopologie $\mathbf{T}(1)$ 
ergibt. Der insensitive topologische Ergodensatz geht nicht auf die sich auf die gesamte Definitionsmenge der
Flussfunktion $\tilde{\Lambda}$ beziehende Kontinuit"at derselben ein;\index{Produkttopologie}
sondern er
besagt
schon
unter Voraussetzung, dass
lediglich die jeweiligen Fl"usse $\Phi^{t}\in \Gamma_{\Phi}$ stetig sind, dass die 
Identit"at (\ref{cloe}) gilt.
\section{Die Diversit"at der Koh"arenzkriterien}
Als das Erstaunlichste an der Konzeption des Begriffes des Attraktors
empfanden wir von Anfang an, dass der 
Begriff des Attraktors den 
Begriff der Integrit"at transzendiert, wie wir in der Einleitung heraustellten. Da
war die Exposition dieses Erstaunes naturgem"ass 
noch mit einem Dunkel behaftet, das unsere 
Untersuchungen nunmehr teilweise erhellten. Allerdings
w"are unsere 
Aufgabe bei weitem allzu selbstreferent,
h"atte sie
haupts"achlich in der exegetischen Erhellung unserer im Eingangsdunkel
formulierten Motivierungen bestanden. Wir untersuchten
das Verh"altnis von Attraktoren und gemeinsamen invarianten Topologien zueinander
unter m"oglichst allgemeinen topologischen Gesichtspunkten.  
\newline
Der verallgemeinernde Entwurf topologischer Attraktoren greift
das
Koh"arenzkriterium 
auf und er unterstreicht, dass der
entscheidende Handgriff bei der Konstituierung des Begriffes des Attraktors im Koh"arenzkriterium
lokalisiert ist.
Und sogar unter dem dominanten Eindruck
des Erstaunens "uber das Transzendieren des Begriffes
der Integrit"at durch den Begriff des Attraktors 
staunten wir bereits am Anfang
auch
"uber die 
Flexibilit"at der Koh"arenzbedingung an einen
Attraktor, die
in verschiedenen Versionen formulierbar ist. 
Und es ist ja auch verwunderlich,
dass die 
Begriffsbildung des Attraktors
gegen"uber der Variierung ausgerechnet
des \glqq entscheidenden Handgriffes\grqq, der
sie formt, Robustheit zeigt -- und zwar gleichsam in der 
N"ahe des Ausgangspunktes weiterer Generalsierungen, den 
der Fall stetiger reeller 
Flussfunktionen mit kompaktem Zustandsraum darstellt. 
Wir
sprachen bereits da
schon von der auf der Ebene der Begriffsverfassung beobachtbaren
\glqq Metainvarianz des Begriffes des Attraktors gegen"uber der
Ab"anderung der
Koh"arenzbedingung.\grqq
\newline
Diese verwunderliche Metainvarianz gegen"uber der konstitutionellen
Ab"anderung der
Koh"arenzbedingung k"onnen wir n"amlich 
im Fall der Attraktoren stetiger reeller 
Flussfunktionen $\Psi$ mit kompaktem Zustandsraum $\zeta\subset\mathbb{R}^{n}$ 
in
besonders gleichermassen evidenter wie dramatischer Form beobachten: Die Menge der
Fl"usse
$\Gamma_{\Phi}=\{\Psi^{t}:t\in\mathbb{R}\}$ dieser Flussfunktion $\Psi$ legt 
deren invariante Topologie 
$\widehat{\mathbf{T}}(\Gamma_{\Psi})$ fest, die
wegen der Stetigkeit der Flussfunktion $\Psi$
das Mengensystem aller abgeschlossener und
flussinvarianter Mengen des Zustandsraumes $\zeta\subset\mathbb{R}^{n}$
enth"alt.
Da $\zeta$ kompakt ist, sind auch alle 
abgeschlossenen 
flussinvarianten Mengen des Zustandsraumes
Kompakta, sodass das Mengensystem
$$\mathbf{cl}(\widehat{\mathbf{T}}(\Gamma_{\Psi}))=
\{\chi\in \widehat{\mathbf{T}}(\Gamma_{\Psi}):\mathbf{cl}(\chi)=\chi\}$$
ein Mengensystem aller Mengen ist, die als Attraktoren
der Flussfunktion $\Psi$ in Frage kommen. Es sei
$\mbox{mon}^{-}(\mathbb{R})$ bzw. $\mbox{mon}^{+}(\mathbb{R})$ die Menge aller 
streng monoton fallender bzw. wachsender Folgen 
$\{t_{j}\}_{j\in\mathbb{N}}$
des Zahlenstrahles und
$$\bigcup\lbrack\lbrack\Psi\rbrack\rbrack_{a}=\bigcup\{\Psi(a,t):t\in\mathbb{R}\}$$
$$=\bigcup\Bigl\{\bigcup\{\Psi(x,t):x\in a\}:t\in\mathbb{R}\Bigr\}=\bigcup\bigcup\Psi(a,\mathbb{R})\ .$$
Betrachten
wir 
die folgenden Auswahlen aus dem 
Mengensystem aller flussinvarianter 
Kompakta $\mathbf{cl}(\widehat{\mathbf{T}}(\Gamma_{\Psi}))$: Es sei
\begin{equation}\label{llkaren}
\begin{array}{c}
\Bigl\{\chi\in\mathbf{cl}(\widehat{\mathbf{T}}(\Gamma_{\Psi})):
a,b\in(\mathbf{T}(n)\cap\chi)\setminus\{\emptyset\})\Rightarrow
\bigcup\lbrack\lbrack\Psi\rbrack\rbrack_{a}\cap\bigcup\lbrack\lbrack\Psi\rbrack\rbrack_{b}\not=\emptyset\Bigr\}\\
=:\mbox{{\bf @}}_{0}(\Psi)
\end{array}
\end{equation}
ferner
\begin{equation}\label{lrkaren}
\begin{array}{c}
\Bigl\{\chi\in\mathbf{cl}(\widehat{\mathbf{T}}(\Gamma_{\Psi})):\forall\
a,b\in(\mathbf{T}(n)\cap\chi)\setminus\{\emptyset\}\ \exists\ t\in \mathbb{R}:\\
\Psi(a,t)\cap b\not=\emptyset\Bigr\}\\
=:\mbox{{\bf @}}(\Psi)
\end{array}
\end{equation}
und
\begin{equation}\label{rlkaren}
\begin{array}{c}
\Bigl\{\chi\in\mathbf{cl}(\widehat{\mathbf{T}}(\Gamma_{\Psi})):\forall\
a,b\in(\mathbf{T}(n)\cap\chi)\setminus\{\emptyset\}\ \exists\ \{t_{j}\}_{j\in\mathbb{N}}\in\mbox{mon}^{+}(\mathbb{R}):\\
\Psi(a,t_{j})\cap b\not=\emptyset\Bigr\}\\
=:\mbox{{\bf @}}_{+}(\Psi) 
\end{array}
\end{equation}
und schliesslich
\begin{equation}\label{rrkaren}
\begin{array}{c}
\Bigl\{\chi\in\mathbf{cl}(\widehat{\mathbf{T}}(\Gamma_{\Psi})):\forall\
a,b\in(\mathbf{T}(n)\cap\chi)\setminus\{\emptyset\}\ \exists\ \{t_{j}\}_{j\in\mathbb{N}}\in\mbox{mon}^{-}(\mathbb{R}):\\
\Psi(a,t_{j})\cap b\not=\emptyset\Bigr\}\\
=:\mbox{{\bf @}}_{-}(\Psi)\ .
\end{array}
\end{equation}
Es gilt offensichtlich die folgende\newline\newline
{\bf Bemerkung 2.3: Diversit"at der Koh"arenzkriterien}\newline
{\em Es sei $\Psi$ eine
stetige reelle 
Flussfunktion $\Psi$ mit kompaktem Zustandsraum $\zeta\subset\mathbb{R}^{n}$ 
und es seien $\mbox{{\bf @}}_{0}(\Psi)$, $\mbox{{\bf @}}(\Psi)$, $\mbox{{\bf @}}_{+}(\Psi)$ und
$\mbox{{\bf @}}_{-}(\Psi)$ die in} (\ref{llkaren})-(\ref{rrkaren}) {\em festgelegten Teilmengensysteme des 
Mengensystemes
flussinvarianter Kompakta $\mathbf{cl}(\widehat{\mathbf{T}}(\Gamma_{\Psi}))$: Dann ist die
Identit"at 
\begin{equation}\label{karenk}
\mbox{{\bf @}}_{0}(\Psi)=\mbox{{\bf @}}(\Psi)=\mbox{{\bf @}}_{+}(\Psi)=\mbox{{\bf @}}_{-}(\Psi)
\end{equation}
wahr.}
\newline\newline
{\bf Beweis:}\newline
Wenn die Aussage
$$\forall\ a,b\in(\mathbf{T}(n)\cap\chi)\setminus\{\emptyset\}\ \exists\ t\in \mathbb{R}:\ \Psi(a,t)\cap b\not=\emptyset$$
wahr ist, dann ist wegen der Stetigkeit der Flussfunktion $\Psi$ auch die Aussage
$$\forall\ a,b\in(\mathbf{T}(n)\cap\chi)\setminus\{\emptyset\}\ \exists\ t\in \mathbb{R}, \vartheta\in \mathbb{R}^{+}:\ 
\Psi(a,t)\cap b\not=\emptyset\ \land$$
$$(\Psi(a,t+\vartheta)\cap b\not=\emptyset\ \lor\ \Psi(a,t+\vartheta)\cap b=\emptyset)\ \land\ 
\Psi(a,t+\vartheta)\in (\mathbf{T}(n)\cap\chi)\setminus\{\emptyset\}$$
wahr. F"ur alle
$a,b\in(\mathbf{T}(n)\cap\chi)\setminus\{\emptyset\}$ und f"ur alle 
$$t_{1}\in\{t\in \mathbb{R}: \Psi(a,t)\cap b\not=\emptyset\}\not=\emptyset$$
gibt es also Zahlen $t_{2}>t_{1}$ in
$$\{t\in \mathbb{R}: \Psi(a,t)\cap b\not=\emptyset\}\setminus]-\infty,t_{1}]\not=\emptyset\ .$$
Daher ist $\mbox{{\bf @}}(\Psi)=\mbox{{\bf @}}_{+}(\Psi)$. 
Analog gilt wegen der Stetigkeit der Flussfunktion $\Psi$ auch die Identit"at 
$\mbox{{\bf @}}(\Psi)=\mbox{{\bf @}}_{-}(\Psi)$. Die Identit"at
$\mbox{{\bf @}}(\Psi)=\mbox{{\bf @}}_{0}(\Psi)$
gilt, weil
$$\mbox{{\bf @}}(\Psi)=[[\Psi]]\in\mathbf{part}(\mathbf{P}_{2}\Psi)$$
und
$$[[\Psi]]\subset\mathbf{cl}(\widehat{\mathbf{T}}(\Gamma_{\Psi}))$$
eine allgemeine Gleichheit spezifiziert: F"ur jede Menge $X$, f"ur jedes Mengensystem ${\rm T}$, das
$X\subset\bigcup{\rm T}$ "uberdeckt und f"ur jede 
im Mengensystem ${\rm Q}\supset{\rm P}$
enthaltene
Partition ${\rm P}\in\mathbf{part}(X)$ der Menge $X$ ist die Gleichung
$${\rm P}=\Bigl\{
\chi\in{\rm Q}:a,b\in {\rm T}\cap\chi\Rightarrow\bigcup{\rm P}_{a}\cap\bigcup{\rm P}_{b}\not=\emptyset\Bigr\}$$
offensichtlich wahr, falls die Mengensysteme ${\rm T}$ und ${\rm Q}$ folgende 
Bedingung erf"ullen:
F"ur jede Menge $\chi\in{\rm Q}\setminus {\rm P}$ gibt es Mengen $a,b$ des Spurmengensystemes
${\rm T}\cap\chi$, f"ur die
$$\bigcup{\rm P}_{a}\cap\bigcup{\rm P}_{b}=\emptyset$$
gilt. Wir haben also zu zeigen, dass die 
nat"urliche Topologie $\mathbf{T}(n)$ und das Mengensystem der 
Abschl"usse $\mathbf{cl}(\widehat{\mathbf{T}}(\Gamma_{\Psi}))$ der Elemente der
invarianten Topologie $\widehat{\mathbf{T}}(\Gamma_{\Psi})$ so beschaffen sind, dass es
f"ur jede Menge $\chi\in\mathbf{cl}(\widehat{\mathbf{T}}(\Gamma_{\Psi}))\setminus [[\Psi]]$ Mengen $a,b$ des Spurmengensystemes
${\rm T}\cap\chi$ gibt, f"ur die
$$\bigcup[[\Psi]]_{a}\cap\bigcup[[\Psi]]_{b}=\emptyset$$
ist.
Die Abschl"usse aus der Menge $\mathbf{cl}(\widehat{\mathbf{T}}(\Gamma_{\Psi}))$ sind 
dabei Kompakta, weil der Zustandsraum der Flussfunktion $\Psi$ als ein Kompaktum 
bez"uglich der nat"urlichen Topologie $\mathbf{T}(n)$ vorausgesetzt ist. Es gibt eine
Bijektion $b$, die eine Indexmenge $J$ eineindeutig auf die Zimmer der Menge $[[\Psi]]$ abbildet, sodass es
f"ur jedes flussinvariante Kompaktum $\chi\in\mathbf{cl}(\widehat{\mathbf{T}}(\Gamma_{\Psi}))\setminus [[\Psi]]$ 
eine Teilmenge $J(\chi)\subset J$
gibt, f"ur die
$$\chi=\bigcup b(J(\chi))$$
ist.
Wenn es
zwei verschiedene Indizes $j(1),j(2)\in J(\chi)$ gibt, f"ur die es 
zwei
Punkte $x_{1}\in b(j(1)),x_{2}\in b(j(2))$ und
zwei Radien $\delta(1),\delta(2)\in\mathbb{R}^{+}$
zweier offener, zueinander disjunkter Kugeln
$$B(1):=\mathbb{B}_{\delta(1)}(x_{1})\ \mbox{bzw.}\ B(2):=\mathbb{B}_{\delta(2)}(x_{2})$$
der nat"urlichen Topologie
$\mathbf{T}(n)$
um $x_{1}$ bzw. $x_{2}$
gibt,
die so beschaffen sind, dass f"ur alle $k\in J(\chi)$ und f"ur beide $p\in\{1,2\}$ genau dann, 
wenn $k\not=j(p)$ ist, die
Disjunktion
$$\mathbb{B}_{\delta(p)}(x_{p})\cap b(k)=\emptyset$$
gilt, 
dann gilt
$$\bigcup[[\Psi]]_{B(1)}\cap\bigcup[[\Psi]]_{B(2)}=\emptyset\ .$$
Andernfalls, wenn es keine zwei solchen 
Indizes $j(1),j(2)\in J(\chi)$ gibt, liegen
die Zimmer aus $[[\Psi]]$, die Teilmengen des invarianten 
Kompaktums $\chi\in\mathbf{cl}(\widehat{\mathbf{T}}(\Gamma_{\Psi}))\setminus [[\Psi]]$
sind,
zueinander dicht; sodass  
es 
wegen der Kompaktheit der Menge $\chi$ 
f"ur jedes Paar zweier verschiedener 
Indizes $(j(1),j(2))\in J(\chi)\times J(\chi)$ nicht
nur zwei Punkte $x_{1}^{+}\in b(j(1)),x_{2}^{+}\in b(j(2))$ und
und zwei Radien $\rho(1),\rho(2)\in\mathbb{R}^{+}$
zweier offener, zueinander disjunkter Kugeln
$$B(1)^{+}:=\mathbb{B}_{\rho(1)}(x_{1}^{+})\ \mbox{bzw.}\ B(2)^{+}:=\mathbb{B}_{\rho(2)}(x_{2}^{+})$$
der nat"urlichen Topologie
$\mathbf{T}(n)$
um $x_{1}^{+}$ bzw. $x_{2}^{+}$
gibt, die sowohl das Zimmer
$b(j(1))$ als auch das Zimmer $b(j(2))$
schneiden. Es gibt wegen der Kompaktheit der Menge $\chi$ sogar zwei 
solche offenen, zueinander disjunkten Kugeln
$B(1)^{+}$ und $B(2)^{+}$, f"ur welche die
f"ur beide $p\in\{1,2\}$ festgelegten Indexmengen 
$$\{k\in J(\chi):B(p)^{+}\cap b(k)\not=\emptyset\}$$
zueinander disjunkt und jeweils nicht leer sind. Dabei sind
$$B(1)^{+},B(2)^{+}\in(\mathbf{T}(n)\cap\chi)\setminus\{\emptyset\}$$ 
Spurtopologieelemente, f"ur welche die
Disjunktion
$$\bigcup[[\Psi]]_{B(1)^{+}}\cap\bigcup[[\Psi]]_{B(2)^{+}}
=\emptyset$$ 
zutrifft: Allemal ist also das Mengensystem der 
Abschl"usse $\mathbf{cl}(\widehat{\mathbf{T}}(\Gamma_{\Psi}))$ der Elemente der
invarianten Topologie $\widehat{\mathbf{T}}(\Gamma_{\Psi})$ so beschaffen, dass es
f"ur jede Menge $\chi\in\mathbf{cl}(\widehat{\mathbf{T}}(\Gamma_{\Psi}))\setminus [[\Psi]]$ Mengen $a,b$ des Spurmengensystemes
${\rm T}\cap\chi$ gibt, f"ur die
$$\bigcup[[\Psi]]_{a}\cap\bigcup[[\Psi]]_{b}=\emptyset$$
ist.\newline 
{\bf q.e.d.}
\newline
\newline
Im allgemeinen, dann, wenn die Flussfunktion
$\Psi$ die Vorrausetzungen der Bemerkung 2.3 nicht erf"ullt,
ist dabei aber 
das Kriterium, das 
die Auswahl $\mbox{{\bf @}}_{0}(\Psi)$ aus
dem Mengensystem aller flussinvarianter 
Kompakta $\mathbf{cl}(\widehat{\mathbf{T}}(\Gamma_{\Psi}))$ formuliert,
sch"acher als das Kriterium, das 
die Auswahl $\mbox{{\bf @}}(\Psi)$ festlegt, welches
wiederum schw"acher ist als die Kriterien, die 
die Auswahlen $\mbox{{\bf @}}_{+}(\Psi)$ bzw. $\mbox{{\bf @}}_{-}(\Psi)$ festlegen.
\newline
Betrachten wir eine Flussfunktion $\Phi$, deren
Zustandsraum $\mathbf{P}_{2}\Phi=\bigcup\mathcal{Z}$ durch
das Mengensystem $\mathcal{Z}\subset 2^{\mathbf{P}_{2}\Phi}$ "uberdeckt ist,
sodass 
$$[[\Phi]]_{\mathcal{Z}}:=\{\mathbf{cl}_{\mathcal{Z}}(\Phi(z,\mathbb{R})):z\in \mathbf{P}_{2}\Phi\}$$
das Mengensystem derjeniger 
Teilmengen
des Zustandsraumes ist, die wir gegen Ende des n"achsten Kapitels
als die Vorzimmer der Flussfunktion $\Phi$ bez"uglich des Mengensystemes $\mathcal{Z}$ 
einf"uhren und charakterisieren werden.
Die Menge
$\widehat{\mathbf{T}}(\Gamma_{\Phi})$ sei die flussinvariante Topologie der Flussfunktion $\Phi$:
Wir k"onnen 
die Setzungen in (\ref{llkaren})-(\ref{rrkaren})
verallgemeinern zu
\begin{equation}\label{kakren}
\begin{array}{c}
\Bigl\{\chi\in\widehat{\mathbf{T}}(\Gamma_{\Phi}):
a,b\in(\mathcal{Z}\cap\chi)\setminus\{\emptyset\})\setminus\{\emptyset\}\Rightarrow
\bigcup(\lbrack\lbrack\Phi\rbrack\rbrack_{\mathcal{Z}})_{a}\cap\bigcup(
\lbrack\lbrack\Phi\rbrack\rbrack_{\mathcal{Z}})_{b}\not=\emptyset\Bigr\}\\
=:\mbox{{\bf @}}_{0}^{\mathcal{Z}}(\Phi)\ ,\\
\Bigl\{\chi\in\widehat{\mathbf{T}}(\Gamma_{\Phi}):\forall\
a,b\in(\mathcal{Z}\cap\chi)\setminus\{\emptyset\}\ \exists\ t\in \mathbb{R}:\ \Phi(a,t)\cap b\not=\emptyset\Bigr\}\\
=:\mbox{{\bf @}}^{\mathcal{Z}}(\Phi)\ ,\\
\Bigl\{\chi\in\widehat{\mathbf{T}}(\Gamma_{\Phi}):\forall\
a,b\in(\mathcal{Z}\cap\chi)\setminus\{\emptyset\}\ \exists\ \{t_{j}\}_{j\in\mathbb{N}}\in\mbox{mon}^{+}(\mathbb{R}):\ \Phi(a,t_{j})\cap b\not=\emptyset\Bigr\}\\
=:\mbox{{\bf @}}_{+}^{\mathcal{Z}}(\Phi)\ ,\\
\Bigl\{\chi\in\widehat{\mathbf{T}}(\Gamma_{\Phi}):\forall\
a,b\in(\mathcal{Z}\cap\chi)\setminus\{\emptyset\}\ \exists\ \{t_{j}\}_{j\in\mathbb{N}}\in\mbox{mon}^{-}(\mathbb{R}):\ \Phi(a,t_{j})\cap b\not=\emptyset\Bigr\}\\
=:\mbox{{\bf @}}_{-}^{\mathcal{Z}}(\Phi)\ .
\end{array}
\end{equation}
Diese Auswahlen $\mbox{{\bf @}}_{0}^{\mathcal{Z}}(\Phi)$, $\mbox{{\bf @}}^{\mathcal{Z}}(\Phi)$, $\mbox{{\bf @}}_{+}^{\mathcal{Z}}(\Phi)$
und $\mbox{{\bf @}}_{-}^{\mathcal{Z}}(\Phi)$
sind 
im Unterschied zu den Auswahlen gem"ass (\ref{llkaren})-(\ref{rrkaren})
keine solchen Auswahlen aus der flussinvarianten Topologie
$\widehat{\mathbf{T}}(\Gamma_{\Phi})$, von denen wir 
die Erf"ulltheit irgendwelcher
Abgeschlossenheitskriterien 
verlangen.
Die gem"ass (\ref{llkaren})-(\ref{rrkaren})
f"ur die stetige reelle 
Flussfunktion $\Psi$ mit kompaktem Zustandsraum $\zeta\subset\mathbb{R}^{n}$ 
festgelegten
Konstruktionen sind 
\begin{equation}\label{karren}
\begin{array}{c}
\mathbf{cl}(\mbox{{\bf @}}_{0}^{\mathbf{T}(n)}(\Psi))=\mbox{{\bf @}}_{0}(\Psi)\ ,\\
\mathbf{cl}(\mbox{{\bf @}}^{\mathbf{T}(n)}(\Psi))=\mbox{{\bf @}}(\Psi)\ ,\\
\mathbf{cl}(\mbox{{\bf @}}_{+}^{\mathbf{T}(n)}(\Psi))=\mbox{{\bf @}}_{+}(\Psi)\ ,\\
\mathbf{cl}(\mbox{{\bf @}}_{-}^{\mathbf{T}(n)}(\Psi))=\mbox{{\bf @}}_{-}(\Psi)\ .
\end{array}
\end{equation}
F"ur die allgemeinen Konstruktionen gem"ass (\ref{kakren}) 
f"ur eine beliebige 
Flussfunktion $\Phi$, deren
Zustandsraum $\mathbf{P}_{2}\Phi=\bigcup\mathcal{Z}$ durch
das Mengensystem $\mathcal{Z}\subset 2^{\mathbf{P}_{2}\Phi}$ "uberdeckt ist,
gilt statt der in der Bemerkung 2.3 behaupteten Identit"at
(\ref{karenk}) die 
Inklusionenkette
\begin{equation}\label{karrenk} 
\mbox{{\bf @}}_{0}^{\mathcal{Z}}(\Phi)\supset\mbox{{\bf @}}^{\mathcal{Z}}(\Phi)\supset
\mbox{{\bf @}}_{\pm}^{\mathcal{Z}}(\Phi)\ .
\end{equation}
Deswegen nennen wir f"ur jede Flussfunktion $\Phi$ das
Kriterium, dass alle Teilmengen $a,b\in(\mathcal{Z}\cap\chi)\setminus\{\emptyset\}$ f"ur eine jeweilige  
flussinvariante Menge $\chi\in\widehat{\mathbf{T}}(\Gamma_{\Phi})$ 
$$\bigcup(\lbrack\lbrack\Phi\rbrack\rbrack_{\mathcal{Z}})_{a}\cap\bigcup(
\lbrack\lbrack\Phi\rbrack\rbrack_{\mathcal{Z}})_{b}\not=\emptyset$$
erf"ullen, das schwache 
Koh"arenzkriterium,\index{schwaches Koh\"arenzkriterium}
wohingegen wir das 
Kriterium, dass f"ur alle Teilmengen $a,b\in(\mathcal{Z}\cap\chi)\setminus\{\emptyset\}$ f"ur eine jeweilige  
flussinvariante Menge $\chi\in\widehat{\mathbf{T}}(\Gamma_{\Phi})$ eine reelle Zahl $t\in \mathbb{R}$ existiert,
f"ur die
$$\Phi(a,t)\cap b\not=\emptyset$$
gilt,
als das herk"ommliche Koh"arenzkriterium\index{herk\"ommliches Koh\"arenzkriterium}
bezeichnen; und das 
Kriterium, dass f"ur alle Teilmengen $a,b\in(\mathcal{Z}\cap\chi)\setminus\{\emptyset\}$ f"ur eine jeweilige  
flussinvariante Menge $\chi\in\widehat{\mathbf{T}}(\Gamma_{\Phi})$ eine monoton wachsende 
bzw. eine monoton fallende
Zahlenfolge $\{t_{j}\}_{j\in\mathbb{N}}$ in $\mbox{mon}^{+}(\mathbb{R})$ 
bzw. $\mbox{mon}^{-}(\mathbb{R})$ existiert,
f"ur die
$$\Phi(a,t_{j})\cap b\not=\emptyset$$
gilt, nennen
wir das monoton wachsende bzw. monoton fallende Koh"arenzkriterium.\index{monoton wachsendes Koh\"arenzkriterium}
\index{monoton fallendes Koh\"arenzkriterium}
\newline
Die Generalisierungen (\ref{kakren}) unterstreichen
die
Bezogenheit jeweiliger Attraktoren auf eine 
Topologisierung $\mathbf{T}(\mathbf{P}_{2}\Psi)$ der jeweiligen Wertemenge $\mathbf{P}_{2}\Psi$
einer Flussfunktion $\Psi$. Die jeweilige Topologisierung $\mathbf{T}(\mathbf{P}_{2}\Psi)$
ist
in zweifacher Hinsicht f"ur den Begriff des Attraktors konstitutiv:
\newline\newline
{\bf 1.} Zun"achst ist der jeweilige Attraktor ein Kompaktum bez"uglich der 
jeweiligen Topologisierung $\mathbf{T}(\mathbf{P}_{2}\Psi)$
der jeweiligen Wertemenge $\mathbf{P}_{2}\Psi$, die
auch als der Zustandsraum der Flussfunktion $\Psi$
bezeichnet wird. 
Dieses Kompaktheitskriterum an einen Attraktor mag aber vielleicht als eine
willk"urliche Festlegung erscheinen.
\newline
{\bf 2.} Aber zweitens beziehen sich ja
auch 
die Koh"arenzkriterien
auf die jeweilige Relativtopologie $\mathbf{T}(\mathbf{P}_{2}\Psi)\cap\chi$
des jeweiligen Attraktors $\chi$:
Ist beispielsweise das herk"ommliche Koh"arenzkriterium f"ur ein flussinvariantes Kompaktum erf"ullt,
so soll die Reichhaltigkeitsaussage (\ref{verm}) gelten. Und zwar f"ur alle 
$$a,b\in(\mathbf{T}(\mathbf{P}_{2}\Psi)\cap\chi)\setminus\{\emptyset\})=\{{\rm U}\cap\chi:{\rm U}\in\mathbf{T}(\mathbf{P}_{2}\Psi)\}\setminus\{\emptyset\}\ .$$
Falls dabei das flussinvariante Kompaktum $\chi=\{c\}$ ein Fixpunkt\footnote{Ein Fixpunkt ist demnach kein Punkt:
Die gel"aufige Redeweise vom Punkt eines topologischen Raumes ist unmissverst"andlich.
Es koexistiert die ebenfalls verfestigte Sprechweise von einen Punkt $x\in X$,
den der 
jeweilige Kontext von einem lediglichen Element $a\in A$ einer Menge $A$ unterscheidet;
der 
jeweilige Kontext n"amlich,
dass $x$ dabei das Element einer Menge $X$ ist, die topologisiert ist.
Wir bezeichnen hier jede einelementige Teilmenge $\{x\}\subset 2^{\mathbf{P}_{1}f}$ der Potenzmenge der Definitionsmenge
$\mathbf{P}_{1}f$
einer Abbildung $f$ genau dann 
als einen Fixpunkt der Abbildung 
$f$, wenn $f(x)=x$ ist.\index{Fixpunkt einer Abbildung}
\newline
Daher ist der Fixpunkt $\{x\}$ von dem Element $x$ verschieden, 
exakt welches wir als Fixelement der Abbildung $f$ bezeichnen.\index{Fixelement einer Abbildung}
Dieses 
Fixelement $x$ ist kein Fixpunkt, jedoch
ein Punkt, falls
$\mathbf{P}_{1}f$ topologisiert ist. Ein Fixpunkt ist demnach 
eine elementige Fixmenge.\index{Fixmenge einer Abbildung} 
Sodass eine Fixmenge nicht etwa eine Menge von Fixpunkten ist, sondern eine Vereinigung "uber ein Mengensystem von 
Fixpunkten, also eine Menge von Fixelementen.\newline
Dabei auferlegen wir uns der Unmissverst"andlichkeit halber
eine vielleicht unelegante Begriffssklerose und Umst"andlichkeit, die zwar 
innerhalb der Konformit"at bleibt,
nicht aber konform geht mit den
verfestigten Sprechweisen, die fl"ussiger sind und die auf bereits vielbewandertem Terrain auch zweckm"assiger sind als
unsere Redeweise.}
ist, ist
$$\mathbf{T}(\mathbf{P}_{2}\Psi)\cap\chi=\{\emptyset,\chi\}$$
die gr"obste Topologie. 
Dann gilt jedes Koh"arenzkriterium.
Dass das schwache Koh"arenzkriterium 
schw"acher als das herk"ommliche 
Koh"arenzkriterium und auch schw"acher als die monotonen Koh"arenzkriterien
ist, zeigt sich
nicht anhand von Fixpunkten, die sowohl dem  
schwachen als auch dem herk"ommlichen Koh"arenzkriterium und die sogar
den monotonen Koh"arenzkriterien
gen"ugen.
Gleich, wie wir den Begriff des Attraktors verfassen: Gleich, ob wir von einem flussinvarianten Kompaktum
verlangen, dass es dem schwachen Koh"arenzkriterium, dem herk"ommlichen Koh"arenzkriterium oder den
monotonen Koh"arenzkriterien
gen"ugen soll, damit es ein Attraktor sei:
Fixpunkte sind demnach allemal Attraktoren. 
\newline
Zyklen, geschlossene Trajektorien eines 
endlichdimensionalen $\mathbb{R}^{n}$ f"ur $n\geq 2$
hingegen, sind keine 
Attraktoren -- bez"uglich der nat"urlichen Topologie $\mathbf{T}(n)$, weil sie das herk"ommliche Koh"arenzkriterium
nicht erf"ullen, wie die unmittelbare Anschauung lehrt.
Sie erf"ullen dabei aber durchaus das schwache Koh"arenzkriterium,
was zeigt, dass 
das schwache Koh"arenzkriterium
ein echt schw"acheres   
Kriterium ist als das herk"ommliche Koh"arenzkriterium.
\newline
Wir tun hier gerade so, als h"atten wir die Freiheit, willk"urlich einen jeweiligen 
Begriff des Attraktors festzulegen.
Und in der Tat: Wir legen zwar als den Attraktor einer Flussfunktion
gerade diejenigen flussinvarianten Kompakta fest, die das herk"ommliche      
Koh"arenzkriterium erf"ullen und wir schliessen
damit Zyklen davon aus, als Attraktoren zu gelten. Wir k"onnten
auch das schwache Koh"arenzkriterium zur
Definition der Attraktoren heranziehen. 
Aber warum sollten wir hier keinen differenzierten Pluralismus einrichten, der
sich naturgem"ass ergibt und dessen 
Differenzierung die
Szenarien der nicht nur stetigen und reellen Flussfunktionen, sondern
der allgemeinen Flussfunktionen
illustrieren?\newline
Wir k"onnen diesen Pluralismus zulassen und weisen auf die Spielart der 
Attraktoren senso latto hin, die
L.Garrido u. C.Sim$\acute{\mbox{o}}$ in \cite{garr}
verfassten.
Wir k"onnen
alle diejenigen flussinvarianten Kompakta einer Flussfunktion, die das 
schwache Koh"arenzkriterium bzw. ein monotones Koh"arenzkriterium erf"ullen, 
die schwachen Attraktoren bzw. die monotonen Attraktoren der jeweiligen Flussfunktion nennen; alle 
diejenigen flussinvarianten Kompakta einer Flussfunktion, die das 
herk"ommliche Koh"arenzkriterium erf"ullen, k"onnen wir hingegen
als die 
herk"ommlichen Attraktoren der jeweiligen 
Flussfunktion bezeichnen.\index{herk\"ommliche Attraktoren}\index{schwache Attraktoren einer Flussfunktion}
\index{Attraktoren einer Flussfunktion}
Zyklen sind demnach zwar keine Attraktoren, sie sind jedoch schwache Attraktoren.
Die Zusammenh"ange zwischen den diversen Attraktorbegriffen wollen
wir in diesem Rahmen allerdings nicht ausforschen.
\newline
Halten wir dies fest: 
Es ergibt sich die logische Diversit"at verschiedener 
Attraktorenbegriffe: Wir sind mit einer eigent"umlichen begrifflichen Verfassungsfreiheit konfrontiert.
Der Begriff des Attraktors erscheint uns dabei bislang in zunehmendem Mass als ein wesentlich topologischer Begriff.
Die Betrachtung des n"achsten Abschnittes verst"arkt diesen Eindruck.
\section{Die Kovarianz der Attraktoren}
Wir ziehen uns 
gerne
auf die Attraktoren endlichdimensionaler reeller R"aume zur"uck,
mithin auf reelle Flussfunktionen -- und zwar nicht nur
wegen des leichten Einstieges in die 
Diskussion einer jeweiligen Frage auf der Grundlage der vertrauten
nat"urlichen Topologien; meistens sollte 
die Er"orterung einer jeweiligen Frage in vertrautem Rahmen
der erste Schritt sein.
Wir arbeiten 
schon alleine deshalb gerne im $\mathbb{R}^{n}$ f"ur 
eine jeweilige endliche Dimension $n\in\mathbb{N}$, weil
von unseren Sprechenweisen  eine Last f"allt:
Wir d"urfen es unterlassen, Flussfunktionen ausdr"ucklich
als reelle
Flussfunktionen zu bezeichnen. Vor allem aber k"onnen
wir einfach, d.h., ohne ausdr"uckliche Erw"ahnung jeweiliger Bezugstopologien, von Stetigkeit, von Kompaktheit und von Offenheit reden. Wir beziehen uns immer
auf die entsprechenden nat"urlichen Topologien $\mathbf{T}(n)$ des $\mathbb{R}^{n}$ f"ur 
eine jeweilige endliche Dimension $n\in\mathbb{N}$. 
\newline
Wir k"onnen uns bei dem Nachweis, dass eine jeweilige Zustandsmenge des Zustandsraumes einer
stetigen Flussfunktion 
einer deren Attraktoren ist, auf die
Betrachtung entsprechender Attraktoren 
stetig differenzierbarer Flussfunktionen $\Psi$ beschr"anken; d.h., auf 
Flussfunktionen, erstens deren Fl"usse\index{stetig differenzierbare Flussfunktion}
$\Psi^{t}$ f"ur alle $t\in\mathbb{R}$ stetig differenzierbar
sind
und f"ur die zweitens f"ur alle $z\in\zeta$ die parzielle Ableitung
$\partial_{2}\Psi(z,\mbox{id})$ stetig ist:\newline
Denn f"ur jeden Hom"oomorphismus $f$ des $\mathbb{R}^{n}$ auf sich oder in einen $\mathbb{R}^{m}$ mit $m\geq n$
ist
$\chi\subset\mathbf{P}_{2}\psi$ offenbar genau dann ein Attraktor der stetigen  
Flussfunktion $\psi$ mit kompaktem Zustandsraum $\zeta\subset\mathbb{R}^{n}$, wenn
$f(\chi)$ ein Attraktor der Flussfunktion $\psi_{f}$
ist.
F"ur jede Flussfunktion
$\Phi$ und jede auf deren 
Zustandsraum $\zeta$ definierter Bijektion $b$
sei dabei
\begin{equation}\label{mamma}
\Phi_{b}:=b\circ\Phi\circ((b^{-1}\circ\mathbf{P}_{1})\oplus\mathbf{P}_{2})\ ,
\end{equation}
wobei $\oplus$ die Kokatenation bezeichnet.
Die Eigenschaft einer Teilmenge eines entsprechenden Zustandsraumes, ein Attraktor einer jeweiligen
Flussfunktion zu sein, ist insofern eine Kovariante gegen"uber allen solchen Hom"oomorphismen $f$.
Es liegt daher nahe, von der topologischen Kovarianz reeller Attraktoren zu 
sprechen.\index{topologische Kovarianz der Attraktoren}
\newline
Die topologische Kovarianz der Attraktoren gilt offensichtlich nicht nur in dem Fall, dass
nat"urliche Topologien zugrundeliegen.
Wir werden zeigen, dass 
wegen des verallgemeinerten insensitiven Ergodensatzes 3.3
die Kovarianz der Attraktoren sehr allgemein gegeben ist.\index{verallgemeinerter insensitiver Ergodensatz} 
Die allgemeine topologische Kovarianz der Attraktoren einer Flussfunktion $\Psi$ fassen wir als deren 
Eigenschaft, den Attraktoren einer 
analog zu (\ref{mamma}) transformierten 
Flussfunktion $\Psi_{b}$ zu entsprechen,
falls die Bijektion $b$ ein Hom"oomorphismus ist.
Wir wollen hier nicht auf die allgemeine topologische Kovarianz der Attraktoren 
eingehen.\newline Betrachten wir eine reelle Flussfunktion
$\Psi$, deren kompakter Zustandsraum $\zeta\subset\mathbb{R}^{n}$
endlichdimensional ist, und bezeichnen wir mit 
$\mbox{{\bf @}}(\Psi)$ das Mengensystem all
deren Attraktoren und mit
$\mbox{{\bf @}}^{S}(\Psi)$ das Mengensystem all
deren sensitiver Attraktoren:
Die topologische Kovarianz der Attraktoren der Flussfunktion
$\Psi$ ist der Sachverhalt, dass f"ur jeden
Hom"oomorphismus $f$ des $\mathbb{R}^{n}$ auf sich oder in einen $\mathbb{R}^{m}$ mit $m\geq n$
\begin{equation}
\mbox{{\bf @}}(\Psi_{f})=f(\mbox{{\bf @}}(\Psi))
\end{equation}
ist, wobei $f(\mbox{{\bf @}}(\Psi))=\{f(\chi):\chi\in\mbox{{\bf @}}(\Psi)\}$ ist. Mehr noch: 
Wir finden, 
dass damit auch die 
sensitiven Attraktoren in dem Sinn
Kovariante gegen"uber jedem
Hom"oomorphismus $f$ des $\mathbb{R}^{n}$ auf sich oder in einen $\mathbb{R}^{m}$ mit $m\geq n$
sind, dass
\begin{equation}\label{mammaz}
\mbox{{\bf @}}^{S}(\Psi_{f})=f(\mbox{{\bf @}}^{S}(\Psi))
\end{equation}
ist.
Denn nach dem Traktat "uber den elementaren
Quasiergodensatz \cite{raab} gilt ja f"ur jede reelle Flussfunktion
$\Psi$, deren kompakter Zustandsraum $\zeta\subset\mathbb{R}^{n}$
endlichdimensional ist, dass nach dem Satz von der Sensitivit"at\index{Sensitivit\"at} nicht trivialer Zimmer $(2.2.2)^{0}$
der Abhandlung des elementaren Quasierdodensatzes \cite{raab}
\begin{equation}\label{mammay}
\mbox{{\bf @}}(\Psi)\setminus[\Psi]=\mbox{{\bf @}}^{S}(\Psi)
\end{equation}
ist.
Allerdings m"ussten wir dazu, die topologische Kovarianz der sensitiven Attraktoren (\ref{mammaz}) einzusehen,
nicht auf das wichtige Resultat (\ref{mammay}) jenes Traktates
zur"uckgreifen: Denn schon unmittelbar aus der
Definition des sensitiven Attraktors ist die topologische
Kovarianz der sensitiven Attraktoren (\ref{mammaz})
ablesbar.\index{topologische Kovarianz der sensitiven Attraktoren}\index{topologische Kovarianz der Attraktoren}
\newline
Es ist f"ur die sensitiven Attraktoren wie f"ur die Attraktoren: 
Weil die Eigenschaft einer Teilmenge eines Zustandsraumes, ein seltamer Attraktor einer jeweiligen
reellen
Flussfunktion zu sein, eine Kovariante gegen"uber den betrachtenen Hom"oomorphismen $f$ ist,
k"onnen wir uns bei der Untersuchung
jeweiliger Gesichtspunkte
der sensitiven Attraktoren 
stetiger reeller 
Flussfunktionen mit kompaktem Zustandsraum  
auf stetig differenzierbare Flussfunktionen beschr"anken.
Die topologische Kovarianz sensitiver Attraktoren\index{topologische Kovarianz} 
stetiger reeller 
Flussfunktionen mit kompaktem Zustandsraum ergibt sich unmittelbar aus
dem Satz von der Sensitivit"at nicht trivialer Zimmer $(2.2.2)^{0}$ des
Traktates
"uber den elementaren Quasierdodensatz \cite{raab}.\newline
Aber auch, wenn wir uns 
bei allgemeineren Gegebenheiten
nicht schon auf die S"atze jenes Traktats berufen k"onnen
und wenn 
wir nicht 
von Attraktoren 
stetiger reeller
Flussfunktionen reden, sondern von 
den Attraktoren 
einer allgemeinen Flussfunktion 
$\Phi\in A^{A\times \mathbb{R}}$ mit der Wertemenge $A$, welche die Topologie $\mathbf{T}(A)$ topologisiere:
F"ur jeden Hom"oomorphismus $\phi$ des Zustandsraumes $A$ auf $Z$
bez"uglich der Topologien $\mathbf{T}(A)$ und $\mathbf{T}(Z)$ gilt, dass
die Eigenschaft einer Menge $\chi$, eine Fixmenge jedes Phasenflusses $\{\Psi^{t}\}_{t\in\mathbb{R}}$
zu sein, genau dann vorliegt, wenn
$\phi(\chi)$ eine Fixmenge jedes Phasenflusses 
$\{(\phi\circ\Phi)^{t}\}_{t\in\mathbb{R}}$ ist; und 
die Eigenschaft einer Menge $\chi$, das Koh"arenzkriterium bez"uglich der Topologie $\mathbf{T}(A)$
zu erf"ullen, liegt genau dann vor, wenn
$\phi(\chi)$ das Koh"arenzkriterium\index{Koh\"arenzkriterium} bez"uglich der Topologie $\mathbf{T}(Z)$
erf"ullt. Insofern ist die Eigenschaft einer Teilmenge eines entsprechenden Zustandsraumes, ein Attraktor einer jeweiligen allgemeinen
Flussfunktion zu sein, eine topologische Kovariante per conceptionem.\newline
Die analoge allgemeine Aussage k"onnen wir hingegen f"ur sensitive Attraktoren nicht aussprechen,
solange 
der Begriff des sensitiven Attraktors ist nicht in der totalen topologischen Generalit"at konstituiert ist.
Wir kennen bislang lediglich sensitive Attraktoren von 
Flussfunktionen $\Phi$ bez"uglich einer jeweiligen Metrik $d$ als
Metrik der Wertemenge $A$ der Flussfunktion
$\Phi\in A^{A\times \mathbb{R}}$.
\newline
Die topologische Kovarianz der Attraktoren
scheint den Begriff des Attraktors als einen 
Begriff zu zeigen, der nicht nur 
in der Sprache der
Topologie verfasst ist. Der Begriff des Attraktors kommt uns 
angesichts der topologischen Kovarianz der Attraktoren
insofern wie ein wesentlich topologischer Begriff vor, als ihn seine topologische Kovarianz 
sachlich
in die
Topologie involviert.\newline

\chapter{Verallgemeinerungen, welche die Abschlusskommutation\\ nahelegt}\label{ddcln}
\section{Der verallgemeinerte insensitive Ergodensatz}
{\small Bald wird uns der Begriff des Attraktors nicht mehr als ein wesentlich topologischer Begriff erscheinen und wir 
werden es der Begriffsbildung des freien Attraktors danken, dass sie im 
Gegensatz zu der Begriffsbildung des topologischen Attraktors den Verallgemeinertheitsgrad hat, in welchem
es m"oglich ist, den verallgemeinerten insensitiven Ergodensatz zu formulieren. 
Dieser liegt auf der Allgemeinheitsstufe der Mengenlehre;
er darf sich zu den Errungenschaften der \glqq Grundz"uge der Mengenlehre\grqq Hausdorffens \cite{haus}
oder der \glqq Th\'eorie des ensembles\grqq \cite{pour}
des Bourbaki-Kreises hinzugesellen. Wir weisen hierbei darauf hin, dass wir 
hier {\em keine Mengenlehre} betreiben:
F"ur die Verfassung des 
verallgemeinerten insensitiven Ergodensatzes und seiner Erschliessung f"ur die Applikation brauchen wir keinerlei Kenntnisse 
desjenigen Teiles der axiomatischen Mengenlehre, der den Aufbau der Mengenlehre selber zum Thema hat.
Wir gehen hier von der Mengenlehre aus, ohne diese zu problematisieren. Dieses Kapitel ist deutlich das 
interessanteste des ersten Teiles der abstrakten Ergodentheorie.} 
\newline\newline
Es l"asst sich 
im Fall der Stetigkeit der Fl"usse $\Phi^{t}\in \Gamma_{\Phi}$ gegen"uber einer beliebigen 
Topologie $\mathbf{T}(\mathbf{P}_{2}\Phi)$ 
die 
f"ur die reelle, stetige Flussfunktion $\Lambda$ mit einem kompakten Zustandsraum 
g"ultige
Identit"at
$$\underline{\mbox{{\bf @}}}(\mathbf{T}(\Gamma_{\Lambda}),\mathbf{T}(\dim\Lambda),\mathbf{C}(\mathbf{T}(\dim\Lambda)))=[[\Lambda]]$$
zu der Identit"at
$$\underline{\mbox{{\bf @}}}(\mathbf{T}(\Gamma_{\Phi}),
\mathbf{T}(\mathbf{P}_{2}\Phi),\mathbf{C}(\mathbf{T}(\mathbf{P}_{2}\Phi)))=[[\Phi]]_{\mathbf{T}(\mathbf{P}_{2}\Phi))}\cap
\mathbf{C}(\mathbf{T}(\mathbf{P}_{2}\Phi))$$
verallgemeinern, wie wir mit dem verallgemeinerten insensitiven Ergodensatz 3.3\index{verallgemeinerter insensitiver Ergodensatz} zeigen werden. Wenn 
die Kompaktheit des Zustandsraumes bez"uglich der Topologie des
Zustandsraumes $\mathbf{T}(\mathbf{P}_{2}\Phi)$ gegeben ist,
bez"uglich welcher die Fl"usse $\Phi^{t}\in \Gamma_{\Phi}$ stetig sind,
gilt insofern vollst"andige Analogie zu dem Fall der reellen, stetigen Flussfunktion $\Lambda$, deren
Zustandsraum bez"uglich der nat"urlichen Topologie $\mathbf{T}(\dim\Lambda)$ kompakt ist, 
als
\begin{equation}
\begin{array}{c}
[\Phi]=\mbox{{\bf @}}(\Phi,\widehat{\mathbf{T}}(\Gamma_{\Phi}))\\
\quad\\
\subset\mbox{{\bf @}}(\Phi,\mathbf{T}(\mathbf{P}_{2}\Phi))=
\underline{\mbox{{\bf @}}}(\Phi,\widehat{\mathbf{T}}(\Gamma_{\Phi}),\mathbf{C}(\mathbf{T}(\mathbf{P}_{2}\Phi)))\\
=[\Phi]\cup[[\Phi]]\ ,
\end{array}
\end{equation} 
\begin{equation}
\begin{array}{c}
\underline{\mbox{{\bf @}}}(\mathbf{T}(\Gamma_{\Phi}),2^{\mathbf{P}_{2}\Phi},\mathbf{C}(\mathbf{T}(\mathbf{P}_{2}\Phi)))
=[\Phi]\cap\mathbf{C}(\mathbf{T}(\mathbf{P}_{2}\Phi))\ ,
\end{array}
\end{equation}
\begin{equation}
\begin{array}{c}
\mbox{{\bf @}}(\Phi,\mathbf{T}(\mathbf{P}_{2}\Phi))=
\underline{\mbox{{\bf @}}}(\Phi,\widehat{\mathbf{T}}(\Gamma_{\Phi}),\mathbf{C}(\mathbf{T}(\mathbf{P}_{2}\Phi)))=[[\Phi]]\ ,
\end{array}
\end{equation} 
f"ur
\begin{equation}
[[\Phi]]_{\mathbf{T}(\mathbf{P}_{2}\Phi)}=\{\mathbf{cl}_{\mathbf{T}(\mathbf{P}_{2}\Phi)}(\gamma):\gamma\in[\Phi]\}
\end{equation} 
gilt. Zum verallgemeinerten insensitiven Ergodensatz 3.3\index{verallgemeinerter insensitiver Ergodensatz} kommen wir auf die folgende Weise:
Wir k"onnen die kollektive Abschlusstransitivit"at stetiger Fl"usse leicht verallgemeinern: Es sei f"ur jedes Mengensystem $\mathcal{A}$
\begin{equation}
\mathcal{A}^{c}:=\Bigl\{(\bigcup\mathcal{A})\setminus {\rm A}:{\rm A}\in\mathcal{A}\Bigr\} 
\end{equation}
das jeweilige zu dem Mengensystem $\mathcal{A}$ komplement"are Mengensystem\index{komplement\"ares Mengensystem}
und 
der Operator 
\begin{equation}\label{cloee}
\mathbf{cl}_{\mathcal{A}}:=\bigcap\{A\in\mathcal{A}^{c}:A\supset\mbox{id}\}
\end{equation} auf seiner jeweiligen 
Definitionsmenge 
\begin{equation}\label{clocee}
\mathbf{P}_{1}\mathbf{cl}_{\mathcal{A}}=2^{\bigcup\mathcal{A}}
\end{equation}
eingef"uhrt, exakt welchen 
wir als den H"ullenoperator des Mengensystemes $\mathcal{A}$ bezeichnen.\index{H\"ullenoperator eines Mengensystemes}
Dabei ist $\bigcap\emptyset=\emptyset$. Daher ist auch das Mengensystem
\begin{equation}\label{cloocee}
\overline{\mathcal{A}}:=\Bigl\{\mathbf{cl}_{\mathcal{A}}(X):X\subset\bigcup\mathcal{A}\Bigr\} 
\end{equation}
aller bez"uglich $\mathcal{A}$ abgeschlossener Mengen festgelegt.\newline
Wir erwarten, dass
dieser Operator wegen der uneingeschr"ankten Allgemeinheit des ihn jeweils
indizierenden Mengensystemes $\mathcal{A}$ nur wenigen allgemein verfassbaren 
Aussagen gen"ugt. Und in der
Tat zeigt sich, dass zu jeder
der Aussagen
\begin{equation}\label{iduwed}
\begin{array}{c}
\mathbf{cl}_{\mathcal{A}}(A\cap_{j} B)\supset_{k}\mathbf{cl}_{\mathcal{A}}(A)\cap_{j}\mathbf{cl}_{\mathcal{A}}(B)\ ,\\
\end{array}
\end{equation}
f"ur alle Indexpaare
$$(j,k)\in\{0,1\}^{2}$$
und
f"ur die Setzungen (\ref{idned})
\begin{equation}\label{idned}
\begin{array}{c}
(\in_{0},\in_{1}):=(\in,\not\in)\ ,\\
(\supset_{0},\supset_{1}):=(\supset,\subset)\ ,\\
(\cap_{0},\cap_{1}):=(\cap,\cup)
\end{array}
\end{equation} einfache
Gegenbeispiele existieren, die
mit Hilfe 
wenige Mengen umfassender 
Mengensysteme
formulierbar sind. Immerhin
gilt, gleich, welches
Mengensystem $\mathcal{A}$ geben ist, die Idempotenz des 
H"ullenoperators eines Mengensystemes $\mathcal{A}$: Es ist
\begin{equation}\label{ydwed}
\mathbf{cl}_{\mathcal{A}}=\mathbf{cl}_{\mathcal{A}}\circ\mathbf{cl}_{\mathcal{A}}\ ,
\end{equation}
denn es gilt f"ur 
jede Menge $A\in \mathbf{P}_{1}\mathbf{cl}_{\mathcal{A}}$ und
f"ur
jedes Element 
$$z\in \mathbf{cl}_{\mathcal{A}}(\mathbf{cl}_{\mathcal{A}}(A))\setminus\mathbf{cl}_{\mathcal{A}}(A)$$
die Aussage
$$Y\in\mathcal{A}^{c}\ \land\  Y\subset\mathbf{cl}_{\mathcal{A}}(A)\Rightarrow z\in Y\ \land$$
$$\exists\ X\in\{Y\in\mathcal{A}^{c}:Y\supset A\}:\ z\not\in X\ ,$$
wobei die Implikation
$$X\in\{Y\in\mathcal{A}^{c}:Y\supset A\}\Rightarrow X\supset\mathbf{cl}_{\mathcal{A}}(A)$$
wahr ist. Das Element $z$ konkretisierte also einen Widerspruch, sodass
$$\mathbf{cl}_{\mathcal{A}}(\mathbf{cl}_{\mathcal{A}}(A))\setminus\mathbf{cl}_{\mathcal{A}}(A)=\emptyset$$
ist. 
Der generische Fall ist so beschaffen, dass es Mengen $A\in\mathbf{P}_{1}\mathbf{cl}_{\mathcal{A}}$
gibt, zu denen
kein
Element $A^{+}$ 
des Mengensystemes $\mathcal{A}^{c}$ existiert, das $A$ enth"alt.
Wenn es keine Menge $A^{+}$ in dem zu $\mathcal{A}$ komplement"aren Mengensystem $\mathcal{A}^{c}$ gibt, die 
$A$ enth"alt, ist die H"ulle $\mathbf{cl}_{\mathcal{A}}(A)$ leer. 
Der Fall, dass f"ur ein besonderes Mengensystem $\mathcal{A}_{\star}$ die Identit"at
\begin{equation}\label{vidwed}
\mathbf{cl}_{\mathcal{A}_{\star}}^{-1}(\{\emptyset\})=\{\emptyset\}
\end{equation}
gegeben ist,
ist in dieser Hinsicht die bezeichnenswerte Ausnahme vom generischen Szenario. 
Es gilt allerdings f"ur jedes Mengensystem $\mathcal{A}$
\begin{equation}\label{iiwed}
\mathbf{P}_{1}\mathbf{cl}_{\mathcal{A}}\setminus\mathbf{cl}_{\mathcal{A}}^{-1}(\{\emptyset\})\supset\mathcal{A}^{c}\ . 
\end{equation}
Genau dann, wenn f"ur ein Mengensystem $\mathcal{A}_{\star}$ die 
Gleichung (\ref{vidwed}) 
gilt, nennen wir es vollst"andig.\index{vollst\"andiges Mengensystem}
Jedes Mengensystem $\mathcal{A}$, f"ur das
$$\emptyset\in\mathcal{A}$$
gilt, ist vollst"andig. Jede Topologie und jede $\sigma$-Algebra  
ist demnach vollst"andig.
Falls ein Mengensystem $\mathcal{A}$ so beschaffen ist, dass
das zu ihm komplement"are Mengensystem $\mathcal{A}^{c}$
gegen"uber 
der Vereinigung je zweier seiner Mengen abgeschlossen 
ist, 
ist $\mathcal{A}$ nicht notwendigerweise vollst"andig. $\mathcal{A}$ ist aber vollst"andig,
wenn $\mathcal{A}$ endlich und $\mathcal{A}^{c}$ dabei abgeschlossen 
gegen"uber 
der Vereinigung je zweier seiner Mengen ist. 
Wenn ein Mengensystem 
$\mathcal{A}_{\star}$
vollst"andig ist und dabei
dessen
komplement"ares Mengensystem $\mathcal{A}_{\star}^{c}$
gegen"uber 
der Vereinigung bzw. dem Schnitt je zweier seiner Mengen
abgeschlossen ist,
gelten f"ur alle 
$A,B\in \mathbf{P}_{1}\mathbf{cl}_{\mathcal{A}_{\star}}$
die Aussagen 
\begin{equation}\label{idwed}
\begin{array}{c}
\mathbf{cl}_{\mathcal{A}_{\star}}(A\cup B)\supset\mathbf{cl}_{\mathcal{A}_{\star}}(A)\cup\mathbf{cl}_{\mathcal{A}_{\star}}(B)\\
\mbox{bzw.}\\
\mathbf{cl}_{\mathcal{A}_{\star}}(A\cap B)\subset\mathbf{cl}_{\mathcal{A}_{\star}}(A)\cap\mathbf{cl}_{\mathcal{A}_{\star}}(B)\ .
\end{array}
\end{equation}
Mit den Setzungen (\ref{idned})
k"onnen wir parallel beweisen:
F"ur alle 
$A,B\in \mathbf{P}_{1}\mathbf{cl}_{\mathcal{A}_{\star}}$ und f"ur alle Elemente $x\in\bigcup \mathcal{A}_{\star}$ 
gilt f"ur beide Indizes $j\in\{0,1\}$ die
"Aquivalenz 
\begin{displaymath}
\begin{array}{c}
x\in_{j}\mathbf{cl}_{\mathcal{A}_{\star}}(A)\cap_{j}\mathbf{cl}_{\mathcal{A}_{\star}}(B)\\ 
\Leftrightarrow\\ 
(\exists\ X\in\{a\in\mathcal{A}_{\star}^{c}:a\supset A\}:\  x\in_{j} X\ \lor\ \bigcap\{a\in\mathcal{A}_{\star}^{c}:a\supset A\}=\emptyset)\  \land\\
(\exists\ Y\in\{a\in\mathcal{A}_{\star}^{c}:b\supset B\}:\ x\in_{j} Y\ \lor\ \bigcap\{b\in\mathcal{A}_{\star}^{c}:b\supset B\}=\emptyset)\ .\\
\end{array}
\end{displaymath}
Wegen der Vollst"andigkeit des Mengensystemes $\mathcal{A}_{\star}^{c}$
ist dabei
$$\bigcap\{a\in\mathcal{A}_{\star}^{c}:a\supset A\}\not=\emptyset\ \land\ \bigcap\{b\in\mathcal{A}_{\star}^{c}:b\supset B\}\not=\emptyset\ .$$
Also gilt die Implikation
\begin{displaymath}
\begin{array}{c}
x\in_{j}\mathbf{cl}_{\mathcal{A}_{\star}}(A)\cap_{j}\mathbf{cl}_{\mathcal{A}_{\star}}(B)\\
\Rightarrow\\
(\exists\ X\in\{a\in\mathcal{A}_{\star}^{c}:a\supset A\}\ \land\
\exists\ Y\in\{b\in\mathcal{A}_{\star}^{c}:b\supset B\}:\ x\in_{j} X\cap_{j} Y)\ .\\
\end{array}
\end{displaymath}
Wobei wegen der Abgeschlossenheit des Mengensystemes $\mathcal{A}_{\star}^{c}$ 
gegen"uber 
dem Schnitt je zweier seiner Mengen im Fall $j=0$ und gegen"uber 
der Vereinigung je zweier seiner Mengen im Fall $j=1$
die Implikationen
\begin{displaymath}
\begin{array}{c}
X\in\{a\in\mathcal{A}_{\star}^{c}:a\supset A\}\ \land\
Y\in\{b\in\mathcal{A}_{\star}^{c}:b\supset B\}
\Rightarrow\\
X\cap _{j}Y\in\mathcal{A}_{\star}^{c}\ \land\ X\cap_{j} Y\supset A\cap_{j} B
\end{array}
\end{displaymath}
und 
\begin{displaymath}
\begin{array}{c}
X\cap_{j}Y\in\mathcal{A}_{\star}^{c}\ \land\ X\cap_{j} Y\supset A\cap_{j} B
\Rightarrow\\
X\cap_{j} Y\in\{c\in \mathcal{A}_{\star}^{c}:c\supset A\cap_{j} B\}
\end{array}
\end{displaymath}
wahr sind,
sodass schliesslich die Implikationen
\begin{displaymath}
\begin{array}{c}
x\in_{j}\mathbf{cl}_{\mathcal{A}_{\star}}(A)\cap_{j}\mathbf{cl}_{\mathcal{A}_{\star}}(B)
\Rightarrow\\
\exists\ Z\in\{c\in \mathcal{A}_{\star}^{c}:c\supset A\cap_{j} B\}:\ x\in_{j} Z 
\end{array}
\end{displaymath}
und
\begin{displaymath}
\begin{array}{c}
\exists\ Z\in\{c\in \mathcal{A}_{\star}^{c}:c\supset A\cap_{j}B\}:\ x\in_{j} Z\Rightarrow\\ 
x\in_{j}\mathbf{cl}_{\mathcal{A}_{\star}}(A\cap_{j} B)
\end{array}
\end{displaymath}
wahr sind, womit die Implikation
\begin{displaymath}
\begin{array}{c}
x\in_{j}\mathbf{cl}_{\mathcal{A}_{\star}}(A)\cap_{j}\mathbf{cl}_{\mathcal{A}_{\star}}(B)
\Rightarrow\\
x\in_{j}\mathbf{cl}_{\mathcal{A}_{\star}}(A\cap_{j} B)
\end{array}
\end{displaymath}
verfiziert ist; und damit ist im Fall $j=0$ die Submultiplikativit"at und im Fall $j=1$ die
Supraadditivit"at des H"ullenoperators des Mengensystemes $\mathcal{A}_{\star}$ best"atigt.
Jede Topologie ist ein
vollst"andiges Mengensystem $\mathcal{A}_{\star}$, das
abgeschlossen 
ist sowohl gegen"uber der Vereinigung als auch gegen"uber dem Schnitt je zweier seiner Mengen,
sodass auch dessen komplement"ares Mengensystem $\mathcal{A}^{c}$
sowohl gegen"uber der Vereinigung als auch gegen"uber dem Schnitt je zweier seiner Mengen
abgeschlossen
ist.
In dem speziellen Fall, dass 
das Mengensystem $\mathcal{A}_{\star}$ eine Topologie ist,
sind die Aussagen (\ref{idwed}) Implikationen
der bekannten Befunde (\ref{ijedwed}) 
f"ur Topologien.
Die Aussagen 
(\ref{ijedwed})
sind mit der selben Argumentation best"atigbar, mit der
wir (\ref{idwed}) zeigten.
Deswegen nennen wir exakt jedes vollst"andige
Mengensystem $\mathcal{A}_{\star}$, das
"uberdies abgeschlossen 
ist sowohl gegen"uber der Vereinigung als auch gegen"uber dem Schnitt je zweier seiner Mengen
eine Quasitopologie.\index{Quasitopologie} Jede $\sigma$-Algebra ist demnach eine
Quasitopologie und jede
endliche Quasitopologie ist demnach eine Topologie.
Das Kriterium, abgeschlossen gegen"uber
der Vereinigung
beliebiger Teilmengen zu sein, das diejenigen Quasitopologien
erf"ullen, die 
Topologien sind,
erlaubt es, weitere Schl"usse zu ziehen:\newline
Wenn $\mathbf{T}$
eine Topologie ist, so gilt f"ur die Wertemenge des Operators $\mathbf{cl}_{\mathbf{T}}$ offensichtlich
\begin{equation}\label{iwwed}
\mathbf{P}_{2}\mathbf{cl}_{\mathbf{T}}\subset  \mathbf{T}^{c}\ .
\end{equation}
F"ur alle Mengensysteme $\mathcal{A}$ ist die Monotonie 
deren jeweiliger H"ullenoperatoren $\mathbf{cl}_{\mathcal{A}}$ in Form der
f"ur alle $A,B\in \subset\bigcup \mathcal{A}$ g"ultigen Aussage
\begin{equation}
\begin{array}{c}
\mathbf{cl}_{\mathcal{A}}(A)\supset A\ \land\\
A\subset B\Rightarrow \mathbf{cl}_{\mathcal{A}}(A)\supset\mathbf{cl}_{\mathcal{A}}(B)
\end{array}
\end{equation}
gegeben.
Wegen dieser Monotonie
und wegen (\ref{idwed}) und (\ref{iwwed})
ist
f"ur alle Teilmengen $A,B\subset\bigcup \mathbf{T}$ die Aussage
\begin{displaymath}
\begin{array}{c}
\underbrace{\mathbf{cl}_{\mathbf{T}}(A)\cup\mathbf{cl}_{\mathbf{T}}(B)}_{\in \mathbf{T}^{c}}\supset A\cup B\\
\land\\
\mathbf{cl}_{\mathbf{T}}(A)\cup\mathbf{cl}_{\mathbf{T}}(B)\\
=\mathbf{cl}_{\mathbf{T}}(\mathbf{cl}_{\mathbf{T}}(A)\cup\mathbf{cl}_{\mathbf{T}}(B))\supset
\mathbf{cl}_{\mathbf{T}}(\mathbf{cl}_{\mathbf{T}}(A))\cup\mathbf{cl}_{\mathbf{T}}(\mathbf{cl}_{\mathbf{T}}(B))\\
=\mathbf{cl}_{\mathbf{T}}(A)\cup\mathbf{cl}_{\mathbf{T}}(B)
\end{array}
\end{displaymath}
wahr.
Wiederum
wegen (\ref{idwed}) gilt also die Inklusionenkette
\begin{displaymath}
\begin{array}{c}
\mathbf{cl}_{\mathbf{T}}(A)\cup\mathbf{cl}_{\mathbf{T}}(B)\\
\subset\mathbf{cl}_{\mathbf{T}}(A)\cup\mathbf{cl}_{\mathbf{T}}(B)\\ \subset
\mathbf{cl}_{\mathbf{T}}(A)\cup\mathbf{cl}_{\mathbf{T}}(B)\ ,
\end{array}
\end{displaymath}
also die Gleichung (\ref{wwwed}) 
\begin{displaymath}
\mathbf{cl}_{\mathbf{T}}(A\cup B)=\mathbf{cl}_{\mathbf{T}}(A)\cup\mathbf{cl}_{\mathbf{T}}(B)\ ,
\end{displaymath}
deren G"ultigkeit wir bereits im Abschnitt (\ref{imit}) behaupteten. 
Ebenda
wiesen wir auch
schon darauf hin,
dass es
offensichtliche Gegenbeispiele selbstdualer Topologien gibt, die 
die G"ultigkeit der
zu der (\ref{wwwed})
entsprechenden Gleichung, in der statt \glqq$\cup$\grqq das Zeichen
\glqq$\cap$\grqq steht, widerlegen.
\newline
Als wichtiger als der Trost, dass der H"ullenoperator eines Mengensystemes $\mathcal{A}$
immerhin idempotent ist, erweist sich f"ur unsere Belange, dass
die folgende Aussage gilt, die wir im Beweis 
sowohl 
des verallgemeinerten insensitiven Ergodensatzes 3.3 
als auch der Bemerkung 3.4 
anwenden:
\newline
\newline
{\bf Lemma 3.1: Reichhaltigkeitslemma}\index{Reichhaltigkeitslemma}\newline
{\em F"ur jedes Mengensystem $\mathcal{A}$ und f"ur jede Menge $B\subset\bigcup \mathcal{A}$ gilt die Implikation}
\begin{equation}\label{loee}
X\in(\mathcal{A}\cap\mathbf{cl}_{\mathcal{A}}(B))\setminus\{\emptyset\}\Rightarrow B\cap X\not=\emptyset\ .
\end{equation}
{\bf Beweis:}\newline
Wenn $X\in \mathcal{A}\cap\mathbf{cl}_{\mathcal{A}}(B)\setminus\{\emptyset\}$ ein nicht leeres Element
des Spurmengensystemes $\mathcal{A}\cap\mathbf{cl}_{\mathcal{A}}(B)$ ist, gibt es eine Menge $A_{0}\in \mathcal{A}$, f"ur die
$$X=A_{0}\cap\mathbf{cl}_{\mathcal{A}}(B)$$
$$=A_{0}\cap\bigcap\{A\in\mathcal{A}^{c}:A\supset B\}$$
und
$$ B\cap X= B\cap A_{0}$$
ist. W"are $B\cap A_{0}$ leer, so w"are das zu der Menge $A_{0}$ disjunkte Komplement
$A_{0}^{c}\in \mathcal{A}^{c}$ in dem
Mengensystem 
$$\{A\in\mathcal{A}^{c}:A\supset B\}\ni A_{0}^{c}\ ,$$
sodass
$$\mathbf{cl}_{\mathcal{A}}(B)=\bigcap\{A\in\mathcal{A}^{c}:A\supset B\}\cap A_{0}\ =\emptyset$$
g"alte, im Widerspruch zu der 
Vorraussetzung, dass $X$
im 
Spurmengensystem
$(\mathcal{A}\cap\mathbf{cl}_{\mathcal{A}}(B))\setminus\{\emptyset\}$ ist.
\newline
{\bf q.e.d.}
\newline\newline
Welch zentrale Rolle 
diese Reichhaltigkeitsaussage (\ref{loee}) 
beim Beweis 
des verallgemeinerten insensitiven Ergodensatzes 3.3 spielt, werden wir sehen.
Deshalb werten wir die Implikation (\ref{loee}) 
als ein Lemma. Das Reichhaltigkeitslemma ist dabei dadurch universalisierbar,\index{Reichhaltigkeitslemma}
dass wir 
f"ur jedes Mengensystem $\mathcal{A}$ 
deren H"ullenoperator $\mathbf{cl}_{\mathcal{A}}$ universalisieren, der f"ur 
jede Menge $B$
den Wert $$\mathbf{cl}_{\mathcal{A}}(B):=\bigcap\{A\in\mathcal{A}^{c}:A\supset B\}$$
haben soll. F"ur 
jedes Mengensystem $\mathcal{A}$ und
deren jeweiligen universalisierten
H"ullenoperator $\mathbf{cl}_{\mathcal{A}}$,
der nun auf der Klasse der Mengen definiert ist, 
gilt dann f"ur jede Menge $B$ ebenfalls\index{universalisierter H\"ullenoperator eines Mengensystemes}
die Implikation (\ref{loee}).
\newline
Ferner sei f"ur jedes Paar zweier Abbildungen
$q,p\in X^{X}$ einer Potenzmenge $X=2^{Q}$ einer Menge $Q$ oder einer endlichen kartesischen Potenz
$$X=2^{Q}\times 2^{Q}\dots \times 2^{Q}$$ der Potenzmenge $2^{Q}$ einer Menge $Q$ in diese
Menge $X$ der Kommutator
\begin{equation}\label{cloeee}
[p,q]:=p\Delta q
\end{equation}
als mengenwertige Abbildung der Potenzmenge $2^{Q}$ in dieselbe festgelegt,
wobei f"ur je zwei endliche $n$-Tupel von Mengen $(A_{1},A_{2},\dots A_{n}),
(B_{1},B_{2},\dots B_{n})$
$$(A_{1},A_{2},\dots A_{n})\Delta(B_{1},B_{2},\dots B_{n})$$
\begin{equation}
:=(A_{1}\Delta B_{1})\Delta(A_{2}\Delta B_{2})\dots\Delta(A_{n}\Delta B_{n})
\end{equation}
die symmetrische Differenz\index{symmetrische Differenz zweier Mengen}
dieser beiden endlichen $n$-Tupel von lauter Mengen $(A_{1},A_{2},\dots A_{n})$ und 
$(B_{1},B_{2},\dots B_{n})$
notiere.\footnote{Dies scheint die 
Festlegung einer Schreibweise zu sein, die ausschliesslich eine suggestive Funktion hat.
Die Schreibweise (\ref{cloeee}) hat aber auch schlicht einen schreibtechnischen Zweck, weil sie f"ur die
Handhabung von
Abbildungen gemacht ist, w"ahrend sich die symmetrische Differenz $\Delta$ der
Darstellung der Verkn"upfung zweier Mengen anpasst.}
Demnach k"onnen wir die Gleichung (\ref{wwwed})
durch die Kommutation
\begin{equation}
\Bigl[\mathbf{cl}_{\mathbf{T}}\mathbf{P}_{1}\oplus\mathbf{cl}_{\mathbf{T}}\mathbf{P}_{2},\ (\mathbf{P}_{1}\cup\mathbf{P}_{2})\oplus (\mathbf{P}_{1}\cup\mathbf{P}_{2})\Bigr]=\emptyset
\end{equation}
auf zwar h"ochst unelegante Weise, jedoch unmissverst"andlich paraphrasieren; und nichtsdestotrotz
sehen wir hierbei auch, dass das Abschlussverhalten von Mengensystemen gegebenfalls in Form von
Kommutationen ausgedr"uckt werden kann.
Es gilt offenbar die folgende Verallgemeinerung des insensitiven topologischen 
Ergodensatzes\index{Abschlusskommutation}\index{insensitiver topologischer Ergodensatz} 2.2:
\newline
\newline
{\bf Bemerkung 3.2: Verallgemeinerte kollektive Abschlusstransitivit"at}\newline
{\em Es sei $\mathcal{A}$ ein Mengensystem "uber dem Zustandsraum $\mathbf{P}_{2}\Phi$ einer Flussfunktion $\Phi$,
das denselben "uberdeckt und 
$\Phi$ eine Flussfunktion, all deren 
Fl"usse $\Phi^{t}\in \Gamma_{\Phi}$ mit dem 
gem"ass} (\ref{cloee}) {\em festgelegten
Operator $\mathbf{cl}_{\mathcal{A}}$ 
in dem Sinn kommutieren, dass
\begin{equation}\label{komma}
[\mathbf{cl}_{\mathcal{A}},\Phi^{t}]=\emptyset
\end{equation}
ist.
Dann gilt
\begin{equation}
\mathbf{cl}_{\mathcal{A}}({\rm P})\in 2^{\widehat{\mathbf{T}}(\Gamma_{\Phi})}\cap\mathbf{part}(\mathbf{P}_{2}\Phi)
\end{equation}
f"ur jede Partition}
\begin{equation}\label{trzeu}
{\rm P}\in 2^{\widehat{\mathbf{T}}(\Gamma_{\Phi})}\cap\mathbf{part}(\mathbf{P}_{2}\Phi)\ .
\end{equation}
\newline
{\bf Beweis:}\newline
Wir erl"autern, dass
$\mathbf{cl}_{\mathcal{A}}({\rm P})$ als das Mengensystem $\{\mathbf{cl}_{\mathcal{A}}(p):p\in {\rm P}\}$
zu lesen
ist. F"ur alle Elemente $p$ der Partition ${\rm P}$ und alle $t\in \mathbb{R}$ gilt analog zur der Argumentation
des Beweises des insensitiven topologischen Ergodensatzes wegen der Kommutativit"at (\ref{komma})
$$\Phi^{t}(\mathbf{cl}_{\mathcal{A}}(p))=\mathbf{cl}_{\mathcal{A}}(p)\ ,$$
sodass das Mengensystem 
$\mathbf{cl}_{\mathcal{A}}({\rm P})$ lauter flussinvariante Mengen als Elemente hat und daher
der Inklusion
$$\mathbf{cl}_{\mathcal{A}}({\rm P})\subset\mathbf{T}(\Gamma_{\Phi})$$ 
gen"ugt; deswegen
sind je zwei verschiedene Mengen des Mengensystemes
$\mathbf{cl}_{\mathcal{A}}({\rm P})$ zueinander disjunkt und wegen
der in
(\ref{trzeu}) vorausgesetzten
Partitivit"at des Mengensystemes ${\rm P}$
ist $\mathbf{cl}_{\mathcal{A}}({\rm P})$
eine Zustandsraumpartition.\newline
{\bf q.e.d.}
\newline
\newline
F"ur jede Kollektivierung $[\Phi]$ einer Flussfunktion $\Phi$ 
und zu jedem "uber dem Zustandsraum $\mathbf{P}_{2}\Phi$ 
gegebenen, denselben "uberdeckenden
Mengensystem $\mathcal{A}$
sei 
das Mengensystem
\begin{equation}\label{tzeu}
[[\Phi]]_{\mathcal{A}}:=\{\mathbf{cl}_{\mathcal{A}}(\gamma):\gamma\in [\Phi]\}
\end{equation}
notiert. Wir stossen auf folgenden
\newline
\newline
{\bf Satz 3.3: Verallgemeinerter insensitiver Ergodensatz}\index{verallgemeinerter insensitiver Ergodensatz}\newline
{\em Es sei $\mathcal{A}$ ein Mengensystem "uber dem Zustandsraum $\mathbf{P}_{2}\Phi$ einer Flussfunktion $\Phi$,
das denselben "uberdeckt und 
$\Phi$ eine Flussfunktion, all deren 
Fl"usse $\Phi^{t}\in \Gamma_{\Phi}$ mit dem Operator $\mathbf{cl}_{\mathcal{A}}$ gem"ass}
(\ref{cloee}) {\em in dem Sinn kommutieren, dass
$$[\mathbf{cl}_{\mathcal{A}},\Phi^{t}]=\emptyset$$
ist.
Dann gilt f"ur das gem"ass} (\ref{ienatt}) {\em und} (\ref{idkoh}) und (\ref{cloocee}) {\em verfasste Mengensystem $\mbox{{\bf @}}(\Phi,\overline{\mathcal{A}})$ freier 
Attraktoren\index{freier Attraktor bez\"uglich eines Mengensystemes} bez"uglich der Zustandsraum"uberdeckung $\mathcal{A}$ und
f"ur das gem"ass} (\ref{dikos}) {\em verfasste Mengensystem topologischer Attraktoren:}\index{topologischer Attraktor}
\begin{equation}
\begin{array}{c}
\mbox{{\bf @}}(\Phi,\overline{\mathcal{A}})=
\underline{\mbox{{\bf @}}}(\Phi,\overline{\mathcal{A}},\mathbf{cl}_{\mathcal{A}}(2^{\mathbf{P}_{2}\Phi}))\supset\\
\lbrack\lbrack\Phi\rbrack\rbrack_{\mathcal{A}}\in\mathbf{part}(\mathbf{P}_{2}\Phi)\ .
\end{array}
\end{equation}
\newline
{\bf Beweis:}\newline 
Wegen der soeben gezeigten verallgemeinerten kollektiven Abschlusstransitivit"at
ist f"ur die bez"uglich des Mengensystemes $\mathcal{A}$ kommutative 
Flussfunktion $\Phi$ mit deren Kollektivierung $[\Phi]$ auch das gem"ass (\ref{tzeu}) festgelegte 
Mengensystem 
$$[[\Phi]]_{\mathcal{A}}\in\mathbf{part}(\mathbf{P}_{2}\Phi)$$
eine Partition des Zustandsraums $\mathbf{P}_{2}\Phi$. Das Reichhaltigkeitslemma, die
Implikation (\ref{loee}) also, spezifiziert sich hier f"ur alle Kollektivelemente $\gamma\in [\Phi]$ zu der
Implikation
\begin{equation}\label{amml}
X\in\mathcal{A}\cap(\mathbf{cl}_{\mathcal{A}}(\gamma))\setminus\{\emptyset\}\Rightarrow \gamma\cap X\not=\emptyset\ ,
\end{equation}
die wahr ist. Also gilt f"ur alle 
$$\mathbf{cl}_{\mathcal{A}}(\gamma)\in \mathbf{cl}_{\mathcal{A}}(\mathbf{T}(\Gamma_{\Phi}))=[[\Phi]]_{\mathcal{A}}$$
die nicht nur die im Hinblick auf die Festlegung 
durch die Implikation
(\ref{idkoh}) angepasste "Aquivalenz
\begin{displaymath}
\begin{array}{c}
A,B\subset\mathbf{cl}_{\mathcal{A}}(\gamma)\ \land\ \emptyset\not\in\{A,B\}\ \land\ A,B\in\mathcal{A}\cap
\mathbf{cl}_{\mathcal{A}}(\gamma)\\ 
\Leftrightarrow\\
A,B\in(\mathcal{A}\cap\mathbf{cl}_{\mathcal{A}}(\gamma))\setminus\{\emptyset\}\ .
\end{array}
\end{displaymath}
Wegen der 
Implikation (\ref{amml}) gilt ausserdem
auch die Spezifizierung des Reichhaltigkeitslemma\index{Reichhaltigkeitslemma} in Form der
Implikation
\begin{equation}
\begin{array}{c}
A,B\subset\mathbf{cl}_{\mathcal{A}}(\gamma)\ \land\ \emptyset\not\in\{A,B\}\ \land\ A,B\in\mathcal{A}\cap
\mathbf{cl}_{\mathcal{A}}(\gamma)\Rightarrow\\
A\cap\gamma\not=\emptyset\not=\gamma\cap B\ \ ,
\end{array}
\end{equation}
wobei $\gamma\in [\Phi]\subset \mathbf{T}(\Gamma_{\Phi})$ ist. Daher beinhaltet die
Aussage, dass f"ur zwei Mengen
$A$ und $B$ die Differenz
$$A\cap\gamma\not=\emptyset\not=\gamma\cap B$$ 
gegeben ist, dass es einen Zustand $x$ gibt, f"ur den die Differenz
$$A\cap\Phi(x,\mathbb{R})\not=\emptyset\not=\Phi(x,\mathbb{R})\cap B$$
vorliegt; sodass es nach dem Lemma 1.3 eine reelle Zahl $t\in\mathbb{R}$ gibt, f"ur welche 
schliesslich die Reichhaltigkeit
$$\Phi(\Phi(x,\mathbb{R})\cap B,t)\cap(A\cap\Phi(x,\mathbb{R}))\not=\emptyset$$
$$\not=\Phi(\Phi(x,\mathbb{R})\cap B,t)\cap A$$
gegeben ist, die sich in der von ihr implizierten Reichhaltigkeit
$$\Phi( B,t)\cap A\not=\emptyset$$
als die Gegebenheit erweist, aus der
sich an der Festlegung 
des Mengensystemes $\mbox{{\bf @}}(\Phi,\overline{\mathcal{A}})$
gem"ass
(\ref{idkoh}) die Behauptung ablesen l"asst.\newline
{\bf q.e.d.}
\newline
\newline
Eine unmittelbare Konsequenz des Satzes 3.3
ist die allgemeine Kovarianz der Attraktoren:
Dann, wenn f"ur eine 
auf dem
Zustandsraum $\mathbf{P}_{2}\Phi$ einer Flussfunktion\index{Kovarianz der Attraktoren}
$\Phi$ definierte Bijektion $f$ die
Kommutatorrelation
\begin{equation}
[\mathbf{cl}_{f(\mathcal{A})},\Phi_{f}^{t}]=\emptyset
\end{equation} f"ur alle $t\in\mathbb{R}$ 
erf"ullt ist, gilt f"ur alle Teilmengen $\chi\subset \mathbf{P}_{2}\Phi$ des Zustandsraumes die Implikation
\begin{equation}
\begin{array}{c}
\chi\in\underline{\mbox{{\bf @}}}(\Phi,\overline{\mathcal{A}},\mathbf{cl}_{\mathcal{A}}(2^{\mathbf{P}_{2}\Phi}))\\
\Rightarrow\\
f(\chi)\in\underline{\mbox{{\bf @}}}(\Phi_{f},\overline{f(\mathcal{A})},\mathbf{cl}_{f(\mathcal{A})}(2^{f(\mathbf{P}_{2}\Phi)}))\  .
\end{array}
\end{equation}
Dabei ist
$$f(\mathcal{A})=\{\{f(a):a\in A\}:A\in\mathcal{A}\}\ ,$$
$$\overline{f(\mathcal{A})}=\{\mathbf{cl}_{f(\mathcal{A})}(y):y\inf(\mathcal{A})\ .$$
Der Satz 3.3 formuliert ein hinreichendes, jedoch kein notwendiges Kriterium
f"ur die Attraktorgenerativit"at eines Flusses:\index{Attraktorgenerativit\"at eines Flusses}
Genau dann, wenn f"ur einen Fluss $\{\Phi^{t}\}_{t\in \mathbb{R}}$ und f"ur ein
dessen Zustandsraum $\mathbf{P}_{2}\Phi$ "uberdeckendes Mengensystem
$\mathcal{A}$ "uber $\mathbf{P}_{2}\Phi$ die Inklusion
\begin{equation}\label{qqert}
\mbox{{\bf @}}(\Phi,\overline{\mathcal{A}})=\underline{\mbox{{\bf @}}}(\Phi,\overline{\mathcal{A}},\mathbf{cl}_{\mathcal{A}}(2^{\mathbf{P}_{2}\Phi}))\supset[[\Phi]]_{\mathcal{A}}
\end{equation}
wahr 
ist, nennen wir den Fluss $\{\Phi^{t}\}_{t\in \mathbb{R}}$
bez"uglich $\mathcal{A}$ attraktogen.\index{attraktogener Fluss}\index{kommutativer Fluss}
Genau dann hingegen, wenn f"ur einen Fluss $\{\Phi^{t}\}_{t\in \mathbb{R}}$ und
ein Mengensystem $\mathcal{A}$ "uber dem Zustandsraum $\mathbf{P}_{2}\Phi$ desselben
die Implikation (\ref{qqert}) 
\begin{equation}\label{qqert}
t\in\mathbb{R}\Rightarrow 
[\mathbf{cl}_{\mathcal{A}},\Phi^{t}]=\emptyset
\end{equation}
gilt, nennen wir $\{\Phi^{t}\}_{t\in \mathbb{R}}$ bez"uglich $\mathcal{A}$ kommutativ. 
\newline
Zun"achst wollen wir uns
die Kommutativit"at eines Flusses bez"uglich eines Mengensystemes
erst einmal als Negativ vor Augen halten: Genau dann, wenn sie nicht vorliegt,
gibt es 
eine Teilmenge $X\subset\mathbf{P}_{2}\Phi$
des
Zustandsraumes und eine reelle Zahl $t$, f"ur die
\begin{equation}\label{qert}
\Phi^{t}(\mathbf{cl}_{\mathcal{A}}(X))\setminus\mathbf{cl}_{\mathcal{A}}(\Phi^{t}(X))\not=\emptyset 
\end{equation}
ist.
Denn, dass 
$$\emptyset\not=[\mathbf{cl}_{\mathcal{A}},\Phi^{t}]=\mathbf{cl}_{\mathcal{A}}(\Phi^{t}(X))
\Delta\Phi^{t}(\mathbf{cl}_{\mathcal{A}}(X))$$
ist, impliziert,
die Reichhaltigkeit (\ref{qert}) oder, dass 
$$\mathbf{cl}_{\mathcal{A}}(\Phi^{t}(X))\setminus\Phi^{t}(\mathbf{cl}_{\mathcal{A}}(X))\not=\emptyset$$
ist. Letzteres gilt genau dann, wenn es 
eine Teilmenge $X\subset\mathbf{P}_{2}\Phi$
des
Zustandsraumes und eine reelle Zahl $t$ gibt, f"ur die
$$\mathbf{cl}_{\mathcal{A}}(\Phi^{t}(X))\setminus
\Phi^{t}(\mathbf{cl}_{\mathcal{A}}((\Phi^{t})^{-1}(\Phi^{t}(X)))\not=\emptyset$$ 
ist, sodass es genau dann die 
Teilmenge
$Y=\Phi^{t}(X)$ und die Zahl $\tilde{t}\in\mathbb{R}$ 
mit $\Phi^{\tilde{t}}:=(\Phi^{t})^{-1}$ gibt, f"ur die
$$\Phi^{\tilde{t}}(\mathbf{cl}_{\mathcal{A}}(Y))\setminus
\mathbf{cl}_{\mathcal{A}}(\Phi^{\tilde{t}}(Y))\not=\emptyset$$
gilt. Die Zahl $\tilde{t}\in\mathbb{R}$ gibt es auch, wenn der 
Fluss $\{\Phi^{t}\}_{t\in \mathbb{R}}$ nicht 
phasisch ist: 
Auch, wenn die Menge der Autobolismen $\{\Phi^{t}:t\in \mathbb{R}\}$ nicht algebraisch ist,
gibt es f"ur jede reelle Zahl $t$ die Inversion $(\Phi^{t})^{-1}\in \{\Phi^{t}:t\in \mathbb{R}\}$.
Demnach gilt die Implikation (\ref{qqert}) auch genau dann nicht, wenn es 
eine Teilmenge $X\subset\mathbf{P}_{2}\Phi$
des
Zustandsraumes und eine reelle Zahl $t$ gibt, f"ur die
\begin{equation}\label{qqqert}
\mathbf{cl}_{\mathcal{A}}(\Phi^{t}(X))\setminus\Phi^{t}(\mathbf{cl}_{\mathcal{A}}(X))\not=\emptyset
\end{equation}
ist. Wenn die Kommutativit"at (\ref{qqert}) erf"ullt ist,
und es die verallgemeinerten Zimmer $\chi\in [[\Phi]]_{\mathcal{A}}$ als
Attraktoren der Menge $\mbox{{\bf @}}(\Phi,\overline{\mathcal{A}})$
gibt, so sind dieselben Invariante auch jeden abgewandelten Flusses $\{\Phi^{t}_{\beta}\}_{t\in \mathbb{R}}$, der sich 
folgendermassen aus dem Fluss $\{\Phi^{t}\}_{t\in \mathbb{R}}$ ergibt:
\newline
Es gibt zu jedem Paar $(\{\Psi^{t}\}_{t\in \mathbb{R}}, \mathcal{Z})$,
das aus einem Fluss $\{\Psi^{t}\}_{t\in \mathbb{R}}$
und aus einem  
"uberdeckenden
Mengensystemes $\mathcal{Z}$ "uber
dessen Zustandsraum besteht,
die Menge aller seiner 
Flussvertauschungsfelder.
Exakt jedes reellwertige auf der Definitionsmenge $\mathbf{P}_{1}\Psi$ der durch den Fluss 
$\{\Psi^{t}\}_{t\in \mathbb{R}}$ festgelegten Flussfunktion definierte Feld 
\begin{equation}\label{qqqerq}
\begin{array}{c}
\beta:\mathbf{P}_{1}\Psi\to\mathbb{R}\ ,\\
(x,t)\mapsto\beta(x,t)
\end{array} 
\end{equation}
bezeichnen wir als ein Flussvertauschungsfeld des Paares $(\{\Psi^{t}\}_{t\in \mathbb{R}},\mathcal{Z})$,
falls es
so beschaffen ist, dass f"ur alle Zust"ande $x\in \mathbf{P}_{2}\Psi$ die
Invertierbarkeit
\begin{equation}
\beta(x,\mbox{id})^{-1}\in \mathbb{R}^{\mathbb{R}}
\end{equation}
gegeben und "uberdies die
die "Aquivalenz
\begin{equation}
\beta(x,t)=\beta(\tilde{x},t)\Leftrightarrow \tilde{x}\in\mathbf{cl}_{\mathcal{Z}}(\Psi(x,\mathbb{R}))
\end{equation} 
f"ur alle $t\in\mathbb{R}$
g"ultg ist. 
Es sei $\beta$ ein Flussvertauschungsfeld des Paares $(\{\Phi^{t}\}_{t\in \mathbb{R}},\mathcal{A})$.
Dann ist der Fluss $\{\Phi^{t}_{\beta}\}_{t\in \mathbb{R}}$, f"ur den
f"ur alle $t\in\mathbb{R}$
\begin{equation}\label{qqqerq}
\begin{array}{c}
\Phi^{t}_{\beta}:\mathbf{P}_{2}\Phi\to\mathbf{P}_{2}\Phi\ ,\\
x\mapsto\Phi^{t}_{\beta}(x):=\Phi^{\beta(x,t)}(x)
\end{array} 
\end{equation}
ist,
von der Art, dass die Invarianz der verallgemeinerten Zimmer in Form der 
G"ultigkeit der Implikation
\begin{equation}
(\chi,t)\in [[\Phi]]_{\mathcal{A}}\times\mathbb{R} \Rightarrow \Phi^{t}_{\beta}\chi=\chi
\end{equation}
vorliegt, wobei die Invarianz der Kollektivelemente $\gamma\in[\Phi]$ gegen"uber
den Autobolismen des Flusses $\{\Phi^{t}_{\beta}\}_{t\in \mathbb{R}}$
im allgemeinen 
keineswegs gegeben ist. Die verallgemeinerten Zimmer $\chi\in[[\Phi]]_{\mathcal{A}}$
stellen aber dennoch offensichtlich die abgeschlossenen Attaktoren des Flusses $\{\Phi^{t}_{\beta}\}_{t\in \mathbb{R}}$
bez"uglich des Mengensystemes $\mathcal{A}$ dar, die
per constructionem invariant sind gegen 
jeweils f"ur einzelne Zimmer gegebene Indizierungspermutation des
Flusses $\{\Phi^{t}\}_{t\in \mathbb{R}}$. 
\newline
Wir nennen zwei Fl"usse $\{\Psi^{t}\}_{t\in \mathbb{R}}$ und $\{\Xi^{t}\}_{t\in \mathbb{R}}$
genau dann bez"uglich eines Mengensystemes $\mathcal{Z}$ "aquivalent, wenn
sie den gleichen Zustandsraum haben und wenn $\mathcal{Z}$ ein
denselben "uberdeckendes Mengensystemes $\mathcal{Z}$ "uber diesem ist, f"ur das
\begin{equation}
[[\Psi]]_{\mathcal{Z}}=[[\Xi]]_{\mathcal{Z}}
\end{equation}
gilt; was keineswegs heisst, dass dabei die bez"uglich $\mathcal{Z}$
abgeschlossenen Mengen des Mengensystemes $[[\Psi]]_{\mathcal{Z}}$ 
Attraktoren sein m"ussen. Offensichtlich
sind zwei Fl"usse $\{\Psi^{t}\}_{t\in \mathbb{R}}$ und $\{\Xi^{t}\}_{t\in \mathbb{R}}$
genau dann zueinander bez"uglich eines Mengensystemes $\mathcal{Z}$ "aquivalent, wenn
es zwei Flussvertauschungsfelder $\beta$ bzw. $\bar{\beta}$
des Paares $(\{\Psi^{t}\}_{t\in \mathbb{R}},\mathcal{Z})$ bzw. $(\{\Xi^{t}\}_{t\in \mathbb{R}},\mathcal{Z})$
gibt, f"ur die
\begin{equation}
\begin{array}{c}
\{\Xi^{t}_{\beta}\}_{t\in \mathbb{R}}=\{\Psi^{t}\}_{t\in \mathbb{R}}\ ,\\ 
\{\Psi^{t}_{\bar{\beta}}\}_{t\in \mathbb{R}}=\{\Xi^{t}\}_{t\in \mathbb{R}}
\end{array}
\end{equation}
ist. Es erhebt sich damit die
Frage, ob
ein Fluss $\{\Psi^{t}\}_{t\in \mathbb{R}}$ bez"uglich
eines seinen Zustandsraum
"uberdeckenden Mengensystemes $\mathcal{Z}$ genau dann
attraktogen ist, wenn er bez"uglich desselben zu einem
Fluss $\{\Phi^{t}\}_{t\in \mathbb{R}}$
"aquivalent ist, der 
bez"uglich $\mathcal{Z}$
kommutativ ist. Der Antwort auf diese Frage bringt die folgende Bemerkung n"aher:
\newline
\newline
{\bf Bemerkung 3.4: "Aquivalente Attraktorgenerativit"at}\newline
{\em Jeder Fluss $\{\Psi^{t}\}_{t\in \mathbb{R}}$ ist bez"uglich
eines seinen Zustandsraum $\mathbf{P}_{2}\Psi$
"uberdeckenden Mengensystemes $\mathcal{Z}$ "uber $\mathbf{P}_{2}\Psi$
genau dann
attraktogen, wenn seine Zimmer 
bez"uglich des Mengensystemes $\mathcal{Z}$
flussinvariant sind: Es gilt die "Aquivalenz}
\begin{equation}
\begin{array}{c}
((\chi,t)\in [[\Psi]]_{\mathcal{Z}}\times\mathbb{R} \Rightarrow \Psi^{t}\chi=\chi)\\
\Leftrightarrow\\
\lbrack\lbrack\Psi\rbrack\rbrack_{\mathcal{Z}}\subset
\mbox{{\bf @}}(\Psi,\overline{\mathcal{Z}})\ .
\end{array}
\end{equation}
\newline
\newline
{\bf Beweis:}\newline
Die Flussinvarianz ist eine Eigenschaft, die jeder freie Attraktor per constructionem hat.
Wir brauchen lediglich zu zeigen, dass f"ur alle Zust"ande $x$ des Zustandsraumes 
$\mathbf{P}_{2}\Psi$
die Implikation
$$A,B\in (\mathcal{Z}\cap\mathbf{cl}_{\mathcal{Z}}(\Psi(x,\mathbb{R})))\setminus\{\emptyset\}\Rightarrow
\exists\ t\in \mathbb{R}:$$
$$\Psi^{t}(A)\cap B\not=\emptyset $$
wahr ist: Falls $A\cap\Psi(x,\mathbb{R})$ 
oder $B\cap\Psi(x,\mathbb{R})$
leer w"are, so k"onnte es sein, dass diese Implikation nicht gilt.
Andernfalls schliesst das Lemma 1.3 diese M"oglichkeit aus.
Wegen des Reichhaltigkeitslemmas allerdings gilt
$$A,B\in (\mathcal{Z}\cap\mathbf{cl}_{\mathcal{Z}}(\Psi(x,\mathbb{R})))\setminus\{\emptyset\}\Rightarrow$$
$$A\cap\Psi(x,\mathbb{R})\not=\emptyset\not=B\cap\Psi(x,\mathbb{R})\ .$$
{\bf q.e.d.}
\newline
\newline
Daher gilt
$$[[\Psi]]_{\mathcal{Z}}\subset\mbox{{\bf @}}(\Psi,\overline{\mathcal{Z}})\subset[[\Psi]]_{\mathcal{Z}}^{\cup}\ ,$$
denn die Mengen $\chi\in [[\Psi]]_{\mathcal{Z}}$ sind die
kleinsten flussinvarianten, bez"uglich des Mengensystemes $\mathcal{Z}$ abgeschlossenen
Attraktoren.\newline
Die folgende, letzte Definition des ersten Teiles der abstrakten Ergodentheorie
ist die entsprechend allgemeine Abstraktion des Begriffes des Zimmers, die
den beschriebenen, f"ur den Attraktorbegriff vollzogenen, Verallgemeinerungen korrespondiert: 
\newline
\newline
{\bf Definition 3.5:\newline Allgemeine Zimmer und Vorzimmer von Wellenfunktionen}\newline
{\em Es seien $X$ und $Y$ zwei Mengen, sodass $$\Xi\in Y^{X\times\mathbb{R}}$$ 
eine Wellenfunktion ist und 
$$\mathcal{Z}\subset 2^{\mathbf{P}_{2}\Xi}\subset2^{Y}$$
ein Mengensystem, das die Wertemenge $\mathbf{P}_{2}\Xi=\bigcup\mathcal{Z}$ der Wellenfunktion 
$\Xi$
"uberdeckt. Exakt die Mengen des Mengensystemes
\begin{equation}\label{koren}
\lbrack\lbrack\Xi\rbrack\rbrack_{\mathcal{Z}}:=\{\mathbf{cl}_{\mathcal{Z}}(\Xi(x,\mathbb{R})):x\in X\}
\end{equation}
nennen wir die Vorzimmer der Wellenfunktion $\Xi$ bez"uglich
des Mengensystemes $\mathcal{Z}$. 
Genau dann, wenn das Mengensystem der Vorzimmer der Wellenfunktion $\Xi$ bez"uglich
des Mengensystemes $\mathcal{Z}$ eine Partition\index{Vorzimmer einer Wellenfunktion bez\"uglich eines
Mengensystemes}
\begin{equation}\label{kooren}
\lbrack\lbrack\Xi\rbrack\rbrack_{\mathcal{Z}}\in\mathbf{part}(\mathbf{P}_{2}\Xi)
\end{equation}
der Wertemenge $\mathbf{P}_{2}\Xi$ der Wellenfunktion $\Xi$
ist, bezeichnen wir 
dessen Elemente als die Zimmer der Wellenfunktion $\Xi$
bez"uglich
des Mengensystemes $\mathcal{Z}$.}\index{Zimmer einer Wellenfunktion bez\"uglich eines
Mengensystemes}
\newline
\newline
W"ahrend die nicht trivialen Zimmer stetiger reeller 
Flussfunktionen $\Psi$ mit kompaktem Zustandsraum $\zeta\subset\mathbb{R}^{n}$
sowohl
mit deren Attraktoren
als auch mit deren sensitiven Attraktoren "uberstimmen,
geht die begriffliche Differenz zwischen Zimmern und 
Vorzimmern in der vollen Generalit"at auf.\newline\newline
{\small Wir wiesen 
bei der Einf"uhrung der Zimmer reeller Zustandsr"aume
in der Abhandlung des elementaren Quasiergodensatzes \cite{raab}
auf eine durchaus um Effekt bem"uhte Weise
auf Leibnizens Monadologie\index{Monadologie} \cite{mona} hin und wir brachten
Leibnizens Monaden mit den Zimmern in Zusammenhang.
Der Sachverhalt, dass es f"ur endliche
und einfache Zustandsraum"uberdeckungen auf eine sehr 
durchsichtige Weise sehr
schlicht beschaffene Zimmer gibt, mag nun den
Eindruck erwecken, dass solche schlichte Zimmer kaum den 
Leibnizschen\index{Leibniz, Gottfried Wilhelm} Monaden entsprechen k"onnen. Leibnizens Monaden liegen jenseits der Mathematik, jedoch innerhalb einer Ontologie,
welche Sein, Reales, meint und die insofern nicht als von der physikalischen Realit"at 
abgetrennt gedacht ist. 
Leibnizens Monadologie spricht also insofern in der Tat nicht von den logisch reduziblen Zimmern, die
wir hier diskutieren. Unseren Zimmern fehlt bislang der Bezug zur physikalischen Zeit, die hier
nur schattenhaft als Parametrisierung des Flusses auftaucht.
Erst wenn es gel"ange, die physikalische Zeit im Modell zu konkretisieren und das, was
sie zur physikalischen Echtzeit macht, im Modell zu erfassen,
k"onnen wir erhoffen, dass wir die Monaden der durchaus dunklen Leibnizschen Schau
treffen.}\newline
\newline
Die Konsequenzen des verallgemeinerten insensitiven Ergodensatzes 3.3 und der Bemerkung 3.4 lassen sich im Rahmen dieses Traktates 
nicht darlegen. Der abstrakte, vom Generalisierungsangebot der Abschlusskommutation\index{Abschlusskommutation} 
induzierte
Gedanke, die Kontinuit"at   
einer Flussfunktion $\Phi$ mittels der Kommutatorrelation f"ur ihre Fl"usse zu
generalisieren, ist kein abstraktes Hirngespinst. Manche der dynamischen Systeme, die
aus dem Raster
der klassischen Kontinuit"at fallen, d"urften
mit Hilfe des verallgemeinerten insensitiven Ergodensatzes 3.3 bearbeitbar sein.\index{verallgemeinerter insensitiver Ergodensatz}
Was vorzuf"uhren w"are.
\newline
Alleine der Ansatz, die Kontinuit"at von Funktionen\index{Kontinuit\"at von Funktionen}  
mittels Kommutatorrelationen zu beschreiben, ist insofern gleichsam die Spitze eines Kegels 
mit weitem "Offnungswinkel, als
innerhalb dieses Kegels 
verallgemeinernde Umfassungen der Theorie der Kontinuit"at zu erwarten sind.
Grob gesagt: So viele Typen von Mengensystem es gibt, so viele 
Kontinuit"atstypen gibt es. Diese Kegelspitze
er"offnet eine generalisierte, mehr als klassische Perspektive auf die Kontinuit"at:\index{Cantor, Georg}\index{Weierstrass, Karl}\index{Newton, Isaac}\index{Leibniz, Gottfried Wilhelm}
Es sei $f\in Y^{X}$ eine Funktion und $\mathcal{W}\subset 2^{X}$ ein Mengensystem, das die 
Definitionsmenge und die Wertemenge "uberdeckt, sodass
$$X\cup Y=\bigcup\mathcal{W}$$ ist. Genau dann, wenn 
die Funktion $f$
f"ur den gem"ass
(\ref{cloee}) und (\ref{clocee}) gegebenen H"ullenoperator $\mathbf{cl}_{\mathcal{W}}$
des Mengensystemes $\mathcal{W}$ 
die gem"ass (\ref{cloeee}) bestimmte Kommutatorrelation  
\begin{equation}\label{cant}
[f,\mathbf{cl}_{\mathcal{W}}]=\emptyset 
\end{equation}
erf"ullt, sagen wir, 
dass $f$ bez"uglich des Mengensystemes $\mathcal{W}$ kommutativ Cantor-stetig ist. Diese Benennung 
der kommutativen Cantor-Stetigkeit\index{kommutative Cantor-Stetigkeit einer Funktion} schlagen wir deshalb vor, 
weil die Idee der Kontinuit"at ein Archetyp\index{Archetyp der Kontinuit\"at} ist, den
Newton und Leibniz erstmals innerhalb des rationalen Diskurses mit Nachdruck
aussprechen und als solchen problematisieren; wobei 
sich Sir Isaac Newton im Gegensatz zu Gottfried Wilhelm Leibniz f"ur dessen Hinterfragung und kritische
Er"orterung weit weniger als f"ur Anwendungen 
der Kontinuit"atsvorstellungen im Calculus interessiert.
Vor allem Karl Weierstra"s
ist es, dem es gelingt,
die Analysis logisch zu b"andigen.
Hier nun formulieren wir 
in Form der kommutativen Cantor-Stetigkeit
auf der Ebene elementarer, handgreiflicher 
Mengenlehre\index{Mengenlehre} einen Kontinuit"atsbegriff, der 
die klassischen Kontinuit"atsvorstellungen objektivierend umfasst.
Die Elementarit"at des vorgestellten postklassischen Kontinuit"atsbegriffes
zeigt uns die Unmittelbarkeit, mit der
wir die M"oglichkeit, diesen postklassischen Kontinuit"atsbegriff zu formulieren, Georg Cantor
verdanken. Deswegen meinen wir, dass es angemessen ist, die kommutative Cantor-Stetigkeit
nach Georg Cantor zu benennen.\newline 
Der Kommutator
$[f,\mathbf{cl}_{\mathcal{W}}]$ ist 
gem"ass (\ref{cloeee})
eine mengenwertige, auf der
Potenzmenge $2^{X}$ definierte Abbildung und die 
Bedingung 
(\ref{cant}), welche die 
kommutative Cantor-Stetigkeit
definiert,
ist als die Bedingung zu lesen, dass die Abbildung $[f,\mathbf{cl}_{\mathcal{W}}]$ auf ihrer Definitonsmenge 
konstant ist und sie 
den Wert $\emptyset$  
hat. Damit ist (\ref{cant}) nicht nur eine Kommutatorrelation, sondern eine solche, die 
im Sinne  
einer Unreichhaltigkeit direkt eine Freiheit oder Reinheit formuliert.
Bez"uglich eines Mengensystemes 
$\mathcal{W}$, das die Reichhaltigkeit
\begin{equation}
\Bigl\{Q\in\mathcal{W}\in\setminus\{\emptyset\}:\exists (x,y)\in 2^{X}\times 2^{Y}: x\cap Q\not=\emptyset\not=Q\cap y\Bigr\}\not=\emptyset
\end{equation}
hat, gibt es demnach insofern keine kommutativ Cantor-stetige Funktion $f\in Y^{X}$, als in diesem Fall
die Funktion $[f,\mathbf{cl}_{\mathcal{W}}]$ nicht
auf der gesamten Potenzmenge $2^{X}$ definiert ist.
\section{Die Explizierung der kommutativen\\ Cantor-Stetigkeit}
Auch, wenn sich unser Interesse an abstrahierenden Reformen der Theorie der Kontinuit"at
an dieser Stelle zur"uckhalten soll, k"onnen wir uns kaum
damit zufrieden geben, es
alleine bei der blossen Einf"uhrung und Benennung der kommutativen
Cantor-Stetigkeit zu belassen.
Obwohl zun"achst noch im Dunkel, dr"angt doch die Frage,
was die kommutative Cantor-Stetigkeit 
darstellt und was durch sie unmittelbar impliziert ist.
Diese Frage l"asst Freiheiten offen und sie entspricht der gleichermassen  
mit Bearbeitungsfreiheitsgraden gestellten Aufgabe, Kriterien 
f"ur die kommutative Cantor-Stetigkeit
zu finden, die kommutative Cantor-Stetigkeit also zu paraphrasieren.
Notwendige 
Kriterien f"ur die kommutative Cantor-Stetigkeit sind das, was durch die letztere 
impliziert ist und daher sagen uns notwendige
Kriterien f"ur die kommutative Cantor-Stetigkeit, wie diese
angewandt werden kann; oder aber, jene Kriterien erlauben uns,
auszuschliessen, dass kommutative Cantor-Stetigkeit "uberhaupt gegeben ist.
Ob solche notwendigen Kriterien f"ur die kommutative Cantor-Stetigkeit
dabei durch sie \glqq unmittelbar\grqq impliziert sind, ist ein Geschmacksurteil.
Notwendige und hinreichende 
Kriterien f"ur die kommutative Cantor-Stetigkeit hingegen stellen sie als Begriff dar und sind Paraphrasen derselben.\newline 
Der Zweck der Phrasierung der kommutativen Cantor-Stetigkeit besteht 
darin, nach dem Vollzug solcher Phrasierungen Gegebenheiten als solche erkennen zu k"onnen,
welche die kommutative Cantor-Stetigkeit ausschliessen oder anzeigen.
Die Phrasierung der kommutativen Cantor-Stetigkeit dient also dem diagnostischen Zweck, 
kommutative Cantor-Stetigkeit 
bei einer jeweils gestellten konkreten Aufgabe zu erkennen, wenn sie
gegeben ist. Damit dient die Phrasierung der kommutativen Cantor-Stetigkeit dem,
den Begriff der
kommutativen Cantor-Stetigkeit der applikativen Arbeit mit  
demselben zu erschliessen: Um mit der 
kommutativen Cantor-Stetigkeit arbeiten zu k"onnen, m"ussen wir sie diagnostizieren k"onnen.
\newline 
Dass dabei die Transparenz einfacher logischer Darstellungen eines Begriffes
nicht in jedem Fall
die argumentative Arbeit, 
geschweige denn die applikative Arbeit
mit dem jeweiligen Begriff erleichtert,
ist eine bekannte Erfahrung. Insbesondere
die notwendigen und hinreichenden 
Kriterien f"ur die kommutative Cantor-Stetigkeit, d.h., deren logische Darstellungen, d.h., deren 
Paraphrasen, empfinden wir nichtsdestotrotz als um so befriedigender, je einfacher sie formuliert sind.
Und zwar, weil sie dann die Denk"okonomie verbessern, deren Optimierung 
ein klassisches Hintergrundprinzip der Mathematik "uberhaupt ist. Eine der grauen Eminenzen 
der Mathematik ist
insofern die 
Pointe.
\newline 
Uns treibt letztlich das Streben nach Reduktion: 
Um einfache Phrasierungen der kommutativen Cantor-Stetigkeit zu finden, orientieren
wir uns zun"achst an dem Fall, in dem sie bereits erkundet ist:
Wenn ein Mengensystem $\mathbf{T}$ eine Topologie ist, so k"onnen 
wir das Pr"adikat der
Faser 
$$\Bigl[\mathbf{cl}_{\mathbf{T}},\mbox{id}_{X^{X}}\Bigr]^{-1}(\{\emptyset\})$$ aller 
Abbildungen $f\in X^{X}$ der Menge 
$X=\bigcup\mathbf{T}$, die bez"uglich $\mathbf{T}$ stetig sind,
explizit angeben, wobei wir
aus Gr"unden der "Ubersichtlichkeit f"ur jedes Mengensystem $\mathcal{A}$ die Setzung
\begin{equation}\label{briest}
\mathring{\mathcal{A}}:=\mathcal{A}\setminus\{\emptyset\}
\end{equation} 
benutzen: 
Die Faser $[\mathbf{cl}_{\mathbf{T}},\mbox{id}_{X^{X}}]^{-1}(\{\emptyset\})$ umfasst die Menge
\begin{equation}
\Bigl\{f\in X^{X}:\forall\ {\rm U}\in \mathring{\mathbf{T}}\ \exists\ {\rm V}\in\mathring{\mathbf{T}}: f({\rm V})\subset{\rm U}\Bigr\}\ ,
\end{equation} 
die Menge der stetigen Abbildungen 
des topologischen Raumes $(X,\mathbf{T})$ in denselben.
Diese Menge stetiger Abbildungen bezeichen wir fortan im Zuge des Folgenden mit $\mathcal{C}_{+}(\mathbf{T})$: 
Denn wir bezeichen
allgemein f"ur 
jedes Mengensystem $\mathcal{A}$ und f"ur die durch es "uberdeckte Menge $Y:=\bigcup\mathcal{A}$ mit
\begin{equation}\label{expl}
\begin{array}{c}
\Bigl\{f\in Y^{Y}:
\forall\ {\rm A}\in \mathring{\mathcal{A}}\ 
\exists\ \bar{{\rm A}}\in\mathring{\mathcal{A}}:\ f(\bar{{\rm A}})\subset{\rm A}\Bigr\}\\
:=\mathcal{C}_{+}(\mathcal{A})\ ,\\
\Bigl\{f\in Y^{Y}:
\forall\ {\rm A}\in \mathring{\mathcal{A}}\ 
\exists\ \bar{{\rm A}}\in\mathring{\mathcal{A}}:\ \bar{{\rm A}}\subset f({\rm A})\Bigr\}\\
:=\mathcal{C}_{-}(\mathcal{A})\\
\end{array}
\end{equation}  
zwei einander korrespondierende Funktionenmengen, welche die Abstraktion des 
Stetigkeitsbegriffes  
der allgemeinen Topologie
"uber deren Allgemeinheitsgrad hinausgehend 
auf die direkte Weise
objektivieren. Die Indizierung der Funktionenmengen $\mathcal{C}_{+}(\mathcal{A})$ und $\mathcal{C}_{-}(\mathcal{A})$
ist dabei dadurch motiviert, dass
\begin{equation}\label{movexa}
\mathcal{C}_{-}(\mathcal{A})=\Bigl\{f\in Y^{Y}:
\forall\ {\rm A}\in \mathring{\mathcal{A}}\ 
\exists\ \bar{{\rm A}}\in\mathring{\mathcal{A}}:\ f^{-1}(\bar{{\rm A}})\subset {\rm A}\Bigr\}
\end{equation}
ist, wobei diese Funktionenmenge im allgemeinen von der Funktionenmenge
\begin{equation}\label{movexaa}
\begin{array}{c}
\Bigl\{f^{-1}\in Y^{Y}:
\forall\ {\rm A}\in \mathring{\mathcal{A}}\ 
\exists\ \bar{{\rm A}}\in\mathring{\mathcal{A}}:\ f^{-1}(\bar{{\rm A}})\subset {\rm A}\Bigr\}\\
=\mathcal{C}_{+}(\mathcal{A})\cap\{f^{-1}\in Y^{Y}\}
\end{array}
\end{equation}
differiert, wobei wir 
$$\{f^{-1}\in Y^{Y}\}=\{f\in Y^{Y}:\exists f^{-1}\in Y^{Y}\}$$
als die Menge aller Bijektionen der Menge $Y$ auf dieselbe auffassen. Das Urbild
$f^{-1}(X)$ existiert dabei f"ur jede Teilmenge $X\subset Y$ und f"ur jede Funktion $f\in Y^{Y}$ im Gegensatz
zu dem Element $\phi^{-1}(x)$, das 
f"ur alle $x\in \mathbf{P}_{2}\phi$ genau dann existiert, wenn $\phi$ eine invertierbare Funktion ist.
F"ur alle invertierbaren Abbildungen $\phi\in \{f^{-1}\in Y^{Y}\}$ gilt die "Aquivalenz
\begin{equation}
(\phi,\phi^{-1})\in \mathcal{C}_{+}(\mathcal{A})\times \mathcal{C}_{+}(\mathcal{A})\Leftrightarrow
(\phi,\phi^{-1})\in \mathcal{C}_{-}(\mathcal{A})\times \mathcal{C}_{-}(\mathcal{A})\ ,
\end{equation}
w"ahrend ausserdem f"ur jede invertierbare Abbildung $\phi\in \{f^{-1}\in Y^{Y}\}$ die "Aquivalenzen
\begin{equation}\label{movexd}
\begin{array}{c}
(\phi,\phi^{-1})\in \mathcal{C}_{+}(\mathcal{A})\times \mathcal{C}_{-}(\mathcal{A})\Leftrightarrow
\phi\in \mathcal{C}_{+}(\mathcal{A})\ ,\\
(\phi,\phi^{-1})\in \mathcal{C}_{-}(\mathcal{A})\times \mathcal{C}_{+}(\mathcal{A})\Leftrightarrow
\phi\in \mathcal{C}_{-}(\mathcal{A})
\end{array} 
\end{equation}
und 
\begin{equation}\label{movexe}
\phi\in \mathcal{C}_{\pm}(\mathcal{A})\Leftrightarrow\phi^{-1}\in \mathcal{C}_{\mp}(\mathcal{A})
\end{equation}
gelten. Die durch 
(\ref{movexa})-(\ref{movexe}) beschriebene
Sachlage motiviert uns dazu, f"ur jedwedes Mengensystem $\mathcal{A}$ exakt die
Funktionenmenge $\mathcal{C}_{+}(\mathcal{A})$ als die Menge der Cantor-stetigen 
Abbildungen bez"uglich des Mengensystemes $\mathcal{A}$ zu bezeichnen und demgegen"uber 
$\mathcal{C}_{-}(\mathcal{A})$ die Menge der konvers Cantor-stetigen 
Abbildungen bez"uglich des Mengensystemes $\mathcal{A}$
zu nennen; und dementsprechend von der Cantor-Stetigkeit bzw. von der 
konversen Cantor-Stetigkeit exakt der Elemente der Mengen $\mathcal{C}_{+}(\mathcal{A})$ bzw.
$\mathcal{C}_{-}(\mathcal{A})$ zu sprechen.\index{Cantor-Stetigkeit}\index{konverse Cantor-Stetigkeit}
\newline
Wir fragen uns, ob einfach
die Menge Cantor-stetiger Funktionen 
$\mathcal{C}_{+}(\mathcal{A})$ die Menge 
kommutativ Cantor-stetiger Funktionen ist;
oder, ob 
die Menge konvers Cantor-stetiger Funktionen $\mathcal{C}_{-}(\mathcal{A})$ oder auch eine andere einfache,
aus der Cantor-Stetigkeit oder der konversen Cantor-Stetigkeit aufgebaute Formulierung
mit der Faser
$$\Bigl\lbrack\mathbf{cl}_{\mathcal{A}},\mbox{id}_{Y^{Y}}\Bigr\rbrack^{-1}(\{\emptyset\})$$
identisch ist.
Um wenigstens diese Frage zu kl"aren, ordnen wir jedem 
Mengensystem $\mathcal{A}$ das Mengensystem
aller bez"uglich demselben abgeschlossener Mengen 
\begin{equation}
\{A\subset Y:\mathbf{cl}_{\mathcal{A}}(A)=A\}=\Bigl(\mbox{id}_{Y^{Y}}\Delta\mathbf{cl}_{\mathcal{A}}\Bigr)^{-1}
(\{\emptyset\})
\end{equation} 
zu, wobei sich die Indizierung der Identit"at $\mbox{id}_{Y^{Y}}$ aus dem
Zusammenhang ergibt, sodass
wir uns diese ersparen k"onnen. F"ur alle
bez"uglich $\mathcal{A}$ kommutativ
Cantor-stetigen Funktionen $f\in\lbrack\mathbf{cl}_{\mathcal{A}},\mbox{id}\rbrack^{-1}(\{\emptyset\})$ 
gilt per definitionem die Form von Abgeschlossenheit, welche die Implikation 
$$A\in(\mbox{id}\Delta\mathbf{cl}_{\mathcal{A}})^{-1}(\{\emptyset\})\Rightarrow
f(A)\in (\mbox{id}\Delta\mathbf{cl}_{\mathcal{A}})^{-1}(\{\emptyset\})$$
behauptet: Es gilt die "Aquivalenz
$$f\in\lbrack\mathbf{cl}_{\mathcal{A}},\mbox{id}\rbrack^{-1}(\{\emptyset\})\Leftrightarrow$$
$$\Bigl(A\in(\mbox{id}\Delta\mathbf{cl}_{\mathcal{A}})^{-1}(\{\emptyset\})\Rightarrow
f(A)\in (\mbox{id}\Delta\mathbf{cl}_{\mathcal{A}})^{-1}(\{\emptyset\})\Bigr)\ .$$
Das als die Faser bez"uglich $\mathcal{A}$ abgeschlossener Mengen 
$(\mbox{id}\Delta\mathbf{cl}_{\mathcal{A}})^{-1}(\{\emptyset\})$
formulierte Mengensystem bestimmt daher die 
Faser 
$\lbrack\mathbf{cl}_{\mathcal{A}},\mbox{id}\rbrack^{-1}(\{\emptyset\})$
all der Funktionen, die bez"uglich des Mengensystemes $\mathcal{A}$ kommutativ
Cantor-stetig sind. Jene Faser bez"uglich $\mathcal{A}$ abgeschlossener Mengen
$(\mbox{id}\Delta\mathbf{cl}_{\mathcal{A}})^{-1}(\{\emptyset\})$
induziert dar"uber hinaus
eine Faserung der Potenzmenge der von dem Mengensystem $\mathcal{A}$
"uberdeckten Menge $Y=\bigcup \mathcal{A}$,
n"amlich die Produktfaserung
$$\Bigl\{{\rm C}\Delta\mathbf{cl}_{\mathcal{A}})^{-1}(\{\emptyset\}):
{\rm C}\in  (\mbox{id}\Delta\mathbf{cl}_{\mathcal{A}})^{-1}(\{\emptyset\})\Bigr\}$$     
$$=\Bigl\{\{X\in 2^{Y}:\mathbf{cl}_{\mathcal{A}}(X)= {\rm C}\}:{\rm C}\in  (\mbox{id}\Delta\mathbf{cl}_{\mathcal{A}})^{-1}(\{\emptyset\})\Bigr\}$$
\begin{equation}\label{expla}
:=\langle Y\rangle_{\mathcal{A}}\ ,
\end{equation}
die eine Partition
$$\langle Y\rangle_{\mathcal{A}}\in\mathbf{part}(2^{Y})$$ 
der Potenzmenge $2^{Y}$ ist, was wir betonen. 
Und zwar ist $\langle Y\rangle_{\mathcal{A}}$ eine Partition
der Potenzmenge $2^{Y}$
in die jeweiligen Mengensysteme abschlussgleicher Mengen
$$\{X\in 2^{Y}:\mathbf{cl}_{\mathcal{A}}(X)= {\rm C}\}=({\rm C}\Delta\mathbf{cl}_{\mathcal{A}})^{-1}(\{\emptyset\})$$
f"ur die 
jeweiligen, 
bez"uglich $\mathcal{A}$ abgeschlossenen Mengen
${\rm C}\in(\mbox{id}\Delta\mathbf{cl}_{\mathcal{A}})^{-1}(\{\emptyset\})$.
Es erweisen sich schliesslich exakt die Funktionen $g\in Y^{Y}$ als kommutativ Cantor-stetig 
bez"uglich des Mengensystemes $\mathcal{A}$, welche die
Integrit"at der Partition
$\langle Y\rangle_{\mathcal{A}}$ in dem Sinn erhalten, dass die Identit"at
\begin{equation}\label{explc}
\begin{array}{c}
\Bigl\{\{\underbrace{g(\vartheta)}_{=\ \{g(x):x\in\theta\}}:\theta\in \Theta\}
:\Theta\in\langle Y\rangle_{\mathcal{A}}\Bigr\}\\
=\langle Y\rangle_{\mathcal{A}}
\end{array}
\end{equation}
gilt. Anhand dieser Partition $\langle Y\rangle_{\mathcal{A}}$ wird uns sichtbar werden,
wie die Menge aller kommutativ Cantor-stetigen Funktionen aus der Menge $Y^{Y}$, die Faser
$\lbrack\mathbf{cl}_{\mathcal{A}},\mbox{id}_{Y^{Y}}\rbrack^{-1}(\{\emptyset\})$,
zu explizieren ist.  
Mit Hilfe der folgenden beiden Bemerkungen 
werden wir schliesslich
die Aussage des Explizierungssatzes\index{Explizierungssatz} 3.8 an dem Aufbau der Produktfaserung $\langle Y\rangle_{\mathcal{A}}$
ablesen k"onnen, die z.B. die bereits gestellte Frage beantwortet,
ob die Menge aller kommutativ Cantor-stetigen Funktionen aus der Menge $Y^{Y}$ mit der gem"ass (\ref{expl}) formulierten 
Funktionenmenge
$\mathcal{C}_{+}(\mathcal{A})$
"ubereinstimmt. Eine wichtige Rolle in den 
folgenden Argumentationen spielt
dabei das jeweilige
Mengensystem $\mathbf{un}(\mathcal{A}^{c})$
der komplementfreien Teilmengen 
eines Mengensystemes $\mathcal{A}$.
Dieses Mengensystem der komplementfreien Teilmengen der von einem Mengensystem $\mathcal{A}$
"uberdeckten Menge $Y=\bigcup\mathcal{A}$ legen
wir f"ur alle Mengensysteme $\mathcal{A}$ als 
\begin{equation}\label{dexpld}
\mathbf{un}(\mathcal{A}^{c}):=\{X\in 2^{Y}:A\in\mathcal{A}^{c}\setminus\{\emptyset\}\Rightarrow A\not\subset X\}
\end{equation} 
fest.\index{komplementfreie Teilmenge eines Mengensystemes}
$\mathbf{un}(\mathcal{A}^{c})\ni\emptyset$ ist nie leer. Wir bemerken, dass
f"ur jedes Mengensystem $\mathcal{A}$, das $Y=\bigcup\mathcal{A}$ "uberdeckt,
die "Aquivalenz 
$$\mathbf{un}(\mathcal{A}^{c})=\{\emptyset\}\Leftrightarrow (x\in Y\Rightarrow
\{x\}\in \mathcal{A}^{c})$$
wahr ist, sodass in dem Fall $\mathbf{un}(\mathcal{A}^{c})=\{\emptyset\}$ die Gleichung
$$\mathbf{cl}_{\mathcal{A}}=\mbox{id}_{Y^{Y}}$$
gilt, sodass die kommutative Cantor-Stetigkeit f"ur alle Abbildungen der Menge $Y^{Y}$
gegeben ist: Es gilt die "Aquivalenz
\begin{equation}
\mathbf{un}(\mathcal{A}^{c})=\{\emptyset\}\Leftrightarrow 
\Bigl\lbrack\mathbf{cl}_{\mathcal{A}},\mbox{id}_{Y^{Y}}\Bigr\rbrack^{-1}(\{\emptyset\})=Y^{Y}\ ,
\end{equation}
die es nahelegt, den Fall der Minimalit"at des Mengensystem $\mathbf{un}(\mathcal{A}^{c})=\{\emptyset\}$
der komplementfreien Teilmengen 
eines jeweiligen Mengensystemes $\mathcal{A}$
als den Fall trivialer kommutativer Cantor-Stetigkeit
zu bezeichnen.\index{triviale kommutative Cantor-Stetigkeit} 
Die
f"ur alle Mengensysteme $\mathcal{A}$, die
$Y=\bigcup\mathcal{A}$ "uberdecken, g"ultige
"Aquivalenz 
\begin{equation}
\mathbf{un}(\mathcal{A}^{c})=\{\emptyset\}\Leftrightarrow \mathcal{A}\supset\{Y\setminus\{x\}:x\in Y\}
\end{equation}
macht dabei die triviale kommutative Cantor-Stetigkeit
im jeweils vorgelegten Einzelfall auf triviale Weise diagnostizierbar.
\newline
\newline
{\bf Bemerkung 3.6: Darstellungslemma}\newline
{\em Es sei $Y=\bigcup\mathcal{A}$ und $\mathcal{A}$ ein Mengensystem und $\mathbf{un}(\mathcal{A}^{c})$
das Mengensystem $\mathbf{un}(\mathcal{A}^{c})$
der komplementfreien Teilmengen des Mengensystemes $\mathcal{A}$
gem"ass} (\ref{dexpld}). 
{\em Ferner sei f"ur jede 
bez"uglich $\mathcal{A}$ abgeschlossene Menge ${\rm C}$  
\begin{equation}\label{dexplo}
{\rm C}_{\cap}:=\bigcap ({\rm C}\Delta\mathbf{cl}_{\mathcal{A}})^{-1}(\{\emptyset\})\ .
\end{equation}
Dann gilt}
\begin{equation}\label{exiplc}
\begin{array}{c}
\Bigr\{\{Q\cup {\rm C}_{\cap}:Q\in\mathbf{un}(\mathcal{A}^{c})\cap{\rm C}\}:{\rm C}\in(\mbox{id}\Delta\mathbf{cl}_{\mathcal{A}})^{-1}(\{\emptyset\})\Bigr\}\\
=\langle Y\rangle_{\mathcal{A}}\ .
\end{array}
\end{equation} 
\newline
\newline
{\bf Beweis:}\newline
Zu jeder Menge $Q\in\mathbf{un}(\mathcal{A}^{c})$, die komplementfrei ist, 
findet sich deshalb
eine bez"uglich $\mathcal{A}$ abgeschlossene Menge 
${\rm C}\in  (\mbox{id}\Delta\mathbf{cl}_{\mathcal{A}})^{-1}(\{\emptyset\})$ und 
ein Paar $X_{1},X_{2}\in\{X\in 2^{Y}:\mathbf{cl}_{\mathcal{A}}(X)= {\rm C}\}$, dessen symmetrische 
Differenz diese komplementfreie Menge 
$$Q=X_{1}\Delta X_{2}$$
identifiziert, weil die Vereinigung "uber das Mengensystem aller 
bez"uglich $\mathcal{A}$ abgeschlossener Mengen 
$$\bigcup(\mbox{id}\Delta\mathbf{cl}_{\mathcal{A}})^{-1}(\{\emptyset\})=Y$$
die Menge $Y$ "uberdeckt, sodass 
es zu jeder Menge $Q\in\mathbf{un}(\mathcal{A}^{c})$ eine 
bez"uglich $\mathcal{A}$ abgeschlossene Menge
${\rm C}(Q)\in  (\mbox{id}\Delta\mathbf{cl}_{\mathcal{A}})^{-1}(\{\emptyset\})$ 
gibt, in der $$Q\subset {\rm C}(Q) $$
liegt. Dabei ist das Mengensystem
$$\Bigl\{X\in 2^{Y}:X\subset ({\rm C}(Q)\setminus Q)\Bigr\}\subset \bigcup\langle Y\rangle_{\mathcal{A}}$$
nicht leer und f"ur jede Menge $X_{0}$ dieses Mengensystemes ist
$$(X_{0}\cup Q)\in  ({\rm C}(Q)\Delta\mathbf{cl}_{\mathcal{A}})^{-1}(\{\emptyset\})$$
ein expizites Beispiel daf"ur, dass die symmetrische Differenz
$$(X_{0}\cup Q)\Delta X_{0}=Q$$
ist.
Das gem"ass (\ref{dexpld}) festgelegte
Mengensystem
$\mathbf{un}(\mathcal{A}^{c})$
all der Teilmengen der Menge $Y$, die keine nicht leere Menge des 
komplement"aren Mengensystemes $\mathcal{A}^{c}$ enthalten, wird also von dem 
Mengensystem 
\begin{equation}\label{dyxpld}
\bigcup\Bigl\{\{X_{1}\Delta X_{2}:X_{1},X_{2}\in\Theta\}:\Theta\in \langle Y\rangle_{\mathcal{A}}\Bigr\}\supset \mathbf{un}(\mathcal{A}^{c})
\end{equation} 
inkludiert. Zu jeder 
bez"uglich $\mathcal{A}$ abgeschlossenen Menge ${\rm C}$  
existiert
der 
gem"ass (\ref{dexplo}) formulierte Schnitt ${\rm C}_{\cap}$
und das Mengensystem all dieser Schnitte 
$\{{\rm C}_{\cap}:{\rm C}\in (\mbox{id}\Delta\mathbf{cl}_{\mathcal{A}})^{-1}(\{\emptyset\})\}$.
Ausserdem gibt es zu jeder 
bez"uglich $\mathcal{A}$ abgeschlossenen Menge ${\rm C}$ das Spurmengensystem $\mathbf{un}(\mathcal{A}^{c})\cap {\rm C}$,
wobei das Mengensystem $(\mathcal{A}^{c})$ offensichtlich eine "Uberdeckung der Menge
$$\bigcup\mathbf{un}(\mathcal{A}^{c})=Y$$
ist.
Es gilt also
die Darstellung (\ref{exiplc}) der 
Produktfaserung $\langle Y\rangle_{\mathcal{A}}$.\newline
{\bf q.e.d.}
\newline
\newline
{\bf Bemerkung 3.7: Integrit"atskriterium}\newline
{\em Es sei $Y=\bigcup\mathcal{A}$ und $\mathcal{A}$ ein Mengensystem und $\mathbf{un}(\mathcal{A}^{c})$ das gem"ass} (\ref{dexpld})
{\em beschaffene Mengensystem.
Jede Funktion $g\in Y^{Y}$ verletzt die Integrit"at der Partition
$\langle Y\rangle_{\mathcal{A}}$ in dem Sinn, dass} (\ref{explc}) {\em gilt, genau dann nicht, wenn $g$ die
Komplementfreiheit in dem Sinn erh"alt, dass die
Implikation
\begin{equation}\label{expld}
Q\in\ \mathbf{un}(\mathcal{A}^{c})\Rightarrow g(Q)\in\ \mathbf{un}(\mathcal{A}^{c})
\end{equation}
wahr ist.}\newline\newline
{\bf Beweis:}\newline
F"ur alle Mengen
${\rm C}\in  (\mbox{id}\Delta\mathbf{cl}_{\mathcal{A}})^{-1}(\{\emptyset\})$, die bez"uglich $\mathcal{A}$ abgeschlossen sind,
gilt die Implikation 
$$X_{1},X_{2}\in\{X\in 2^{Y}:\mathbf{cl}_{\mathcal{A}}(X)= {\rm C}\}\ \land\ A\in\mathcal{A}^{c}\setminus \{\emptyset\}$$
$$\Rightarrow\ A\not\subset X_{1}\Delta X_{2}\ .$$
Daher gilt:
Genau dann, wenn eine abgeschlossene Menge ${\rm C}\in  (\mbox{id}\Delta\mathbf{cl}_{\mathcal{A}})^{-1}(\{\emptyset\})$
existiert, f"ur die 
es zwei Mengen
$$X_{1},X_{2}\in\{X\in 2^{Y}:\mathbf{cl}_{\mathcal{A}}(X)= {\rm C}\}$$
und eine Menge $A\in\mathcal{A}^{c}\setminus\{\emptyset\}$
gibt, f"ur die
$$A\subset g(X_{1}\Delta X_{2})$$
gilt,
bewahrt die Abbildung $g$ die Integrit"at der Partition $\langle Y\rangle_{\mathcal{A}}$ im Sinne der 
Gleichung (\ref{explc})
nicht. 
Jede Funktion $g\in Y^{Y}$ verletzt demnach die Integrit"at der Partition
$\langle Y\rangle_{\mathcal{A}}$ nicht, wenn f"ur sie die
Implikation (\ref{expld}) gilt. Soweit die eine Implikationsrichtung,
die das  Integrit"atskriterium behauptet.\newline
Andernfalls aber, wenn die
Implikation (\ref{expld}) nicht gilt, verletzt jede Funktion $g\in Y^{Y}$
die thematisierte Integrit"at:
Denn dann existiert eine komplementfreie Menge 
$Q_{\star}\in \mathbf{un}(\mathcal{A}^{c})$, f"ur die
$$g(Q_{\star})\not\in\mathbf{un}(\mathcal{A}^{c})$$ gilt.
Wenn es jene Menge komplementfreie $Q_{\star}$ gibt,
dann findet sich nach der Bemerkung 3.6
auch eine Menge ${\rm C}\in  (\mbox{id}\Delta\mathbf{cl}_{\mathcal{A}})^{-1}(\{\emptyset\})$ und 
ein Paar $$(X_{1},X_{2})\in({\rm C}\Delta\mathbf{cl}_{\mathcal{A}})^{-1}(\{\emptyset\})\times({\rm C}\Delta\mathbf{cl}_{\mathcal{A}})^{-1}(\{\emptyset\})\ ,$$
f"ur das
$$X_{1}\Delta X_{2}=Q_{\star}$$
ist.
\newline
{\bf q.e.d.}
\newline
\newline
Nun k"onnen wir die Aussage des folgenden Satzes 
ablesen:
\newline
\newline
{\bf Explizierungssatz 3.8: Explizierung der kommutativen Cantor-Stetigkeit}\index{Explizierungssatz}\newline
{\em F"ur jedes Mengensystem $\mathcal{A}$ gilt die Identit"at}
\begin{equation}\label{explb}
\begin{array}{c}
\lbrack\mathbf{cl}_{\mathcal{A}},\mbox{id}\rbrack^{-1}(\{\emptyset\})=\mathcal{C}_{+}(\mathcal{A}^{c})\cap
\mathcal{C}_{-}(\mathcal{A}^{c})\\=\mathcal{C}_{+}(\mathcal{A})\cap
\mathcal{C}_{-}(\mathcal{A})=\lbrack\mathbf{cl}_{\mathcal{A}^{c}},\mbox{id}\rbrack^{-1}(\{\emptyset\})\ ,
\end{array}
\end{equation}
{\em wobei 
$\mathcal{C}_{\pm}(\mathcal{A})$ und
$\mathcal{C}_{\pm}(\mathcal{A}^{c})$ die gem"ass} (\ref{expl}) {\em verfassten Funktionenmengen sind.} 
\newline
\newline
{\bf Beweis:}\newline
Wir k"onnen die kommutativ Cantor-stetigen Funktionen $g\in Y^{Y}$
als gerade diejenigen Funktionen
bestimmen,
f"ur welche die f"ur die
gem"ass (\ref{expla}) verfasste Partition
$\langle Y\rangle_{\mathcal{A}}$ die Invarianz im Sinne der 
Gleichung (\ref{explc}) gegeben ist: 
Jede Funktion $g\in Y^{Y}$, welche die Integrit"at der Partition bewahrt,
bildet die Mengen des Mengensystemes 
$$\Bigl\{\bigcup X:X\in\langle Y\rangle_{\mathcal{A}}\Bigr\}=(\mbox{id}\Delta\mathbf{cl}_{\mathcal{A}})^{-1}(\{\emptyset\})$$   
auf Mengen desselben ab und ist daher kommutativ Cantor-stetig bez"uglich $\mathcal{A}$. 
Wegen der 
Monotonie jeder Funktion $g\in Y^{Y}$, dass f"ur alle $Q,P\subset Y$
die 
Implikation 
$$Q\subset P\Rightarrow\ g(Q)\subset g(P)$$
wahr ist,
bildet dagegen jede Funktion $g\in Y^{Y}$, welche die Integrit"at der Partition
$\langle Y\rangle_{\mathcal{A}}$ verletzt, mindestens eine
bez"uglich $\mathcal{A}$ abgeschlossene Menge nicht auf eine solche ab.
Es sind also exakt die bez"uglich $\mathcal{A}$ kommutativ Cantor-stetigen Funktionen
der Funktionenmenge $Y^{Y}$, welche die Integrit"at der Partition
$\langle Y\rangle_{\mathcal{A}}$ erhalten. Daher sind die bez"uglich $\mathcal{A}$
kommutativ Cantor-stetigen Funktionen gerade die Elemente aus
$Y^{Y}$, welche die Komplementfreiheit erhalten, d.h., f"ur welche die Implikation
(\ref{explc}) wahr ist.\newline 
Allein diejenigen Funktionen $g\in Y^{Y}$, f"ur welche die
Implikation
\begin{equation}\label{explee}
{\rm A}\in \mathcal{A}^{c}\setminus\{\emptyset\}\Rightarrow (
\exists\ \bar{{\rm A}}\in\mathcal{A}^{c}\setminus\{\emptyset\}:\ g(\bar{{\rm A}})\subset{\rm A})
\end{equation}
gilt, sind also 
Kandidaten daf"ur,
bez"uglich $\mathcal{A}$
kommutativ Cantor-stetig zu sein. Denn die Implikation (\ref{explc}) impliziert diese Implikation (\ref{explee}).
Es folgt  
zun"achst lediglich die Inklusion 
\begin{equation}\label{exple}
\lbrack\mathbf{cl}_{\mathcal{A}},\mbox{id}\rbrack^{-1}(\{\emptyset\})\subset\mathcal{C}_{+}(\mathcal{A}^{c})\ .
\end{equation}
Exakt jede Funktion $g\in Y^{Y}$, welche die Komplementfreiheit erh"alt, erh"alt 
aber offensichtlich auch
die $\mathcal{A}$-Haltigkeit
in dem Sinn, dass f"ur $g$ die 
Implikation
\begin{equation}\label{expld}
Q\in\ (\mathbf{un}(\mathcal{A}^{c}))^{c}\Rightarrow g(Q)\in\  (\mathbf{un}(\mathcal{A}^{c}))^{c}
\end{equation}
gilt, wobei wir verdeutlichen, dass
$$(\mathbf{un}(\mathcal{A}^{c}))^{c}=2^{Y}\setminus(\mathbf{un}(\mathcal{A}^{c})\cup\{\emptyset\})$$
sei. Es zeigt sich uns, dass
die Erhaltung der Komplementfreiheit ein "Aquivalent der Erhaltung der $\mathcal{A}^{c}$-Haltigkeit
ist, die ein
"Aquivalent der Erhaltung der $\mathcal{A}$-Haltigkeit ist, die wiederum ein
"Aquivalent der Erhaltung der $\mathcal{A}$-Freiheit ist, die wir exakt den Mengen des 
Mengensystemes
\begin{equation}\label{dexpdd}
\mathbf{un}(\mathcal{A}):=\{X\in 2^{Y}:A\in\mathring{\mathcal{A}}\Rightarrow A\not\subset X\}
\end{equation}
zuschreiben. Dabei ist
$$(\mathbf{un}(\mathcal{A}^{c}))^{c}=\{X\in 2^{Y}:\exists\ A\in\mathring{\mathcal{A}}:A\subset X\}$$
\begin{equation}\label{dixpdd}
=:\mathbf{ov}(\mathcal{A})
\end{equation}
das Mengensystem der $\mathcal{A}$-haltigen Teilmengen
des Mengensystemes $\mathcal{A}$. 
Die Funktionen $g\in Y^{Y}$, welche die $\mathcal{A}$-Haltigkeit im Sinn der 
G"ultigkeit der Implikation (\ref{expld})
bewahren, sind exakt diejenigen, die sowohl Cantor-stetig als auch konvers Cantor-stetig sind im Sinn der 
Festlegungen (\ref{expl}). Wir d"urfen also die
Identifizierung der kommutativ Cantor-stetigen Funktionen aus der Menge $Y^{Y}$
\begin{equation}\label{exply}
\lbrack\mathbf{cl}_{\mathcal{A}},\mbox{id}\rbrack^{-1}(\{\emptyset\})=\mathcal{C}_{+}(\mathcal{A})
\cap\mathcal{C}_{-}(\mathcal{A})
\end{equation} 
formulieren. Dar"uber hinaus erkennen wir auch, dass die 
Erhaltung der $\mathcal{A}^{c}$-Haltigkeit die 
Identifizierung
\begin{equation}\label{eyply}
\lbrack\mathbf{cl}_{\mathcal{A}},\mbox{id}\rbrack^{-1}(\{\emptyset\})=\mathcal{C}_{+}(\mathcal{A}^{c})
\cap\mathcal{C}_{-}(\mathcal{A}^{c})
\end{equation}
erlaubt. Da die Identit"aten (\ref{eyply}) und (\ref{exply}) f"ur alle Mengensysteme $\mathcal{A}$ gelten 
ist $$\lbrack\mathbf{cl}_{\mathcal{A}},\mbox{id}\rbrack^{-1}(\{\emptyset\})= \lbrack\mathbf{cl}_{\mathcal{A}^{c}},\mbox{id}\rbrack^{-1}(\{\emptyset\})\ .$$
{\bf q.e.d.}
\newline
\newline
Der Satz 3.8 "uber die Explizierung der kommutativen Cantor-Stetigkeit
als gleichzeitige Cantor-Stetigkeit und konverse Cantor-Stetigkeit
zeigt uns, dass sich
im Fall der gleichzeitigen Gegebenheit 
der Cantor-Stetigkeit und der konversen Cantor-Stetigkeit
f"ur eine jeweilige Funktion die Differenz zwischen Nichtkomplementarit"at und 
Komplementarit"at aufhebt:
F"ur bez"uglich eines Mengensystemes $\mathcal{A}$
kommutativ Cantor-stetige Funktionen findet sich kein Unterschied zwischen den
komplement"aren Mengen aus $\mathcal{A}^{c}$ und den Lokalisierungen aus $\mathcal{A}$: 
F"ur jedes Mengensystemes $\mathcal{A}$ gibt es das jeweilige
Mengensystem dessen komplement"arer Hybride
\begin{equation}\label{eyplay}
\mathbf{hy}(\mathcal{A}):=\{\mathcal{Z}\subset\mathcal{A}\cup\mathcal{A}^{c}:\mathcal{Z}^{c}\subset\mathcal{A}\cup\mathcal{A}^{c}\}\ .
\end{equation}
Mit Hilfe des Mengensystemes komplement"arer Hybride formulieren wir, dass
f"ur jedes Mengensystem $\mathcal{A}$ die "Aquivalenz 
\begin{equation}\label{yplay}
\begin{array}{c}
\mathcal{B}\in\mathbf{hy}(\mathcal{A})\Leftrightarrow
\lbrack\mathbf{cl}_{\mathcal{A}},\mbox{id}\rbrack^{-1}(\{\emptyset\})=\lbrack\mathbf{cl}_{\mathcal{B}},\mbox{id}\rbrack^{-1}(\{\emptyset\})\ ,
\end{array}
\end{equation} gilt,
welche keineswegs behauptet, dass die 
in ihr angeschriebenen
H"ullenoperatoren $\mathbf{cl}_{\mathcal{A}}$ und $\mathbf{cl}_{\mathcal{B}}$ "ubereinstimmen.
Es besteht
offensichtlich die Eineindeutigkeit des 
universellen, auf der Klasse der Mengensysteme definierten  
Operators $\mathbf{cl}_{\mbox{id}}$ in Form der
f"ur
alle Mengensysteme $\mathcal{A}$ und $\mathcal{B}$  
g"ultigen "Aquivalenz
$$\mathbf{cl}_{\mathcal{A}}\not=\mathbf{cl}_{\mathcal{B}}\Leftrightarrow\mathcal{A}\not= \mathcal{B}\ .$$
Der Satz 3.8 "uber die Explizierung der kommutativen Cantor-Stetigkeit\index{Explizierungssatz}
erlaubt 
es uns, 
bez"uglich einer jeweiligen Zustandsraum"uberdeckung $\mathcal{A}\subset 2^{\mathbf{P}_{1}\xi^{0}}$
kommutativ Cantor-stetige Phasenfl"usse $\{\xi^{t}\}_{t\in\mathbb{R}}$
sehr einfach zu charakterisieren:
\newline
\newline   
{\bf Korollar 3.9:\newline Die kommutative Cantor-Stetigkeit des Phasenflusses}\newline
{\em Es sei $\{\xi^{t}\}_{t\in\mathbb{R}}$ ein Phasenfluss und $\mathcal{A}\subset 2^{\mathbf{P}_{1}\xi^{0}}$
eine "Uberdeckung dessen Zustandsraumes 
$$\mathbf{P}_{1}\xi^{0}=\bigcup\mathcal{A}\ .$$
Dann gilt die "Aquivalenzenkette\footnote{Es sei f"ur $n\in\mathbb{N}$ eine endliche Sequenz 
$A_{1},A_{2},\dots A_{n}$ von Ausagen gegeben. Wir notieren eine "Aquivalenzenkette gem"ass der Festlegung 
\begin{equation}\label{yyeply}
\begin{array}{c}
\Bigl((A_{1}\Leftrightarrow A_{2})\land (A_{2}\Leftrightarrow A_{3})\dots (A_{n-1}\Leftrightarrow A_{n})\Bigr)\\
\Leftrightarrow:\\
A_{1}\Leftrightarrow A_{2}\Leftrightarrow A_{3}\dots A_{n-1}\Leftrightarrow A_{n}\ ,
\end{array}
\end{equation}
falls $n>4$ ist und entsprechend, falls $n\leq 4$ ist.}
\begin{equation}\label{eyeply} 
\begin{array}{c}
(t\in\mathbb{R}\Rightarrow[\mathbf{cl}_{\mathcal{A}},\xi^{t}]=\emptyset)\\
\Leftrightarrow\\
\{\xi^{t}:t\in\mathbb{R}\}\subset\mathcal{C}_{\pm}(\mathcal{A})\ \Leftrightarrow\ 
\{\xi^{t}:t\in\mathbb{R}\}\subset\mathcal{C}_{\pm}(\mathcal{A})\\ 
\Leftrightarrow\\
\{\xi^{t}:t\in\mathbb{R}\}\subset\mathcal{C}_{\pm}(\mathcal{B})\ \Leftrightarrow\ 
\{\xi^{t}:t\in\mathbb{R}\}\subset\mathcal{C}_{\pm}(\mathcal{B})
\end{array}
\end{equation}
f"ur alle Mengensysteme $\mathcal{B}$ der 
gem"ass} (\ref{eyplay}) {\em verfassten
Menge $\mathbf{hy}(\mathcal{A})$.} \newline
\newline 
{\bf Beweis:}\newline
Die "Aquivalenzenkette (\ref{eyeply}) stellt im Kern die Behauptung dar, dass 
die kommutative Cantor-Stetigkeit des Phasenflusses
zu dessen Cantor-Stetigkeit "aquivalent ist: Wenn
die Cantor-Stetigkeit des Phasenflusses zu dessen 
kommutativer Cantor-Stetigkeit "aquivalent ist, gelten offenbar auch die
Extensionen dieser "Aquivalenz, welche die "Aquivalenzenkette (\ref{eyeply})
aufz"ahlt.\newline
Sei $Y=\bigcup\mathcal{A}$.
Nach dem Explizierungssatz 3.8 gilt die
Implikation 
$$t\in\mathbb{R}\Rightarrow[\mathbf{cl}_{\mathcal{A}},\xi^{t}]=\emptyset$$
genau dann, wenn die Inklusion
$$\{\xi^{t}:t\in\mathbb{R}\}\subset\mathcal{C}_{+}(\mathcal{A})\cap\mathcal{C}_{-}(\mathcal{A})$$
wahr ist.
Der Sachverhalt, dass f"ur jeden Phasenfluss $\{\xi^{t}\}_{t\in\mathbb{R}}$ gilt,
dass das Paar
$$(\{\xi^{t}:t\in\mathbb{R}\},\circ)$$
eine Gruppe ist, beinhaltet, dass 
f"ur jede reelle Zahl $t\in\mathbb{R}$ die
Abbildung $\xi^{t}$
invertierbar ist. 
Genau dann, wenn $f\in \{\xi^{t}:t\in\mathbb{R}\}\cap\mathcal{C}_{+}(\mathcal{A})$ ein 
bez"uglich $\mathcal{A}$ Cantor-stetiger Fluss ist,
ist, ist die Aussage
$$\forall\ {\rm A}\in \mathring{\mathcal{A}}\ 
\exists\ \bar{{\rm A}}\in\mathring{\mathcal{A}}:\ f(\bar{{\rm A}})\subset {\rm A}$$
wahr und genau dann gilt auch die Aussage
$$\forall\ {\rm A}\in \mathring{\mathcal{A}}\ 
\exists\ \bar{{\rm A}}\in\mathring{\mathcal{A}}:\ \bar{{\rm A}}\subset f^{-1}({\rm A})\ ,$$
exakt welche die konverse Cantor-Stetigkeit
der Inversion $f^{-1}$ ist, sodass
$f^{-1}\in \mathcal{C}_{-}(\mathcal{A})$ ist: Es gilt also
die Implikation
$$f\in\{\xi^{t}:t\in\mathbb{R}\}\cap\mathcal{C}_{+}(\mathcal{A})\Rightarrow  
f^{-1}\in\{\xi^{t}:t\in\mathbb{R}\}\cap \mathcal{C}_{-}(\mathcal{A})\ .$$
Auf analoge Weise folgt, dass
genau dann, wenn $f\in \{\xi^{t}:t\in\mathbb{R}\}\cap\mathcal{C}_{+}(\mathcal{A})$ ein
bez"uglich $\mathcal{A}$ konvers Cantor-stetiger Fluss ist,
dessen Inversion
$f^{-1}$ bez"uglich $\mathcal{A}$ Cantor-stetig ist und in $\mathcal{C}_{+}(\mathcal{A})$.
Es gilt also die "Aquivalenz
$$f\in\{\xi^{t}:t\in\mathbb{R}\}\cap\mathcal{C}_{+}(\mathcal{A})\Leftrightarrow  
f^{-1}\in\{\xi^{t}:t\in\mathbb{R}\}\cap \mathcal{C}_{-}(\mathcal{A})\ ,$$
wegen der auch deren kollektive Version 
$$\{\xi^{t}:t\in\mathbb{R}\}\subset\mathcal{C}_{+}(\mathcal{A})\Leftrightarrow 
\{\xi^{t}:t\in\mathbb{R}\}\subset\mathcal{C}_{-}(\mathcal{A})$$
gilt.
\newline
{\bf q.e.d.}
\newline
\newline
Wir wollen in diesem Zusammenhang gleich
auf
eine Aussage hinweisen, die der Aussage des Korollars 3.9 
"ahnelt:
\newline
\newline   
{\bf Bemerkung 3.10: $\mathcal{C}_{+}$-$\mathcal{C}_{-}$-Koinzidenz im endlichen Zustandsraum}\newline
{\em Es sei $\mathcal{A}$ ein Mengensystem, $Y=\bigcup\mathcal{A}$ und 
\begin{equation}
\mathbf{ab}(Y):=\{f\in Y^{Y}:f^{-1}\in Y^{Y}\}
\end{equation}
die Menge aller Autobolismen der Menge $Y$.
Die Implikation 
\begin{equation}\label{yeply} 
\begin{array}{c}
\mathbf{card}(Y)\in\mathbb{N}
\Rightarrow\\
\mathcal{C}_{+}(\mathcal{A})\cap\mathbf{ab}(Y)=\mathcal{C}_{-}(\mathcal{A})\cap\mathbf{ab}(Y)
\end{array}
\end{equation}
ist wahr.}
\newline
\newline 
{\bf Beweis:}\newline
Die Cantor-Stetigkeit einer Abbildung
$f\in \mathcal{C}_{+}(\mathcal{A})$ bez"uglich $\mathcal{A}$ liegt genau dann vor, wenn
$f$ jede nicht $\mathcal{A}$-haltige Menge $X\in\mathbf{un}(\mathcal{A})$ auf
eine ebensolche abbildet. Ist $X_{0}\in\mathbf{un}(\mathcal{A})$ eine nicht $\mathcal{A}$-haltige Menge, so ist 
die durch die Festlegung f"ur alle $\nu\in\mathbb{N}$
$$X_{\nu+1}:=f(X_{\nu})\in\mathbf{un}(\mathcal{A})$$
rekursiv definierte Sequenz
$$X_{0},X_{1},X_{2}\dots\in\mathbf{un}(\mathcal{A})$$
eine solche, dass all deren Glieder $X_{\nu}$ keine $\mathcal{A}$-haltigen Mengen sind.
Dabei ist
$$\mathbf{card}(\mathbf{un}(\mathcal{A}))\leq 2^{\mathbf{card}(Y)}\in\mathbb{N}\ ,$$ 
sodass es 
zwei nicht negative ganze Zahlen 
$\mu$ und $\nu<\mu$ gibt, f"ur die
$$X_{\nu}=X_{\mu}$$
ist. Ist
$f\in \mathbf{ab}(Y)$ speziell ein 
Autobolismus, so gibt es eine kleinste nat"urliche Zahl $n(X_{0})$, f"ur die f"ur alle $\nu\in\mathbb{N}$
$$X_{\nu}=X_{n(X_{0})+\nu}$$
ist. F"ur jede nicht $\mathcal{A}$-haltige Menge $X_{0}\in\mathbf{un}(\mathcal{A})$
gibt es also einen Zyklus
$$Z_{f}(X_{0}):=\{X_{0},X_{1},X_{2}\dots X_{n(X_{0})-1}\}\subset \mathbf{un}(\mathcal{A})\ ,$$
der ein Element des Mengensystemes  
$$\{Z_{f}(X_{0}):X_{0}\in\mathbf{un}(\mathcal{A})\}\in\mathbf{part}(\mathbf{un}(\mathcal{A}))$$
ist, welches das Mengensystem $\mathbf{un}(\mathcal{A})$ partioniert: Es gilt die
Implikation
$$X\in \Xi\ \land\ \Xi\in\{Z_{f}(X_{0}):X_{0}\in\mathbf{un}(\mathcal{A})\}\Rightarrow $$
$$f(X),f^{-1}(X)\in \Xi$$
f"ur alle Autobolismen $f\in \mathbf{ab}(Y)$.
Es ist daher auch 
die Inversion $f^{-1}$ jedes jeweiligen Autobolismus $f\in \mathbf{ab}(Y)$ von der Art, dass
$f^{-1}$ jede nicht $\mathcal{A}$-haltige Menge $X\in\mathbf{un}(\mathcal{A})$ auf
eine ebensolche abbildet. Wegen der "Aquivalenz (\ref{movexe}) gilt demnach
die Implikation 
(\ref{yeply}).\newline
{\bf q.e.d.}
\newline
\newline
Die 
Cantor-Stetigkeit ist
eine wichtige Begriffsbildung der abstrakten Ergodentheorie.
Es versteht sich von selbst, dass wir hier die 
Cantor-Stetigkeit, die konverse Cantor-Stetigkeit und die kommutative Cantor-Stetigkeit nur in den gr"obsten Z"ugen beschreiben,
deren 
Er"orterung wohl mindestens des Rahmens einer eigenen Monographie bed"urfte.
Unsere mehrteilige Darstellung der abstrakten Ergodentheorie\index{abstrakte Ergodentheorie} geht daher auch in ihrem zweiten Teil,
der die abstrakte Theorie der Sensitivit"at darstellt, in dem Sinn nicht weiter auf die Cantor-Stetigkeit ein,
dass da die Cantor-Stetigkeit nicht als eigenes Thema weiter erforscht wird. Wir verraten indess, dass
die Theorie der Sensitivit"at\index{Sensitivit\"at} die Theorie ist, deren Gegenst"ande im Bereich negierter 
Cantor-Stetigkeit liegen.
Hier lag uns vor allem daran, die kommutative Cantor-Stetigkeit zu explizieren, um sie der Anwendung
zu erschliessen, was sich als eine n"utzliche Vorbereitung auf den 
zweiten Teil erweisen wird.\newline
Wir wollen am Schluss dieses Abschnittes auf
weiterf"uhrende Abstraktionen ausblicken, die wir aber in keinem der beiden Teile
der Konzepte der abstrakten Ergodentheorie untersuchen werden:
Wir finden n"amlich schliesslich, dass wir die Betrachtung des Lokalisierungsverhaltens eines Phasenflusses
$\{\xi^{t}\}_{t\in\mathbb{R}}$
sogar noch allgemeiner als unter der Perspektive der kommutativen Cantor-Stetigkeit desselben
anschauen k"onnen: Es sei $\mathbf{P}_{1}\xi^{0}$
dessen Zustandsraum und
\begin{equation}\label{ripley} 
{\rm X}:2^{\mathbf{P}_{1}\xi^{0}}\to 2^{\mathbf{P}_{1}\xi^{0}}
\end{equation}
eine Abbildung der Zustandsmengen auf ebensolche und f"ur alle $t\in\mathbb{R}$ sei
der Kommutator 
$$[{\rm X},\xi^{t}]:2^{\mathbf{P}_{1}\xi^{0}}\to 2^{\mathbf{P}_{1}\xi^{0}}$$ dadurch
festgelegt, dass
f"ur alle $Z\subset\mathbf{P}_{1}\xi^{0}$ sein jeweiliger Wert
$$[{\rm L},\xi^{t}](Z):=\xi^{t}({\rm L}(Z))\Delta{\rm L}(\xi^{t}(Z))$$
sei. Ferner sei
\begin{equation}
[[\xi]]^{{\rm L}}:=\{{\rm L}(\xi(x,\mathbb{R})):x\in\mathbf{P}_{1}\xi^{0}\}\ ,
\end{equation}
sodass f"ur jede Zustandsraum"uberdeckung $\mathcal{A}$
\begin{equation}
[[\xi]]_{\mathcal{A}}=[[\xi]]^{\mathbf{cl}_{\mathcal{A}}}
\end{equation}
ist.
Die Implikation 
\begin{equation}
(\forall\ t\in\mathbb{R}\  [{\rm L},\xi^{t}]=\emptyset)\Rightarrow[[\xi]]^{{\rm L}}\in\mathbf{part}(\mathbf{P}_{1}\xi^{0}) 
\end{equation}
ist gleichermassen trivial wie allgemein:
Nichtsdestotrotz
abstrahiert sie den verallgemeinerten insensitiven Ergodensatz noch weiter; und zwar insofern sogar "uber die
abstrakte Ergodentheorie hinausgehend, als die mengenwertige 
Abbildung (\ref{ripley}) an sich zun"achst
in keinem Verh"altnis zu einer vorgegebenen Zustandsraum"uberdeckung $\mathcal{A}$
steht, die f"ur ein jeweiliges Ergodenproblem als relevant angesehen wird.
Es stellt sich allerdings auf der einen Seite die Frage, welche 
der mengenwertigen Operatoren ${\rm L}$
gem"ass (\ref{ripley}) so beschaffen sind, dass sie 
f"ur ergodentheoretische Fragen 
insofern auf vermittelte Weise 
von Interesse sind, als
es eine relevante Zustandsraum"uberdeckung $\mathcal{L}$ gibt, f"ur die
\begin{equation}
\mathbf{cl}_{\mathcal{L}}={\rm L}
\end{equation}
ist.
Andererseits stellt sich aber auch die Frage, welche Rolle andere Operatoren $\Gamma_{\mathcal{A}}$
als die H"ullenoperatoren $\mathbf{cl}_{\mathcal{A}}\not=\Gamma_{\mathcal{A}}$ bez"uglich einer jeweiligen Zustandsraum"uberdeckung $\mathcal{A}$
in der Ergodentheorie
spielen. $\Gamma_{\mbox{id}}\not=\mathbf{cl}_{\mbox{id}}$
sei hierbei
eine auf der Klasse der Mengensysteme definierte
universelle Konstruktion, sodass $\Gamma_{\mathcal{A}}$
deren operatorwertiger Wert f"ur eine
jeweilige Zustandsraum"uberdeckung $\mathcal{A}$ ist. Diese Frage stellt sich, zumal
andere universelle Konstruktionen $\Gamma_{\mbox{id}}\not=\mathbf{cl}_{\mbox{id}}$
als zu 
$\mathbf{cl}_{\mbox{id}}$ gleichwertig erscheinen:
\newline
Nicht nur, dass es zu jeder Menge 
Menge $P\subset Y$
die Mengensysteme 
$$\mathbf{H}^{-}_{\mathcal{A}}(P):=\{{\rm A}\in \mathcal{A}:{\rm A}\supset P\}$$
und
$$\mathbf{H}^{+}_{\mathcal{A}}(P):=\{{\rm A}\in \mathcal{A}:{\rm A}\subset P\}$$
gibt: F"ur jede Menge $Q$ und f"ur jedes Mengensystem $\mathcal{A}$ existieren 
die Mengensysteme 
\begin{equation}
\begin{array}{c}
\mathbf{H}^{-}_{\mathcal{A}}(Q):=\{{\rm A}\in \mathcal{A}:{\rm A}\supset Q\}\ ,\\
\mathbf{H}^{+}_{\mathcal{A}}(Q):=\{{\rm A}\in \mathcal{A}:{\rm A}\subset Q\}\ ,
\end{array}
\end{equation} 
sodass mit der Festlegung dieser beiden Mengensysteme $\mathbf{H}^{\pm}_{\mathcal{A}}(Q)$
f"ur jede Menge $Q$ und f"ur jedes Mengensystem $\mathcal{A}$
die universellen Operatoren
\begin{equation}
\mathbf{H}^{\pm}_{\mathcal{A}}\ \ \mbox{bzw.}\ \  \mathbf{H}^{\pm}_{\mbox{id}}(Q)\ \mbox{bzw.}\ \ 
\mathbf{H}^{\pm}_{\mathbf{P}_{2}}\circ\mathbf{P}_{1}
\end{equation} 
definiert sind: Und zwar ist $\mathbf{H}^{\pm}_{\mathcal{A}}$
auf der Klasse der Mengen erkl"art, die oft auch als 
Cantorsches Universum\index{Cantorsches Universum} oder einfach als 
Universum\index{Universum} bezeichnet wird. $\mathbf{H}^{\pm}_{\mbox{id}}(Q)$ ist hingegen auf der 
der Klasse der Mengensysteme definiert, w"ahrend $\mathbf{H}^{\pm}_{\mathbf{P}_{2}}\circ\mathbf{P}_{1}$  auf der 
der Klasse der Paare operiert, deren erste Komponente 
eine Menge und deren zweite Komponente
ein Mengensystem ist.
Wir erkennen nicht nur,
dass $$\bigcap\mathbf{H}_{\mathcal{A}^{\star}}^{-}=\mathbf{cl}_{\mathcal{A}}$$
der universalisierte H"ullenoperator des Mengensystemes $\mathcal{A}$ ist.
Offensichtlich ist auch, dass der Operator
\begin{equation}
\bigcup\mathbf{H}_{\mathcal{A}}^{+}=:\mathbf{int}_{\mathcal{A}}
\end{equation} 
insofern die Generalisierung desjenigen universellen Operators ist, der
dem Begriff des offenen Kernes einer Menge entspricht, welcher auch als 
deren Inneres bezeichnet wird:\footnote{Begriffe
entsprechen universellen Konstruktionen, die durch universelle Operatoren in Form 
von Abbildungen objektiviert sind.}
Wenn das Mengensystem $\mathbf{T}$ eine Topologie ist,
ist
die Restriktion des Operators
$\bigcup\mathbf{H}_{\mathbf{T}}^{+}=\mathbf{int}_{\mathbf{T}}$
auf die Potenzmenge der Menge $\bigcup\mathbf{T}$ gerade
diejenige Abbildung, 
die jeder Teilmenge
des topologischen Raumes $(\bigcup\mathbf{T},\mathbf{T})$
deren 
jeweiliges Innere
bez"uglich der Topologie $\mathbf{T}$
zuordnet. Daher liegt es nahe, die jeweiligen Werte
des Operators $\mathbf{int}_{\mathcal{A}}(Q)$
f"ur jedes Mengensystem $\mathcal{A}$ an jeder Stelle $Q$  
des Universums als das jeweilige Innere der Menge $Q$ bez"uglich des 
Mengensystemes $\mathcal{A}$ zu bezeichnen.
Setzen wir ferner f"ur alle $(j,k,l)\in\{0,1\}^{3}$ und f"ur 
jedes Mengensystem $\mathcal{A}$
\begin{equation}
\begin{array}{c}
\mathcal{A}^{(0)}:=\mathcal{A}\ ,\\
\mathcal{A}^{(1)}:=\mathcal{A}^{c}\ ,\\
({\rm O}_{0},{\rm O}_{1}):=(\bigcup,\bigcap)\ ,\\
(s(0),s(1)):=(+,-)\ ,
\end{array}
\end{equation} 
so k"onnen wir die acht Kombinationen
m"oglicher Konstruktionen von Operatoren
\begin{equation}
\mathbf{H}_{jkl}^{\mathcal{A}}:={\rm O}_{j}\mathbf{H}^{s(k)}_{\mathcal{A}^{(l)}} 
\end{equation} 
bin"ar indizieren, wobei sichtbar wird, dass
\begin{equation}
\begin{array}{lll}
\mathbf{H}_{000}^{\mathcal{A}} &= &\mathbf{int}_{\mathcal{A}}\ ,\\
\mathbf{H}_{111}^{\mathcal{A}} &= & \mathbf{cl}_{\mathcal{A}}
\end{array}
\end{equation} 
ist. Es zeigt sich, dass der Begriff des Abschlusses bez"uglich eines 
Mengensystemes dem Begriff des Inneren bez"uglich eines 
Mengensystemes in dreierlei, also in allen m"oglichen Hinsichten
gegen"uberliegt.  

\section{Die Emanzipation der Ergodentheorie von der Topologie}
Wir haben nun eine Allgemeinheitsstufe jenseits 
der der allgemeinen Topologie erreicht.
Wir k"onnen Ergodentheorie auf allgemeiner mengentheoretischer Basis betreiben.
Denken wir nun zur"uck an 
den Anfang des zweiten Kapitels, an
die Einf"uhrung 
freier Attraktoren, an die bez"uglich eines Mengensystemes $\mathcal{A}$
gem"ass (\ref{ienatt}) und (\ref{idkoh}) verfassten Mengensysteme 
$\mbox{{\bf @}}(\Phi,\mathcal{A})$,\index{freier Attraktor bez\"uglich eines Mengensystemes}
deren kennzeichende Variable $\mathcal{A}$ wir dabei das jeweilige, die freien
Attraktoren $a\in \mbox{{\bf @}}(\Phi,\mathcal{A})$ 
relativierende Mengensystem\index{relativierendes Mengensystem eines freien Attraktors} nennen:
Uns vergegenw"artigt sich wieder,
dass in der Allgemeinheit der Wahl des 
Mengensystemes $\mathcal{A}$ diejenige 
Verallgemeinerung des herk"ommlichen Koh"arenzkriteriums besteht, die
sogar
noch "uber die Allgemeinheitsstufe der allgemeinen Topologie hinausgeht.
Und die Frage, die im Hintergrund blieb und die uns vielleicht denndoch umtrieb,
die Frage,
was es sei, was die spezielle Wahl ausgerechnet einer 
Topologie
als das relativierende Mengensystem $\mathcal{A}$, 
auszeichne,
diese Frage steht nun vor 
uns da.
\newline
Der verallgemeinerte insensitive Ergodensatz 3.3\index{verallgemeinerter insensitiver Ergodensatz} gibt uns aber jetzt eine 
lakonische Antwort darauf, inwieweit  
Topologien als 
jene relativierenden Mengensysteme
bei der Verfassung jeweiliger
freier Attraktoren im Hinblick auf die bis hierher untersuchte Ergodentheorie 
ausgezeichnet
sind.
Der verallgemeinerte insensitive Ergodensatz sagt uns:
Im Hinblick auf ihn sind Topologien als solche gar nicht
ausgezeichnet.
Es ist lediglich 
die 
Relevanz 
stetiger Flussfunktionen f"ur andere als f"ur ergodentheoretische Belange, welche die 
Bevorzugung von Topologien
bedingt, die 
von der Ergodentheorie bis zum verallgemeinerten insensitiven Ergodensatz losgel"ost ist.
\newline
Sogleich erkennen wir diesen Sachverhalt, wenn wir uns ausmalen,
dass sich auch die Theorie der Sensitivit"at\index{Sensitivit\"at} von Phasenfl"ussen 
analog zu der bis hierher untersuchten Ergodentheorie auf allgemeiner 
mengentheoretischer Basis 
statt auf topologischer Basis  
konstituieren lassen k"onnte: Es k"onnte sich zeigen,
dass es auch durch die 
Theorie der Sensitivit"at
keine Bevorzugung von Topologien gibt.
Die Theorie der Sensitivit"at stellt dabei den Kern der 
sogenannten
Chaostheorie\index{Chaostheorie}
dar. Wir rechnen die Theorie der Sensitivit"at zur Ergodentheorie.\index{Ergodentheorie}
Im zweiten Teil der abstrakten Ergodentheorie wollen wir diese Theorie der Sensitivit"at
darstellen. Wir werden da sehen, dass die Cantor-Stetigkeit des Phasenflusses 
insofern als die vollst"andige Negation der Sensitivit"at zu verstehen ist, als es verschiedene 
Formen der Sensitivit"at ist, deren schw"achste exakt die 
Negation der Cantor-Stetigkeit ist. Der entsprechend abstrakt
auf der Grundlage der Mengenlehre verfasste Sensitivit"atsbegriff, der "uber die allgemeine
Topologie hinausgehend 
generalsisiert formuliert werden wird,
kennt verschiedene Grade und Formen der Sensitivit"at. Es gibt sogar noch st"arkere Formen als die maximale 
Sensitivit"at.
Was wir hierbei unter maximaler Sensitivit"at stetiger
reeller Phasenfl"usse verstehen, haben wir im Traktat "uber den elementaren Quasiergodensatz \cite{raab}
bereits gesagt, wo wir maximale Sensitivit"at f"ur glatte Phasenfl"usse vermuteten.
Die Vermutung 
maximaler Sensitivit"at
der Abhandlung des elementaren Quasiergodensatzes wird hierbei auch im 
zweiten Teil der abstrakten Ergodentheorie noch nicht gezeigt werden.
Dass es die Negation der Cantor-Stetigkeit des Phasenflusses ist, die 
f"ur jede Form der
Sensitivit"at notwendig ist,
k"onnte glauben machen, dass Sensitivit"at einfach Diskontinuit"at als Negation jeweiliger Formen der Kontinuit"at der 
Phasenfl"usse ist. F"ur dynamische Systeme hoher Kontinuit"at, wie sie beispielsweise
in \cite{wirr} dargestellt werden, 
gibt es indess sehr wohl Sensitivit"at und chaotisches Verhalten.
\newline
Damit sich 
auch durch diesen Teil der Ergodentheorie, der die 
Sensitivit"at betrifft,
keine Auszeichnung der Topologien
erg"abe,
m"usste allerdings der 
Sensitivit"atsbegriff nicht nur auf der
Allgemeinheitsstufe der allgemeinen 
Mengenlehre formulierbar sein.
Wir m"ussten "uberdies noch feststellen, dass
wir mit diesem mengentheoretischen
Sensitivit"atsbegriff abstrakte 
"Aquivalente zu denjenigen 
Aussagen finden, die
wir von der bislang bekannten   
Theorie der Sensitivit"at her kennen
und die wir als f"ur diese 
Theorie charakteristische Aussagen werten.
\newline
In diesem Fall
g"abe es 
sowohl durch die Ergoden- als auch durch die Chaostheorie
keine Bevorzugung von Topologien und dementsprechend
w"are es dann auch angemessen, den
Begriff des Attraktors in der entsprechenden
mengentheoretischen Allgemeinheit zu verfassen;
und dabei einzusehen, dass der 
Begriff des Attraktors lediglich
als ein wesentlich topologischer Begriff erschien.
Dass es dabei m"oglich ist, den Begriff des
Attraktors in
mengentheoretischer Allgemeinheit
zu formulieren,
wissen wir bereits, weil wir dies
in Gestalt freier 
Attraktoren vorf"uhrten. 
\newline
Wir k"onnen aber
den herk"ommlichen Attraktorbegriff durchaus abstrahieren und 
dabei dennoch
insofern innerhalb der Topologie
verbleiben, als wir \glqq Attraktoren als rein topologische
Epikonstrukte\grqq
verfassen k"onnen; wie wir dies,
der "Uberschrift des 
zweiten Kapitels gem"ass, ebenda
zeigten. Es gibt dabei, wie wir nun wissen, allerdings keinen Grund daf"ur,
sich auf Topologien einzuschr"anken,
der sich
aus der allgemeinen mengentheoretischen Ergodentheorie ergibt -- soweit sie uns bis jetzt bekannt wurde.
\newline 
Die Spannweite des Begriffes des topologischen Attraktors ist gross:
Die Klasse der Werte des Operators $\underline{\mbox{{\bf @}}}$ umfasst
beispielsweise die Klasse der Topologien.   
Es gibt also Mengen topologischer Attraktoren,
die in keinem fachspezifischen Verh"altnis mehr zu dynamischen
Systemen stehen. Der einfachste Fall einer Menge topologischer Attraktoren beispielsweise ist ein
Mengensystem inkoh"arenter Attraktoren 
einer einelementigen Topologienmenge $\{\mathbf{T}\}$
des Mengensystemes 
$$\underline{\mbox{{\bf @}}}(\{\mathbf{T}\})$$
\begin{equation}
\begin{array}{c}
=\{{\rm U}\in\mathbf{T}:\forall\ a,b\subset{\rm U}\ \exists\ \vartheta\in\mathbf{T}:\\
a\cap\vartheta\not=\emptyset\not=\vartheta\cap b \}=\mathbf{T}\ ,
\end{array}
\end{equation}
sodass also 
$$\underline{\mbox{{\bf @}}}(\{\mbox{id}\})=\{\mbox{id}\}$$
ist. F"ur jede Menge topologischer Attraktoren einer einelementigen Topologienmenge $\{\mathbf{T}\}$
gilt daher erst recht f"ur jede Koh"arenz ${\rm X}^{\star}$, 
die auf der Menge (\ref{aegy}) gegeben
und die
gem"ass der Definition 2.1 
zul"assig ist, 
f"ur die also $\bigcup{\rm X}^{\star}=\bigcup\mathbf{T}$ ist, dass
$\underline{\mbox{{\bf @}}}(\{\mathbf{T}\},{\rm X}^{\star})=\mathbf{T}$ ist;
sodass f"ur jede zul"assige, auf der Menge (\ref{aegy}) gegebene Koh"arenz ${\rm X}^{\star}$ und
f"ur jede zul"assige
Kadenz ${\rm X}_{\star}$ die 
Menge koh"arenter und nicht freier Attraktoren
\begin{equation}
\underline{\mbox{{\bf @}}}(\{\mathbf{T}\},{\rm X}^{\star},{\rm X}_{\star})=
{\rm X}_{\star}
\end{equation}
die Kadenz ${\rm X}_{\star}$ selbst ist: Das heisst, dass jedes
Mengensystem ${\rm X}_{\star}$ ein topologischer Attraktor beliebiger 
Koh"arenz ${\rm X}^{\star}$ einer einelementigen Topologienmenge ist. Jede
Menge ist ein topologischer Attraktor einer einelementigen Topologienmenge:
Wir k"onnen das Cantorsche Universum\index{Cantorsches Universum} als Universum topologischer Attraktoren 
auffassen.\newline
Es g"alte nun vielleicht, 
die durch die Definition 2.1 
erzeugten Szenarien
in den F"allen zu studieren, die weder der triviale Fall 
einer einelementigen Topologienmengen sind, noch die F"alle, in denen
der Begriff des topologischen Attraktors mit dem herk"ommlichen Begriff des Attraktors
in eins f"allt: Es ist zu erwarten, dass zwischen den Topologien als den Werten
des Operators $\underline{\mbox{{\bf @}}}$ f"ur  
einelementige Topologienmengen und den 
Attraktoren im herk"ommlichen Sinn Welten liegen. Nun hier die Aufgabe anzugehen,
diese Welten auszuforschen, w"urde per se unserem Thema untreu. Denn, 
folgten wir dieser Aufgabe, so
studierten wir zwar durchaus
Attraktoren, n"amlich topologische Attraktoren und wir
blieben dabei insofern auf der Bahn einer Generalisierung, die von dem Ursprung unserer Untersuchungen
ausstrahlt, der 
unsere thematische Mitte ist, die der Attraktor im herk"ommlichen Sinn objektiviert.
Wir blieben wohl auf jener Bahn, indess, 
reden wir mit Metaphern der Astronomie, es w"are folgendermassen: 
Da diese Generalisierungsbahn gleichsam 
hyperbolisch ist, f"uhrt sie doch, wie wir soeben sahen, durch das Cantorsche Universum
und immer weiter, beliebig weit von unsere thematischen Mitte fort.
Wir verl"oren unser Thema aus den Augen.
\newline 
Wenn wir die der Himmelsmechanik entlehnte Metaphorik fortsetzen: Wir beobachten, dass 
die Generalisierungsbahn
freier Attraktoren 
anders als die quasi hyperbolische
Generalisierungsbahn der gem"ass der Defintion 2.1 konstruierten topologischen Attraktoren\index{topologischer Attraktor} 
verl"auft:
Die Generalisierungsbahn
freier Attraktoren $a\in \mbox{{\bf @}}(\Phi,\mathcal{A})$
f"ur die bez"uglich eines relativierenden Mengensystemes $\mathcal{A}$
gem"ass (\ref{ienatt}) und (\ref{idkoh}) verfassten Mengensysteme 
$\mbox{{\bf @}}(\Phi,\mathcal{A})$\index{freier Attraktor bez\"uglich eines Mengensystemes}
f"uhrt uns n"amlich nicht von unserem thematischen Zentrum fort.
Diese Bahn, deren Parametersierung quasi die
das jeweilige, die freien
Attraktoren $a\in \mbox{{\bf @}}(\Phi,\mathcal{A})$ relativierende
Mengensystem $\mathcal{A}$ ist, f"uhrt uns nicht
von der ergodischen Relevanz der Mengensysteme 
$\mbox{{\bf @}}(\Phi,\mathcal{A})$ weg, "uber die der verallgemeinerte insensitive Ergodensatz 3.3 seine Aussage macht;
obgleich wir hier ebenfalls in die Allgemeinheit extendieren, wenn wir 
das relativierende Mengensystem $\mathcal{A}$ entsprechend allgemein einstellen.
\newline
Dass wir dabei an das relativierende Mengensystem $\mathcal{A}$
nicht den Anspruch zu haben brauchen, dass 
es eine Topologie ist, dies
demonstriert auf eklatante Weise, dass der Begriff des Attraktors nicht wesentlich
topologisch verfasst ist!
Gewiss, wir trachten danach, von dem Komfort Gebrauch zu machen, den uns 
die Befunde der Topologen
zur Verf"ugung stellen, wenn wir das 
relativierende Mengensystem $\mathcal{A}$ 
bei der Untersuchung jeweiliger freier
Attraktoren $a\in \mbox{{\bf @}}(\Phi,\mathcal{A})$
als eine Topologie
w"ahlen. Denn seit es zuvorderst wohl
S. Smale war, der in den sp"aten sechziger
Jahren des zwanzigsten Jahrhunderts
z.B. in seiner klassischen Arbeit \cite{male}
die Verbindung zwischen der Theorie
dynamischer Systeme und der mengentheoretischen Topologie 
etablierte, hat sich
die Anwendung topologischer 
Methoden auf dynamische Systeme bew"ahrt.
Wir wollen daher, wenn uns bei der Untersuchung jeweiliger freier
Attraktoren $a\in \mbox{{\bf @}}(\Phi,\mathcal{A})$
eine Topologie als das relativierende Mengensystem $\mathcal{A}$
zur Verf"ugung steht,
jenen Komfort maximieren, den es darstellt, auf Befunde der 
Topologie zur"uckgreifen zu k"onnen.
Sodass wir
bestrebt sind,
sofern bei einer jeweils gestellten Aufgabe
diese Freiheit besteht, 
m"oglichst symmetrische Topologien zugrunde zu legen. Wir versuchen einer jeweiligen Aufgabe
etwa auf der Hausdorffschen Skala der ${\rm T}_{0}$-Topologien bis hin zu den ${\rm T}_{4}$-Topologien
m"oglichst eine ${\rm T}_{4}$-Topologie zu Grunde zu legen.
Die Option 
der Wahl des relativierenden Mengensystemes 
besteht nichtsdestotrotz unabh"angig vom Blick auf die Einrichtungen der Topologie:
Der erste Teil der abstrakten Ergodentheorie hat sich insoweit
von der Topologie emanzipiert.
\newline 
Diese Emanzipation heisst aber keineswegs, dass die Anwendung topologischer 
Methoden in dem Teil der Theorie der dynamischen Systeme, den die
Ergodentheorie darstellt, unfruchtbar w"are.
Deshalb haben wir es an der 
Objektivierung des Begriffes des dynamischen Systemes
als Phasenflussgruppe vermisst,
dass
keine Topologisierung des jeweiligen
Zustandsraumes, 
der Definitionsmenge der Gruppenelemente
der jeweiligen Phasenflussgruppe, festgelegt wird. Dem wurde in dem Traktat "uber den elementaren
Quasiergodensatz \cite{raab} Abhilfe geschafft:
Das Paar
$((\Xi,\circ),\mathbf{T}(\zeta))$
nennen wir genau dann
ein topologisches dynamisches System\index{topologisches dynamisches System} oder ein Smale-System,\index{Smale-System}
wenn $(\Xi,\circ)$ eine Phasenflussgruppe von Mitgliedern des Phasenflusses ist,
die Bijektionen des Zustandsraumes $\zeta$ auf sich selbst sind, welchen
die Topologie $\mathbf{T}(\zeta)$ topologisiert.
Der Begriff des kontinuierlichen dynamischen Systemes differiert demnach von dem des topologischen dynamischen Systemes.
Die Benennung 
topologischer dynamischer Systeme 
nach S. Smale lag nahe, weil vor allem er die
Topologie und die Theorie dynamischer Systeme
zusammenf"uhrte.\newline 
Wir wollen aber nun den zu der Einf"uhrung der Smale-Systeme 
analogen Schritt gehen und exakt die Paare 
\begin{equation}\label{kurzzf} 
\Bigl((\Xi,\circ),\mathcal{A}\Bigr)
\end{equation} 
als Cantor-Systeme\index{Cantor-System} bezeichnen, deren erste Komponente $(\Xi,\circ)$ eine Phasenflussgruppe ist
und deren zweite Komponente $\mathcal{A}\subset 2^{Z}$ den
gemeinsamen Zustandsraum 
\begin{equation}
Z=\bigcup \mathcal{A}=\bigcup_{\xi\in \Xi}\mathbf{P}_{1}\xi=\bigcap_{\xi\in \Xi}\mathbf{P}_{1}\xi 
\end{equation}
der Phasenflussgruppe $\Xi$ "uberdeckt. Die Benennung nach Georg Cantor\index{Cantor, Georg}
kommentieren wir genauso wie die Benennung der Cantor-Stetigkeit.
Die Emanzipation der Ergodentheorie von der Topologie ist f"ur den ersten Teil der 
abstrakten Ergodentheorie vollzogen. \newline
Diese Emanzipation stellt aber keinen Bruch mit der 
Topologie dar!
Es k"onnte sein, dass es selbst innerhalb der 
Ergodentheorie im weiteren Sinn, gem"ass welchem die Theorie der Sensitivit"at
zur Ergodentheorie z"ahlt, Topologien nicht ausgezeichnet sind.
Um so mehr tritt hervor, dass offenbar die Diskussion nat"urlicher Topologien
bevorzugt wird.
Dass es sich bei der anhand 
klassischer 
Lehrb"ucher wie z.B. \cite{anos}, \cite{oran}, \cite{walt} oder \cite{deva}
beobachtbaren
Bervorzugung der Diskussion nat"urlicher Topologien letztendlich um einen 
exemplarischen
Anthropomorphismus\index{Anthropomorphismus}
handele, ist keine unanfechtbare
Hintergrundeinsicht, sondern eine durchaus anfechtbare Meinung.
Inwiefern diese Meinung angefochten werden k"onnte?
Dass f"ur jede nat"urliche Topologie $\mathbf{T}(n)$ eines reellen
Raumes $\mathbb{R}^{n}$
endlicher Dimension $n\in\mathbb{N}$ die beiden 
Kommutationen
\begin{equation}
\begin{array}{c}
\Bigl\lbrack\mathbf{cl}_{\mathbf{T}(n)}\mathbf{P}_{1}\oplus\mathbf{cl}_{\mathbf{T}(n)}\mathbf{P}_{2},\ (\mathbf{P}_{1}\cup\mathbf{P}_{2})\oplus (\mathbf{P}_{1}\cup\mathbf{P}_{2})\Bigr\rbrack=\emptyset\ ,\\
\quad\\
\Bigl\lbrack\mathbf{cl}_{\mathbf{T}(n)}\mathbf{P}_{1}\oplus\mathbf{cl}_{\mathbf{T}(n)}\mathbf{P}_{2},\ (\mathbf{P}_{1}\cap\mathbf{P}_{2})\oplus (\mathbf{P}_{1}\cap\mathbf{P}_{2})\Bigr\rbrack=\emptyset
\end{array}
\end{equation}
gelten, ist 
eine vom Menschen und seiner perspektivischen Voreingenommenheit g"anzlich
unabh"angige Gegebenheit der Logik, die nat"urliche Topologien auszeichnet.
Die Logik ist hierbei bislang nicht als Anthropomorphismus abgetan.\newline
Es ist vielleicht gerade die Logik, die sogar 
als der Schauplatz des innigsten Kontaktes des Menschen zu dem angesehen werden kann,
was jenseits des Anthropomorphismus liegt.
Es ist denkbar, dass sich innerhalb einer Theorie der 
Robotik\index{Robotik} zeigt,
weshalb und
unter welchen Umst"anden -- in welcher Welt, in welcher {\em Physik} -- 
es von Vorteil ist, sich auf der Grundlage nat"urlicher Topologien zu orientieren.
Wie weit wir von einer solchen Sichtweise entfernt sind, die
den Anthropomorphismus implementiert? Wir sind hier mir "uberschaubareren Fragen befasst.
\newline
Die Literatur zu Attraktoren {\em kontinuierlicher} dynamischer Systeme er"ortert vorzugsweise die Attraktoren endlichdimensionaler reeller oder komplexer R"aume 
bez"uglich deren jeweiliger nat"urlicher Topologie; das Ausmass dieser Bevorzugung nat"urlicher Topologien liegt 
vielleicht nicht alleine in deren Elementarit"at des leichten Einstieges begr"undet; vermutlich kommt hinzu,
dass die 
Relevanz
der nat"urlichen Topologie f"ur die empirischen Wissenschaften hoch ist, was an deren dominanter Appliziertheit
in den empirischen Wissenschaften liegt.
Die Thermodynamik, 
das Herkunftsgebiet der Ergodentheorie, k"onnte sehr von den Einsichten durchaus auch der abstrakteren 
Mathematik dynamischer Systeme profitieren. Indess wird beispielsweise in dem Klassiker der 
Lehrbuchliteratur der Thermodynamik \cite{reich} weder der Begriff der Topologie 
noch der Begriff der $\sigma$-Algebra zu finden sein.\footnote{Bei den dominant diskutierten Attraktoren endlichdimensionaler reeller oder komplexer R"aume, den
Attraktoren reeller oder komplexer Flussfunktionen also,
braucht deren Bezug auf
die jeweilige nat"urliche Topologie meistens nicht ausdr"ucklich erw"ahnt zu werden; es ist da in etwa so, wie 
bei der einfachen Rede von der Stetigkeit:
Auch
der Bezug der Stetigkeit auf die jeweiligen nat"urlichen Topologien
wird nicht ausdr"ucklich erw"ahnt, wenn von der Stetigkeit 
einer
zwischen endlichdimensionalen reellen oder komplexen R"aumen vermittelnden 
Funktion die Rede ist. Bei \cite{reich} findet sich aber auch nicht der Begriff des Attraktors.}
Ob sich hierbei geltend machen l"asst, dass Reichls mittlerweile seit bald 
drei Dekaden erschienener \glqq Modern Course in Statistical Physics\grqq 
in die Jahre gekommen sei, darf durchaus bezweifelt werden, wenn
wir einen Blick auf j"ungere Literatur der Thermodynamik werfen.
Gewiss, das volle Ausmass fruchtbarer Anwendung beispielsweise der 
abstrakten Ergodentheorie auf die Thermodynamik gilt es, 
vorzuf"uhren.\newline
Wir fragen uns daher,
ob die Bervorzugung der Diskussion nat"urlicher Topologien
nicht schlicht auch daran liegt,
dass der empirische Wissenschaftler, wenn er nicht zugleich Mathematiker ist,
das Instrumentarium benutzt, dessen Gebrauch ihm vertraut ist, was 
sich dadurch auf die Mathematik "ubertr"agt, welcher der empirische Wissenschaftler mitteilt, was f"ur ihn relevant ist. Ob er dabei oft
einfach
nicht davon Kunde hat, dass es angemessener w"are, andere Mittel einzusetzen -- w"ahrend dem 
Mathematiker gleichzeitig manche Anwendungsm"oglichkeit seiner Errungenschaften entgeht?
Wenn diese Frage zu bejahen ist,
gibt es zwei Komponenten, die 
die Bevorzugung nat"urlicher Topologien bewirken.
Dann kommt zu der direkten Komponente der Elementarit"at des leichten Einstieges in die nat"urlichen Topologien
die sekund"are Komponente der elementarmathematischen Verfasstheit der Modelle empirischer Wissenschaften hinzu,
wo sich doch auch die abstrakte Ergodentheorie recht elementar formulieren l"asst. 
\newline
\newpage
{\Large {\bf Symbolverzeichnis}}\newline\newline
Neben den v"ollig etablierten Bezeichnungen benutzen wir in dieser Abhandlung
auch weniger gel"aufige oder spezielle, f"ur dieses Traktat eigent"umliche Schreibweisen. 
Diese weniger gel"aufigen Schreibweisen f"uhren wir
im fortlaufenden
Text ein und stellen sie hier zusammen.\newline
Bei der Erkl"arung der jeweiligen idiomatischen Notation kommt es 
naturgem"ass "ofter vor, dass diese Erkl"arung eine spezielle Begriffsbildung 
dieses Traktates benutzt. Diese spezielle Begriffsbildung findet sich aber im Index vermerkt,
sodass sich im Text die Stelle auffinden l"asst, wo die jeweilige 
spezielle Begriffsbildung formuliert wird.
\newline
Es sei $n,m\in\mathbb{N}$, $r\in \mathbb{R}^{+}$, $t\in \mathbb{R}$, $(X,d)$ ein metrischer Raum, $\mathbf{T}$
eine Topologie, $\mathbb{T}$ eine Menge von Topologien, 
${\rm E}=\{e\}$
eine einelementige Menge, $j\in\{1,2,\dots n\}$,
$\Theta=(\theta_{1},\theta_{2},\dots \theta_{n})$
ein $n$-Tupel, 
$\Lambda=(\lambda_{1},\lambda_{2},\dots \lambda_{m})$ ein $m$-Tupel,
$\Psi$ eine Flussfunktion, $\psi$ eine reelle Flussfunktion, $\xi$ ein Autobolismus, $\Xi$ 
eine Menge von Autobolismen,
$\phi$ eine Funktion  und es seien ${\rm A}$ und ${\rm B}$ zwei Mengen, $\mathcal{A}$, ${\rm X}_{\star}$, ${\rm X}^{\star}$ 
Mengensysteme, wobei
$$\mathcal{A}^{c}:=\Bigl\{Y\setminus A:A\in\mathcal{A}\Bigr\}$$
das zu $\mathcal{A}$ komplement"are Mengensystem und $Y=\bigcup\mathcal{A}$ sei.
Es ist dann
\newline
\begin{tabbing}
Wir benutzen zudem folgende \= Etwas \kill
${\rm A}^{\mathbb{N}}$ \> die Menge aller Folgen $\{a_{j}\}_{j\in\mathbb{N}}$,\\ 
\quad \> deren Glieder Elemente der Menge ${\rm A}$ sind,\\
$\mathbf{ab}({\rm A})$ \> die Menge aller auf ${\rm A}$ definierten Autobolismen,\\  
$\mbox{{\bf @}}(\psi)$\> die Menge der Attraktoren von $\psi$, \\
$\mbox{{\bf @}}(\Phi,\mathcal{A})$\> die Menge freier Attraktoren bzgl. $\mathcal{A}$, \\
$\underline{\mbox{{\bf @}}}(\mathbb{T},{\rm X}^{\star},{\rm X}_{\star})$ \> die Menge topologischer Attraktoren\\
\quad \> des Tripels $(\mathbb{T},{\rm X}^{\star},{\rm X}_{\star})$,\\
$\mathcal{A}^{\cup}$ \> das Mengensystem aller Vereinigungen\\ \quad \> "uber Teilmengen von $\mathcal{A}$,\\
$\mathcal{A}_{{\rm A}}$ \> die ${\rm A}$-Auswahl aus $\mathcal{A}$, d.h\\
\quad \> $\{a\in \mathcal{A}:a\cap{\rm A}\not=\emptyset\}\ ,$\\
${\rm A}\Delta{\rm B}$\>  die symmetrische Differenz von ${\rm A}$ und ${\rm B}$, \\
$\mathbb{B}_{r}(x)$ \> die offene Kugel des $(X,d)$, in dem $x$ ist, deren\\  
\quad \>  Mittelpunkt $x$ ist und die den Radius $r$ hat,\\
$\mathcal{C}_{+}(\mathcal{A})$\> die Menge bzgl. $\mathcal{A}$ Cantor-stetiger Funktionen,\\
$\mathcal{C}_{-}(\mathcal{A})$\> die Menge bzgl. $\mathcal{A}$ konvers Cantor-stetiger \\
\quad \>  Funktionen,\\
$\Psi^{t}$ \> der durch $\Psi$ und $t$ festgelegte Phasenfluss,\\   
\quad \> d.h die Bijektion $\Psi(\mbox{id},t)$\\ 
\quad \>  des Zustandsraumes auf sich,\\ 
$[\Psi]$  \> die durch $\Psi$ festgelegte trajektorielle Partition,\\ 
\quad \> d.h $[\Psi]:=\{\Psi(x,\mathbb{R}):x\in\mathbf{P}_{2}\Psi\}\ ,$\\ 
$[[\Psi]]_{\mathcal{A}}$ \> die durch $\Psi$ festgelegte Partition \\ 
\quad \>  $[[\Psi]]_{\mathcal{A}}:=\{\mathbf{cl}_{\mathcal{A}}(\tau):\tau\in[\Psi]\}\ ,$\\ 
\quad \> wobei\\
$\mathbf{cl}_{\mathcal{A}}$ \> der H"ullenoperator bzgl. $\mathcal{A}$ ist,\\ 
$:=\bigcap\{A\in\mathcal{A}:A\supset\mbox{id}\}$\> \quad \\
$\mathbf{cl}_{\mathbf{T}}$ \>und\\ 
$:=\bigcap\{A\in\mathcal{A}:A\supset\mbox{id}\}$\>  der H"ullenoperator bzgl. $\mathbf{T}$ ist. Es sei \\
$\mathbf{E}(\mathbf{T})$\> das Mengensystem integerer Mengen\\
\quad \> der Topologie $\mathbf{T}$, \\
$\downarrow {\rm E} =e$ \> das Element der einelementige Menge ${\rm E}$, \\
$\mathbf{or}_{\xi}(z)$ \> der Orbit des Autobolismus $\xi$ durch $z\in\mathbf{P}_{1}\xi$ \\
$\mbox{mon}^{+}(\mathbb{R})$ \> Menge aller 
streng monoton wachsender Folgen,  \\
$\mbox{mon}^{-}(\mathbb{R})$ \> Menge aller 
streng monoton fallender Folgen\\ 
\quad \> des Zahlenstrahles,  \\
$\mathbf{part}({\rm A})$  \> die Menge aller Partitionen der Menge ${\rm A}$,\\
$\mathbf{P}_{j}\Lambda=\lambda_{j}$ , \> \quad\\
$\mathbf{P}_{1}\phi$ \> die Definitionsmenge von $\phi$,\\
$\mathbf{P}_{2}\phi$ \>  die Wertemenge von $\phi$,\\
$\oplus$ \> die Kokatenation,\\ \quad \> sodass $\Theta\oplus\Lambda=(\theta_{1},\theta_{2},\dots \theta_{n},
\lambda_{1},\lambda_{2},\dots \lambda_{m})$ ist,\\
$\mathbf{T}(n)$ \quad \> die nat"urliche Topologie des $\mathbb{R}^{n}$,\\
$\widehat{\mathbf{T}}(\Xi)$ \> die gemeinsame invariante Topologie \\
\quad \> der Autobolismenmenge $\Xi$, \\
$\mathbf{un}(\mathcal{A})$ \quad \> die Menge $\mathcal{A}$-haltiger Teilmengen.\\
$:=\{X\in 2^{Y}:A\in\mathcal{A}\setminus\{\emptyset\}\Rightarrow A\not\subset X\}$\> \quad\\
\end{tabbing}

\renewcommand{\indexname}{Begriffsregister}
\printindex

\end{document}